\documentclass[a4paper,11pt]{article}

\usepackage[a4paper, left=2cm, right=2cm, top=3cm, bottom=3cm]{geometry}

\usepackage[utf8]{inputenc}
\usepackage{amsmath, amssymb, amstext, amsfonts, amsthm, bbm, array, enumerate,scalerel}
\usepackage{mathrsfs, dsfont,mathabx}

\usepackage{graphicx}
\usepackage{caption}

\usepackage{nicefrac}
\usepackage[none]{hyphenat}
\usepackage[dvipsnames]{xcolor}
\usepackage{tikz}
\usepackage{stmaryrd}
\usepackage{tikz-3dplot}
\usepackage{float}
\usepackage{appendix}
\usepackage{longtable}
\usepackage[round]{natbib}
\usepackage{shuffle}
\usepackage{ stmaryrd }
\usepackage{framed}

\usepackage{graphicx}
\usepackage{caption}

\usepackage{tocloft}

\usetikzlibrary{patterns}
\usetikzlibrary{plotmarks}

\newtheoremstyle{boldremark}
    {\dimexpr\topsep/2\relax} 
    {\dimexpr\topsep/2\relax} 
    {}          
    {}          
    {\bfseries} 
    {.}         
    {.5em}      
    {}          

\RequirePackage[colorlinks,citecolor=blue, urlcolor=blue, linkcolor=blue ]{hyperref}

\theoremstyle{plain}
\newtheorem{theorem}{Theorem}[section]
\newtheorem{lemma}[theorem]{Lemma~}
\newtheorem{corollary}[theorem]{Corollary}
\newtheorem{proposition}[theorem]{Proposition}

\theoremstyle{definition}
\newtheorem{definition}[theorem]{Definition}
\newtheorem{example}[theorem]{Example}

\newtheorem*{assumption*}{\assumptionnumber}
\providecommand{\assumptionnumber}{}
\makeatletter

\theoremstyle{boldremark}
\newtheorem{remark}[theorem]{Remark}
\newtheorem{notation}[theorem]{Notation}

\newtheorem{assumpt}[theorem]{Assumption}

\newtheorem{property}[theorem]{Property}

\numberwithin{equation}{section}

\newcommand{\bxi}{\boldsymbol{\xi}}

\def\R{\mathbb{R}}
\def\Rset{\mathbb{R}}
\def\Nset{\mathbb{N}}
\def\N{\mathbb{N}}

\newcommand{\cadlag}{c\`adl\`ag~}
\newcommand{\caglad}{c\`agl\`ad~}
\newcommand{\Cadlag}{C\`adl\`ag~}
\newcommand{\Ito}{It\^{o}}

\newcommand{\Follmer}{Föllmer~}

\newcommand{\tXX}{\widetilde{\XX}}

\newcommand{\e}{\epsilon}

\newcommand{\Mcal}{\mathcal M}
\newcommand{\Mcalp}{\mathcal{M}_{[p]}}

\def\X{\mathbb{X}}
\def\wX{\widehat{\mathbb{X}}}

\def\XX{\mathbf{X}}
\def\wXX{\widehat{\mathbf{X}}}

\def\Y{\mathbb{Y}}
\def\wY{\widehat{\mathbb{Y}}}

\def\YY{\mathbf{Y}}
\def\wYY{\widehat{\mathbf{Y}}}

\def\ZZ{\mathbf{Z}}

\def\dcc{d_{CC}}

\newcommand{\x}{\mathbf x}
\newcommand{\y}{\mathbf y}

\renewcommand{\u}{\mathbf u}
\renewcommand{\v}{\mathbf v}

\newcommand{\fdot}{{\,\cdot\,}}

\usepackage{caption}

\usepackage{scalerel,stackengine}
\stackMath
\newcommand\reallywidehat[1]{%
\savestack{\tmpbox}{\stretchto{%
  \scaleto{%
    \scalerel*[\widthof{\ensuremath{#1}}]{\kern.1pt\mathchar"0362\kern.1pt}%
    {\rule{0ex}{\textheight}}
  }{\textheight}%
}{2.4ex}}%
\stackon[-6.9pt]{#1}{\tmpbox}%
}
\usepackage{scalerel,stackengine}

\newcommand\reallywidetilde[1]{%
  \ThisStyle{%
    \setbox0=\hbox{$\SavedStyle#1$}%
    \stackengine{.8\LMpt}{\SavedStyle#1}{%
      \stretchto{\scaleto{\SavedStyle\mkern.2mu\sim}{.5467\wd0}}{.5\ht0}%
    }{O}{c}{F}{T}{S}%
  }%
}

\title{Functional \Ito-formula and Taylor expansions\\ for non-anticipative maps of \cadlag rough paths}
\author{Christa Cuchiero \thanks{University of Vienna, Department of Statistics and Operations Research, Data Science @ Uni Vienna, Kolingasse 14-16, 1090 Wien, Austria, christa.cuchiero@univie.ac.at
} 
\and Xin Guo  \thanks{University of California, Berkeley, Department of Industrial Engineering and Operations Research, Etcheverry Hall 4141, 94720, Berkeley, CA, USA, xinguo@berkeley.edu}
\and Francesca Primavera  \thanks{University of California, Berkeley, Department of Industrial Engineering and Operations Research, Etcheverry Hall 4141, 94720, Berkeley, CA, USA, francesca.primavera@berkeley.edu
\newline
This research was funded in whole or in part by the Austrian Science Fund (FWF) Y 1235.\newline
We thank Bruno Dupire, Josef Teichmann and Valentin Tissot-Daguette for many fruitful discussions that improved the paper.}
\date{}}

\begin{document}

\setlength{\belowcaptionskip}{1pt}
\maketitle

\begin{abstract}
We derive a functional \Ito-formula for non-anticipative maps of rough paths, based  on the approximation properties of the signature of càdlàg rough paths. This result is a functional extension of the \Ito-formula for càdlàg rough paths (by Friz and Zhang (2018)),  which coincides with the change of variable formula formulated by Dupire (2009) whenever the functionals' representations, the notions of regularity, and the integration concepts can be matched. Unlike  these previous works, we treat the vertical (jump) pertubation via the Marcus transformation,
which allows for  incorporating path functionals where the second order vertical derivatives do not commute, as is the case for typical signature functionals. As a byproduct, we show that sufficiently regular non-anticipative maps admit a functional Taylor expansion in terms of the path's signature, leading to an important generalization of the recent results by Dupire and Tissot-Daguette (2022).
\end{abstract}


\tableofcontents

\section{Introduction}

 The \Ito-formula is one of the major tools of stochastic calculus, extending the chain rule from classical calculus to stochastic analysis and asserting that any sufficiently smooth function of a semimartingale is itself a semimartingale. 
Recognizing 
that the functional dependence on a stochastic system does not only occur through its current state but also through its entire history, 
the seminal papers by \cite{D:09, DUP:09} and \cite{CRF:10,CRF:10b} extend the \Ito-formula to non-anticipative path functionals, which may depend on the past trajectory of a $\Rset^d$-valued semimartingale up to the present time. These works have led to many subsequent contributions, including  \cite{CRF:13},
\cite{LO:11}, ~\cite{K:16}, \cite{CRA:17},~\cite{CRP:18},~\cite{VJ:19}, \cite{HV:24}, and the references therein.

The original formulation of the functional \Ito-formula incorporates vertical and horizontal derivatives, reflecting the functional's dependence on the trajectory and current time, respectively (see Section 2 in \cite{D:09}). This formula is derived using a Taylor expansion of the functional along an approximation of the path, leading to a second-order approximation, which is quadratic in the path increments and thus involves the semimartingale's quadratic variation in the limit.  
It is therefore not surprising that this technique may not only be applied to functionals of semimartingales but also to more general 
paths of finite quadratic variation. This is in fact the finding of the  paper \cite{CRF:13},  which establishes a pathwise functional \Ito-formula for non-anticipative functionals of \cadlag paths of finite quadratic variation in the sense of  \cite{F:81}.
Other forms of generalizations can also be found in \cite{O:26}, which builds on \cite{Bic81}.

Even though these results cover a wide range of pratically relevant path functionals, certain important classes
 of interest are still excluded. Indeed, their framework is limited to functionals that depend continuously on the trajectories of the path with respect to (some variants of) the uniform topology (see e.g., Section 1.2 in \cite{CRF:13}). This excludes several crucial examples, such as  the \Ito-map, which describes the correspondence
between the solution of a stochastic differential equation and the driving signal, or standard \Ito/Stratonovich integrals. To see this, consider, for instance, the sequence of real-valued paths $X^{1,n},X^{2,n}:[0,2\pi]\rightarrow\Rset$ defined via
\begin{align*}
    X^{1,n}_t:=-n^{-\frac{1}{3}}\cos(nt ), \qquad X^{2,n}_t:=n^{-\frac{1}{3}}\sin(nt), \qquad \text{for all }t\in [0,2\pi],
\end{align*}
for $n\in \N$. Then, 
\begin{align}\label{eq: controesempio uniform topology}
    \int_0^{2\pi} X^{2,n}_sdX^{1,n}_s=n^{\frac{1}{3}}\pi,
\end{align}
implying that $X^{1,n},X^{2,n}\rightarrow 0$ in the uniform topology, but $ \int_0^{2\pi} X^{2,n}_sdX^{1,n}_s\rightarrow \infty$, (see \cite{AA:21}), showing that iterated integrals are not continuous 
with respect to the uniform topology.
These topological limitations raise the question of defining a suitable topology on the space of paths to enable a more general applicability of a (pathwise) functional Itô-formula. 

Our approach is inspired  by the theory of rough paths, pioneered by Terry Lyons (see \cite{L:98, L:07}), where one of the fundamental results 
consists of identifying a natural family of topologies on path space so that iterated Stratonovich and \Ito-integrals, and more generally solutions to controlled differential equations, are continuous maps with respect to the driving signal. Relying on this idea, 
we thus equip the path spaces with a (stronger) variation topology (Definition~\ref{eq: q variation distance}) and consider functionals of \emph{weakly geometric~\cadlag $p$-rough paths} for $p\in[1,3)$ (Definition~\ref{3def: roughpathsGROUP}). This means that when dealing with paths of finite $p$-variation  for $p\geq2$, which is the case for sample paths of a semimartingale, we consider functionals that depend on a lifted path, i.e., the $\Rset^d$-valued path itself  and (some of) its rough path lift. In this setup
 we establish  a pathwise \textit{functional \Ito-formula for \cadlag rough paths} that covers the above examples. 
With this approach we can  also treat functionals that just depend on the path itself (and not necessarily on its lift) as in the settings of \cite{D:09} and \cite{CRF:10,CRF:10b}.

Our proof technique is based on a \textit{density approach building on linear functions of the signature} (Section~\ref{3sec: linear functions as Marcus}), which are in fact specific non-anticipative functionals of weakly geometric rough paths that exhibit powerful approximating properties (see e.g., Section 3 in \cite{CPS:22} as well as \cite{KLP:20}, \cite{BHRS:21}, \cite{CSTuat:23}). Specifically, we first establish an \Ito-formula for these specific functionals and then extend it to more general maps via a density argument. This  proof technique has also been used in classical stochastic calculus (see, e.g., Theorem 5.7 in \cite{MS:17}), as well as for measure-valued processes (see, e.g., \cite{GHX:22}, \cite{SR:03_measure}), relying on  density results for certain classes of cylindrical functions.

The key component to make this approach work in the current setup is an appropriate
form of a Nachbin-type universal approximation theorem (UAT) (see \cite{N:49}) for functionals of weakly geometric \cadlag $p$-rough paths for $p\in[1,3)$ and certain derivatives thereof  (Theorem~\ref{3th: Nachbin for path functionals}).
In this respect there are two crucial concepts that need to be developed, namely,  \textit{non-anticipative Marcus canonical} path functionals, and their \emph{vertical derivatives.}

First, inspired by  \cite{CF:19} (see also \cite{M:78,M:81} and \cite{Williams2001}) 
and in particular by the construction of the signature of \cadlag paths (Section~\ref{sec: signature}), we introduce the class 
of \textit{non-anticipative Marcus canonical} path functionals. Roughly speaking, these are maps
\begin{align*}
    F:[0,1]\times D([0,1],G^{[p]}(\Rset^d))\rightarrow \Rset,
\end{align*}
which depend on weakly geometric \cadlag rough paths (denoted by $D([0,1],G^{[p]}(\Rset^d)$) in a non anticipative way 
(Definition~\ref{def: non anticipative path functional}) and which are invariant under the Marcus transformation, i.e.,
\begin{align*}
    F(t,\XX)=F(\psi_{R}^{-1}(\tau_{\XX,R}(t)),\tXX),
\end{align*}
for $\tXX$ denoting the Marcus transformed path of $\XX$ with respect to some pair $(R, \psi_R)$. More precisely, $\tXX$ denotes the continuous path obtained by interpolating the states before and after each jump time of $\XX$ via the log-linear path-function and satisfying   $ \widetilde \XX_{\psi_R^{-1}(\tau_{\XX,R}(\cdot))}=\XX_\cdot $ (see Section~\ref{3sec: Marcus transformation} and Definition~\ref{3def: Marcus canonical path functional} for the precise definition of the involved quantities, and also ~\cite{FK:85}, \cite{AT:92}, \cite{KP:95}, and \cite{CFKM:20}
where a similar approach is used). 

Second, viewing  weakly geometric rough paths as (free step-$[p]$ nilpotent)
Lie group valued paths, we introduce a notion of vertical 
derivatives for Marcus canonical  paths functionals, inspired by the (Euclidean) one introduced in \cite{D:09}. 
We call the quantity 
\begin{align}\label{intro3: vertical derivative}
   U^iF(t,\XX):=\frac{d}{dh}F(\psi_{R}^{-1}(\tau_{\XX,R}(t)),\widetilde{\XX}\otimes \exp^{([p])}(h\e_i)1_{\{\cdot \geq \psi_{R}^{-1}(\tau_{\XX,R}(t)\}})|_{h=0},
\end{align}
vertical derivative of $F$ in the direction $i=1,\dots,d$ at $(t,\XX)$, whenever it exists 
(Definition~\ref{def: vertical differentiable path functional}). 
Notice that~\eqref{intro3: vertical derivative} might be interpreted as a directional derivative of the functional 
in the direction determined by the vertical perturbation $\exp^{([p])}(h\e_i)1_{\{\cdot \geq \psi_{R}^{-1}(\tau_{\XX,R}(t)\}}$, which is designed  to stay in the Lie group and where $\e_i$ denotes the $i$-th canonical basis vector of $\mathbb{R}^d$.
This  is conceptually consistent with~\cite{QT:11}, where a first attempt 
for studying a differential structure of the (non-linear) space of rough paths has been made (see also \cite{S:22} for a discussion on this topic).

A crucial aspect {when inserting a \cadlag path in~\eqref{intro3: vertical derivative} is the following: the Marcus transformation 
of the path $\XX$ needs to be computed before the vertical perturbation. On  one hand, 
this  preserves the Marcus property of the functional also on the level of the (functional) derivative 
(Proposition~\ref{3prop: F Marcus canonical UiF Marcus canonical}). On the other hand, if the original path $\XX$ admits a jump at time $t$, 
the Marcus property of the functional $F$ allows to interpret the derivative in~\eqref{intro3: vertical derivative} as a derivative involving 
a delayed perturbation of the original path
\begin{align*}
      \XX^t\otimes \exp^{([p])}(h\e_i)1_{\{\cdot\geq t+\varepsilon\}},
\end{align*}
for some $\varepsilon>0$ independent of the specific pair $(R,\psi_R)$ 
(Remark~\ref{3rem: vertical derivative indipendent on the pair}\ref{3rem: vertical derivative indipendent on the pair item ii}). This specification 
is particularly relevant when computing the \textit{higher-order vertical derivatives}, which is through an iterative application of the procedure in~\eqref{intro3: vertical derivative} (Definition~\ref{def: higher order derivatives}). In this case, by definition, the perturbations are always computed at a jump time, 
resulting therefore in a notion of  vertical derivatives of order $k\in \N$, $k\geq2$, which involve delayed perturbations of the form
\begin{align*}
      \XX^t\otimes \exp^{([p])}(h_k\e_{i_{k}})1_{\{\cdot\geq t+\varepsilon_k\}}\otimes\dots\otimes \exp^{([p])}(h_1\e_{i_{1}})1_{\{\cdot\geq t+\varepsilon_k+\dots+\varepsilon_1\}},
\end{align*}
for some $h_k,\dots,h_1\in \R$, $\varepsilon_k,\dots,\varepsilon_1>0$, and $i_k,\dots,i_1=1,\dots,d$. 

Relying on these two notions,  we  then establish the first main result of the paper, a universal approximation theorem (UAT) for vertically differentiable path functionals: 
 any $C^K$-non-anticipative Marcus canonical  path functionals $F$ (Definition~\ref{3def: C^K path functionals}) 
evaluated at some \textit{tracking jumps}-extended path $\wXX$ (Definiton~\ref{3def: tracking-jumps paths})  
can be approximated uniformly in time together with its derivatives by linear functionals of the signature and their derivatives 
(Theorem~\ref{3th: Nachbin for path functionals}). Due to the non-linear structure of the vertical (Lie) derivatives, the proof of this  result is highly  delicate.
Indeed, it is a tricky combination of  Nachbin-type theorems (see \cite{N:49}) and some key concepts from Lie group theory.

 With the above notions of vertical derivatives, a functional \Ito-formula for linear functions of the signature follows  by 
the definition of the signature of weakly geometric~\cadlag rough paths.
Furthermore,  the UAT for functionals of weakly geometric~\cadlag $p$-rough paths, combined with some interpolation arguments, yields our second main result: a (rough) functional \Ito-formula for the class of
$C^{[p]+1}$ non-anticipative Marcus canonical path functional (Theorems~\ref{3th: Ito formula p 12},~\ref{3th: Ito formula p 23}). Here, it is required that the functional itself, its derivatives, 
and  a certain remainder (Definition~\ref{3def: 2 parameter functional}) 
are  continuous with respect to the above mentioned variation norms
(Definition~\ref{3def: var continuous funcitonals}).
We then also show that our It\^o-formula
matches some standard as well as functional \Ito-formulas in the literature (see Section~\ref{sec: connections litarature}).

The last main result of the paper is the \textit{functional Taylor expansion} in terms of the \textit{signature} (Theorems~\ref{th: taylor 12},~\ref{3th: taylor p23 continuous}). 
To the best of our knowledge, this is the first purely deterministic (rough) Taylor expansion of functionals of weakly geometric \cadlag $p$-rough paths. 
It nevertheless shares similarities with the work by \cite{BKZtaylorrough:15}, where a rough Taylor expansion is derived by 
identifying the vertical derivatives with the abstract notion of Gubinelli derivatives, 
with expansions coming from control theory e.g., \cite{FL:81, F:83, F:86}, \cite{BLM:23}, as well as stochastic Taylor expansions 
(see e.g., \cite{LO:11}, \cite{PK:92}, \cite{BGA:Stochasticexp}), and  recent results in  \cite{DT:23}.

Let us reiterate that the treatment via the Marcus transformation, which leads to non commutative higher order derivatives,
is crucial for the Taylor expansion in terms of the signature components, which are non-symmetric tensors due to the non-commutativity of the iterated integrals.
These results would not follow from a direct application of the 
differential calculus introduced in \cite{D:09} and \cite{CRF:10} as the higher order functional derivatives always take values in the space of symmetric tensors over $\Rset^d$ 
(see e.g., \cite{CRP:18}, \cite{B:24} as well as Remarks 
\ref{rem: non commutativity derivate} and~\ref{remark: taylor p23}\ref{rem: valentin dupire}).

\paragraph{Organization of the paper.} In Section~\ref{sec: funcionals of RP}, we introduce the space of non-anticipative Marcus canonical  path functionals and the corresponding differential calculus, as well as the considered $p$-variation topologies}. In Section~\ref{sec: uat differential}, we introduce the set of ``tracking-jumps-extended paths'' and present the UAT for vertically differentiable path functionals. In Section~\ref{3sec: functional Ito formula} and Section~\ref{sec: functional taylor}, we prove the functional \Ito-formula and the Taylor expansion, respectively. In the Appendix,  we collect the technical proofs and some remarks on rough integration theory as well as 
different pathwise integration approaches to which we compare the current rough setting.

\section{Preliminaries}
\subsection{Algebraic setting}\label{3sec: algebra}

Fix $d\in \mathbb N$ and let $\Rset^d$  be the Euclidean space. The extended tensor algebra over $\Rset^d$ is defined by
	\[ T((\Rset^d)):=\{\mathbf{u}=(\mathbf{u}^{(0)},\mathbf{u}^{(1)},\dots,\mathbf{u}^{(n)},\dots) \ | \ \mathbf{u}^{(n)}\in \mathbb{R}^{\otimes n}\},
    \]
	where $(\Rset^d)^{\otimes n}$ denotes the $n$-fold tensor product of $\Rset^d$ with the convention $(\Rset^d)^{\otimes 0} := \Rset$. We equip
	$T((\Rset^d))$ with the standard addition $+$, tensor multiplication $\otimes$, and scalar multiplication.
	For $N\in\Nset$, the truncated tensor algebra is defined by
	\[ T^N(\Rset^d):=\{\mathbf{u}=(\mathbf{u}^{(0)},\mathbf{u}^{(1)},\dots,\mathbf{u}^{(N)}) \ | \ \mathbf{u}^{(n)}\in \mathbb{R}^{\otimes n} \ \text{for }n\leq N\},
    \]
    and the tensor algebra via $T(\Rset^d):=\bigcup_{N\in \Nset} T^N(\Rset^d)$.
 Let $\pi_n:T((\Rset^d))\rightarrow(\Rset^d)^{\otimes n}$ and  $\pi_{\leq N}:T((\Rset^d))\rightarrow T^N(\Rset^d)$ be the maps such that for $\u\in T((\Rset^d)),$ 
\begin{align*}
    \pi_n(\u):=\u^{(n)}, \qquad \pi_{\leq N}(\u):=(\u^{(n)})_{n=0}^N.
\end{align*}
  For $c\in \Rset$, set 
$$T^N_c(\Rset^d):=\{\u\in T^N(\Rset^d)\colon \u^{(0)}=c\}.$$
 The space $T_1^N(\Rset^d)$ is a Lie group under the tensor multiplication $\otimes$, truncated beyond level $N$. The neutral element with respect to $\otimes$ is $\mathbf{1}:=(1,0,\dots,0)\in T_1^N(\Rset^d)$. Moreover, for any $\u=(\mathbf{1}+\mathbf{b})\in T_1^N(\Rset^d)$, with $\mathbf{b}\in T_0^N(\Rset^d)$, its inverse is given by
 \begin{align}\label{3eq:inv}
     \u^{-1}=\sum_{k=0}^N(-1)^k\mathbf{b}^{\otimes k}.
 \end{align}
The exponential and logarithm maps are defined as follows:
\begin{equation}\label{eq: expN}
\begin{aligned}
			\exp^{(N)}:T&_{0}^N(\Rset^d)\rightarrow T_{1}^N(\Rset^d)\qquad\quad&	\log^{(N)}:T&_{1}^N(\Rset^d)\rightarrow T_{0}^N(\Rset^d)\\
		&\mathbf{b}\mapsto\mathbf{1}+\sum_{k=1}^N\frac{\mathbf{b}^{\otimes k}}{k!},\qquad&
		&\mathbf{1}+\mathbf{\mathbf{b}}\mapsto\sum_{k=1}^N(-1)^{k+1}\frac{\mathbf{b}^{\otimes k}}{k},
	\end{aligned}\end{equation}
	where the tensor multiplication is again always truncated beyond level $N$. 
We furthermore introduce the (non-truncated) exponential map, which is given by 
\begin{align}\label{eq: exponential map}
\exp({\u}):=1+\sum_{k=1}^\infty\frac{\u^{\otimes k}}{k!}\in T((\Rset^{d+1})),
\end{align}
for each $ \u\in T((\mathbb R^{d+1}))$ such that  $\pi_0(\u)=0$. Let $\mathfrak{g}^N(\Rset^d)$ be the \emph{free step-$N$ nilpotent Lie algebra} over $\Rset^d$, i.e.,
\begin{align}\label{3eq: Lie algebra N}
\mathfrak{g}^N(\Rset^d):=\{0\} \oplus \Rset^d\oplus[\Rset^d,\Rset^d]\oplus\dots\oplus \underbrace{[\Rset^d,[\Rset^d,\dots ,[\Rset^d,\Rset^d]]]}_{(N-1)\text{ brackets}}\subseteq
T_0^N(\Rset^d),
\end{align}
where, for $\u\in T_0^{M}(\Rset^d)$, $1\leq M\leq N-1$, $\mathbf{b}^{(1)}\in \Rset^d$, $[\mathbf{b}^{(1)},\u]:=\mathbf{b}^{(1)}\otimes \u-\u\otimes \mathbf{b}^{(1)}$. 

The image of $\mathfrak{g}^N(\Rset^{d})$ through the exponential map is a subgroup of $T_1^N(\Rset^d)$ with respect to $\otimes$. It is called \emph{free step-$N$ nilpotent Lie group} and is denoted by
\begin{align}\label{3eq:GN}
    G^N(\Rset^d):=\exp^{(N)}(\mathfrak{g}^N(\Rset^d)).
\end{align}
We equip it with the so-called Carnot-Caratheodory (CC) norm $\|\cdot\|_{CC}$ (see Definition and Theorem 7.32 in~\cite{FV:10}) and the induced (left-invariant) metric, denoted by $d_{CC}$ (see Definition 7.41 in~\cite{FV:10}).
Finally, we introduce the set of so-called \textit{group-like elements}, defined via
\begin{align}\label{3eq: G group inf}
    G((\Rset^d)):=\{\x\in T((\Rset^d)) \ | \ \pi_{\leq N }(\x)\in G^N(\Rset^d) \text{ for all }N \}.
\end{align}
We refer to Chapter 7 in~\cite{FV:10} for more details on these algebraic aspects and the specific group $G^N(\Rset^d)$ (see also Section 2 in~\cite{L:07}), and to~\cite{BLU:07} and~\cite{S:22} for a more general treatment
of Lie groups.

Let $I = (i_1,\dots, i_n)$ be a multi-index with entries in $\{1,\dots,d\}$. Denoting by $\epsilon_1,\ldots,\epsilon_d$ the canonical basis of $\mathbb R^d$, we use the notation $|I|:=n$ and $\epsilon_I:=\epsilon_{i_1}\otimes \epsilon_{i_2}\otimes\dots\otimes  \epsilon_{i_n}$. Observe that $(\epsilon_I)_I$ is the canonical orthonormal basis of  $(\Rset^d)^{\otimes n}$.
	Furthermore, we denote by $\e_\emptyset$ the basis element of $(\R^d)^{\otimes 0}$ and set $|\emptyset|:=0$.  We also set $I':=(i_1,\dots, i_{n-1})$ for $n>1$, ${I'}=\emptyset$ for $n=1$,   $I^{''}:=(I')'$ for $n>1$, and use the convention $\e_{I''}=0$ for $n=1$. Given $\x\in T((\Rset^d)),$ we write $\x_{I}:=\langle\x,\e_I\rangle$ and for each ${\u\in T(\Rset^d)}$, we set
\begin{align*}
    \langle \u, \x\rangle:=\sum_{|I|\geq0}\u_I\x_I\in \Rset.
\end{align*}
For $k\in \Nset$, we denote by $(\fdot)^{(k)}:T((\mathbb{R}^{d}))\to( T((\mathbb{R}^{d})))^{\otimes k}$ the shifts given by
\begin{align}\label{3eq: shifts}
\u^{(k)}_I:=\sum_{|J|\geq0}\u_{JI}\e_J,
\end{align}
for each $|I|=k$, where $JI$ denotes the concatenation of the multi-indices $J$ and $I$. 
We also write  $\u^{(0)}:=\u$ for notational convenience.

For two multi-indices $I\in \{1,\dots,d\}^{|I|}$, $J\in \{1,\dots,d\}^{|J|}$, and $a,b\in\{1,\dots,d\}$, the shuffle product $\shuffle$ is defined recursively by
\begin{align*}
	&I\shuffle \emptyset=\emptyset\shuffle I=I,\\
	&(I,a)\shuffle (J,b)=((I\shuffle (J,b)),a)+(((I,a)\shuffle J),b),
\end{align*}
where $(I, a)$ denotes the concatenation of the multi-index $I$ with the element $a$. 

Via the shuffle product, every polynomial on the set of group-like elements $G((\Rset^d))$ admits a linear representative. More precisely, for $\x\in G((\Rset^d))$ and two multi-indices $I\in \{1,\dots,d\}^{|I|}$, $J\in \{1,\dots,d\}^{|J|}$, it holds that
\begin{align}\label{3eq:shuffle}
    \langle \epsilon_{I},\x\rangle \langle \epsilon_{J},\x\rangle =\langle \epsilon_{I}\shuffle \epsilon_{J},\x\rangle,
\end{align}
where $\epsilon_{I}\shuffle \epsilon_{J}:=\sum_{k=1}^K\epsilon_{I_k}$ with $K, I_k$ for $k=1, \ldots, K$ determined via $I\shuffle J=\sum_{k=1}^K I_k$.

\subsection{Weakly geometric \cadlag rough paths}\label{sec: rough paths and rough integration}

Throughout, we denote by
$C([0, 1],E)$ and $D([0, 1],E)$ the space of continuous and \cadlag maps
(paths), respectively, from the interval $[0, 1]$ into a metric space $(E, d)$ equipped with metric $d$.
For $t\in (0,1]$, we denote a partition of $[0,t]$ by $\pi_{[0,t]}=\{0=t_0<t_1<\dots<t_k=t\}$, and write $\sum_{t_i\in \pi_{[0,t]}}$ for the summation over all points in $\pi_{[0,t]}$. The mesh size of $\pi_{[0,t]}$ is given by $|\pi_{[0,t]}|:=\max\{t_{i+1}-t_i \ : \ i=0,\dots,k-1\}$. For $p > 0$, we define the $p$-variation of a path $X\in D([0, 1],E)$ by
\begin{align}\label{3pvarpaths}
    \|X\|_{p\text{-}var}:=\sup_{\pi_{[0,1]}\subset[0,1]}\left(\sum_{t_i\in\pi_{[0,t]}}d(X_{t_i},X_{t_{i+1}})^p\right)^{\frac{1}{p}}.
\end{align}
If $X$ takes values in a vector space, for $s,t\in [0,1]$, $s\leq t$, we use the shortcut $X_{s,t}:=X_t-X_s$ and denote the jumps by $\Delta X_t:=\lim_{s \nnearrow t} X_{s,t}$.~The space of continuous and \cadlag  paths of finite $p$-variation are denoted respectively by $C^p([0, 1],E)$ and $D^p([0, 1],E)$. We endow these spaces with the $p$-variation pseudometric, defined via
\begin{align}\label{eq: q variation distance}
d_p(X,Y):=\sup_{\pi_{[0,1]}\subset[0,1]}\left(\sum_{t_i\in\pi_{[0,1]}}d(X_{t_i,t_{i+1}},Y_{t_i,t_{i+1}})^p\right)^{\frac{1}{p}},
\end{align}
for all $X,Y\in D^p([0,1],E)$. Additionally, we also consider two-parameter functions $A:\Delta_1\rightarrow V$, where $(V,\|\cdot\|)$ is a normed vector space and $\Delta_1:=\{(s,t)\in [0,1]\times [0,1] \ : \ s\leq t\}$. Analogously to paths, the notion of $p$-variation is valid for such functions and defined as 
\begin{align*}
\|A\|_{p\text{-}var}:=\sup_{\pi_{[0,1]}\subset[0,1]}\left(\sum_{t_i\in\pi_{[0,1]}}\|A_{t_i,t_{i+1}}\|^p\right)^{\frac{1}{p}}.
\end{align*}
Similarly, we set
$d_p(A,A'):=\|A-A'\|_{p\text{-}var}$
for all $A,A':\Delta_1\rightarrow V$ for which $\|A\|_p,\|A'\|_p<\infty$. We write $\|\cdot \|$ to denote the norm on any vector space $V$ that may differ from case to case. 

We say that a path $X\in D([0,1],E)$ is a time-reparametrization of some $Y\in D([0,1],E)$ if $Y=X_\phi$, for some $\phi$ time-reparametrization, i.e., $\phi:[0,1]\rightarrow[0,1]$ is increasing and bijective.

\noindent Let $C([0, 1],G^N(\Rset^d))$ and $D([0, 1],G^N(\Rset^d))$ be the space of continuous and \cadlag maps
(paths), respectively, from the interval $[0, 1]$ into  $(G^N(\Rset^d), d_{CC})$.
For $\mathbf{X}\in D([0,1],G^N(\Rset^d))$, $s,t\in [0,1]$, $s\leq t$, the path increments  are defined  via
\begin{align}\label{eq: increments group valued path}
\mathbf{X}_{s,t}:=\mathbf{X}_{s}^{-1}\otimes \mathbf{X}_t
\end{align}
and the jumps by
$\Delta\mathbf{X}_{t}:= \lim_{s \nnearrow t} \mathbf{X}_{s,t}$.
 For $p > 0$, we denote by $[p]$ its entire part.
\Cadlag paths of finite $p$-variation with values in the specific group $G^{[p]}(\Rset^d)$ are called weakly geometric $p$-rough path. We formalize this notion in the definition below.

We restrict the presentation only to paths of finite $p$\text{-}variation for $p\in [1,3)$, which are the most relevant in the settings of stochastic analysis.

\begin{definition}\label{3def: roughpathsGROUP}
	Let $p\in[1,3)$. A \cadlag path $\mathbf{X}:[0,1]\rightarrow G^{[p]}(\Rset^d)$ is a \emph{weakly geometric~\cadlag $p$-rough path over $\mathbb R^d$} if $\|\mathbf{X}\|_{p\text{-}var}<\infty.$
We denote the space of such paths by $D^{p}([0,1],G^{[p]}(\Rset^{d}))$ and its subspace consisting of continuous paths by $C^{p}([0,1],G^{[p]}(\Rset^{d}))$.  

For $p\in[2,3)$, $\XX\in D^{p}([0,1],G^{[p]}(\Rset^{d}))$ is \emph{Marcus-like} if for all $t\in [0,1]$,
 \begin{align*}
     \log^{(2)}(\Delta\mathbf{X}_t)\in \{0\}\oplus \Rset^d\oplus \{0\}.
 \end{align*}

 \end{definition}
\begin{assumpt}
    Unless otherwise specified, in this paper, we always assume $p\in [1,3)$.
\end{assumpt}

Finally, we introduce the notion of a controlled rough path with respect to $X:=\pi_1(\XX)$, for some $\XX\in D^{p}([0,1],G^{2}(\Rset^{d}))$ with $p\in [2,3)$.  Let $\mathcal{L}(\Rset^d,\Rset^m)$ denote the space of linear maps from $\R^d$ into $\R^m$, for some $m\in \Nset$. 
\begin{definition}\label{3def: controlled path qr}
  Fix $p\in [2,3)$  and $\XX\in D^{p}([0,1],G^{2}(\Rset^{d}))$. Let $p'\geq p$ such that $  \frac{2}{p}+\frac{1}{p'}>1$  and define $r\geq 1$ by the relation $ \frac{1}{r}=\frac{1}{p}+\frac{1}{p'}.$  A pair $(Y,Y')$ is called a \textit{controlled rough path} (with respect to $X:=\pi_1(\XX)$) if $Y\in D^{p}([0,1],\mathcal{L}(\Rset^d,\Rset^m))$,  $ Y'\in D^{p'}([0,1],\mathcal{L}(\Rset^d,\mathcal{L}(\Rset^d,\Rset^m)),$ and $R:\Delta_1\rightarrow\mathcal{L}(\Rset^d,\Rset^m)$, defined by $ R_{s,t}:=Y_{s,t}-Y_s'X_{s,t}$, for $ (s,t)\in \Delta_1,$ has finite $r$-variation. We denote the space of such controlled paths by $\mathcal{V}^{p',r}_{\XX}$.
\end{definition}
\begin{remark}
    Notice that a pair of controlled rough paths $(Y,Y')$ as in Definition~\ref{3def: controlled path qr} is controlled with respect to $X:=\pi_1(\XX)$. Nevertheless, with some abuse of notation, we denote the space of such paths by $\mathcal{V}^{p',r}_{\XX}$.
\end{remark}

 \subsection{Time-stretching of continuous weakly geometric rough paths}

We introduce the notion of a \textit{time-stretched version} of a continuous weakly geometric rough path.

For $p\in [1,3)$, let $\XX\in C([0,1],G^{[p]}(\Rset^d))$ and fix $t\in (0,1]$. Observe that there exist  $N\in \Nset\cup \{\infty\} $ and some sequences 
$(s_k)_{k=1}^N$, $(\tilde{s}_k)_{k=0}^N$, such that 

\begin{align*}
    0=\tilde{s}_0\leq s_1< \tilde{s}_1< \dots< s_N\leq \tilde{s}_N=t,
\end{align*}
and the following conditions hold true. For $k=1,\dots,N$,
\begin{enumerate}
    \item \label{i stretched}if $u\in [\tilde{s}_{k-1},s_k]$, then for every sufficiently small $\varepsilon>0$, there exists $u'\in [\tilde{s}_{k-1},s_k]$ such that $|u-u'|\leq \varepsilon$ and $\XX_u\neq \XX_{u'}$;
    \item \label{ii stretched}if $u\in [s_k,\tilde{s}_{k}]$, then $\XX_u=\XX_{s_k}$.
\end{enumerate}

Notice that these sequences are designed to capture the intervals where the path remains constant.
Moreover, if $N=1$ and $s_1=\tilde{s}_1=t$, then there is no subinterval of $[0,t]$ where the path $\XX$ remains constant. Similarly,  if $s_1=\tilde{s}_0=0$, and $\tilde{s}_1=t$, then $\XX_u=\XX_0$ for all $u\in [0,t]$.

\begin{definition}\label{def: time-stretched version}
  Let ${{\mathbf{X}}}\in {C}^p([0,1],G^{[p]}(\Rset^{d})))$ and fix $t\in (0,1]$. Let 
  $N\in \Nset\cup \{\infty\} $, 
$(s_k)_{k=1}^N$, $(\tilde{s}_k)_{k=0}^N$ be the sequences such that conditions \ref{i stretched},\ref{ii stretched} are satisfied. We say that a  continuous path ${\mathbf{X}^{t,\rhd}}:[0,1]\rightarrow G^{[p]}(\Rset^{d})$ is a \textit{time-stretched version of $\XX$ on $[0,t]$} if   for all $u\in [0,1]$,
      $${\mathbf{X}}^{t,\rhd}_u:={\mathbf{X}_u^{t,\rhd}}1_{\{u<t\}}+{\mathbf{X}_t}1_{\{u\geq t\}},$$
       where,
  \begin{enumerate}
      \item  if $\XX$ is non-constant on $[0,t]$,  for $u\in [0,t)$,
   \begin{equation*}
{\mathbf{X}}^{t,\rhd}_u:=
\begin{cases}
\sum_{k=1}^N\XX_{\theta^k(u)}1_{\{\tilde{s}_{k-1}\leq u\leq \tilde{s}_{k}\}},& \text{if }s_1>0,\\ 
   
   \XX_{\theta^{1,2}(u)}1_{\{  \tilde{s}_0\leq u\leq \tilde{s}_2\}}+\sum_{k=3}^N\XX_{\theta^k(u)}1_{\{\tilde{s}_{k-1}\leq u\leq \tilde{s}_{k}\}}1_{\{N\geq 3\}},& \text{if }s_1=0,
\end{cases}
\end{equation*}
and  for $k=1,\dots,N$, $\theta^k: [\tilde{s}_{k-1},\tilde{s}_{k}]\rightarrow [\tilde{s}_{k-1},{s}_{k}]$,  $\theta^{1,2}:[\tilde{s}_0,\tilde{s}_2]\rightarrow [\tilde{s}_1,{s}_2]$ denotes some increasing continuous bijection;  
\item  if $\XX$ is constant on $[0,t]$,   ${\mathbf{X}}^{t,\rhd}_u:={\mathbf{X}_0}$  for $u\in [0,t)$.
  \end{enumerate}

\end{definition}

\begin{figure}[h]
 \captionsetup{skip=1pt}
 \captionsetup{position=above}
    \centering
    \begin{minipage}{0.45\textwidth}
        \centering
         \caption*{Continuous 1-dimensional path\\
         \vspace{0.1cm}$\XX$}\includegraphics[width=\textwidth]{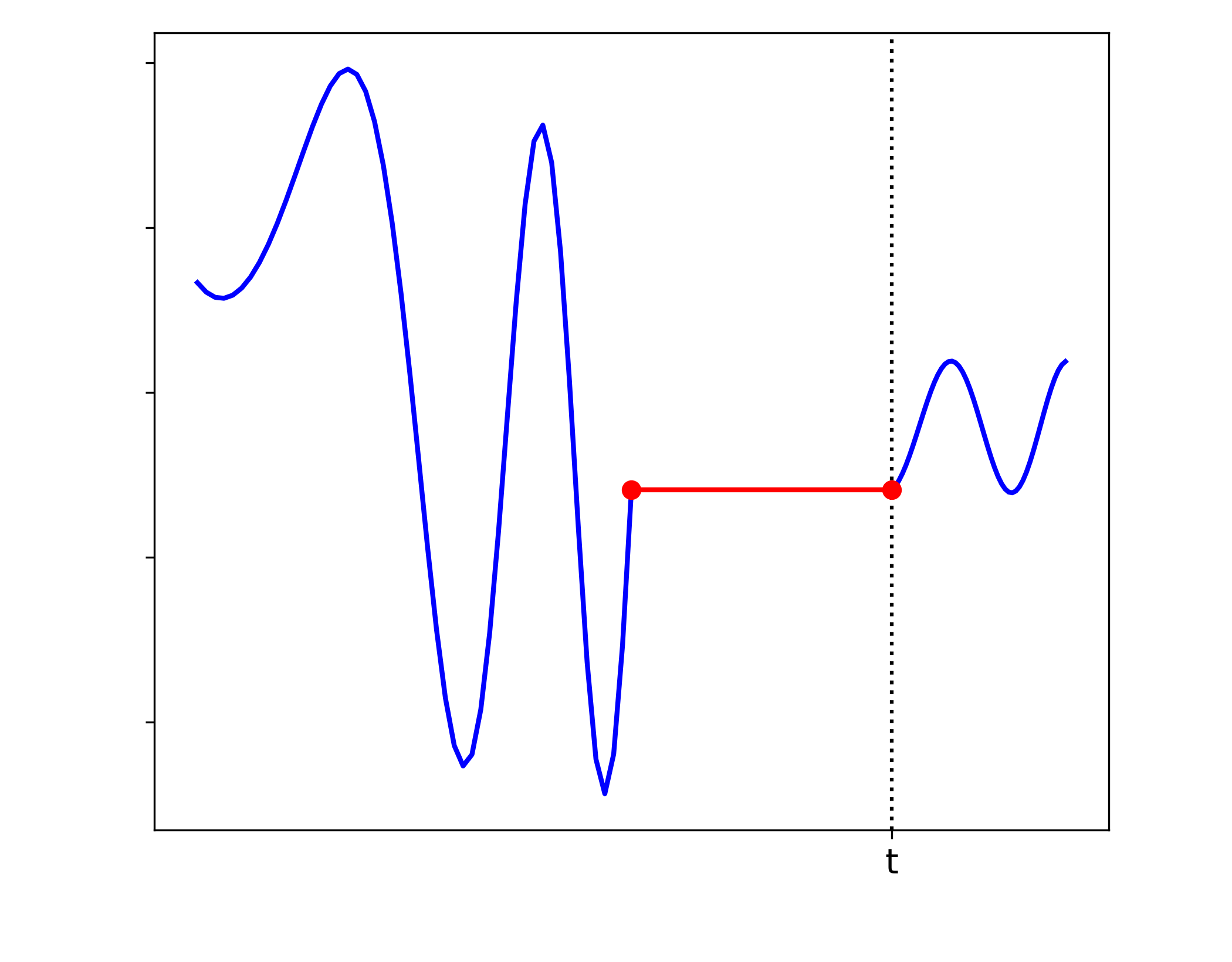}    
    \end{minipage}
   \hspace{-3mm}
    \begin{minipage}{0.45\textwidth}  
        \centering
         \caption*{Stretched path on $[0,t]$\\ \vspace{0.1cm}$\mathbf{X}^{t,\rhd}$}\includegraphics[width=\textwidth]{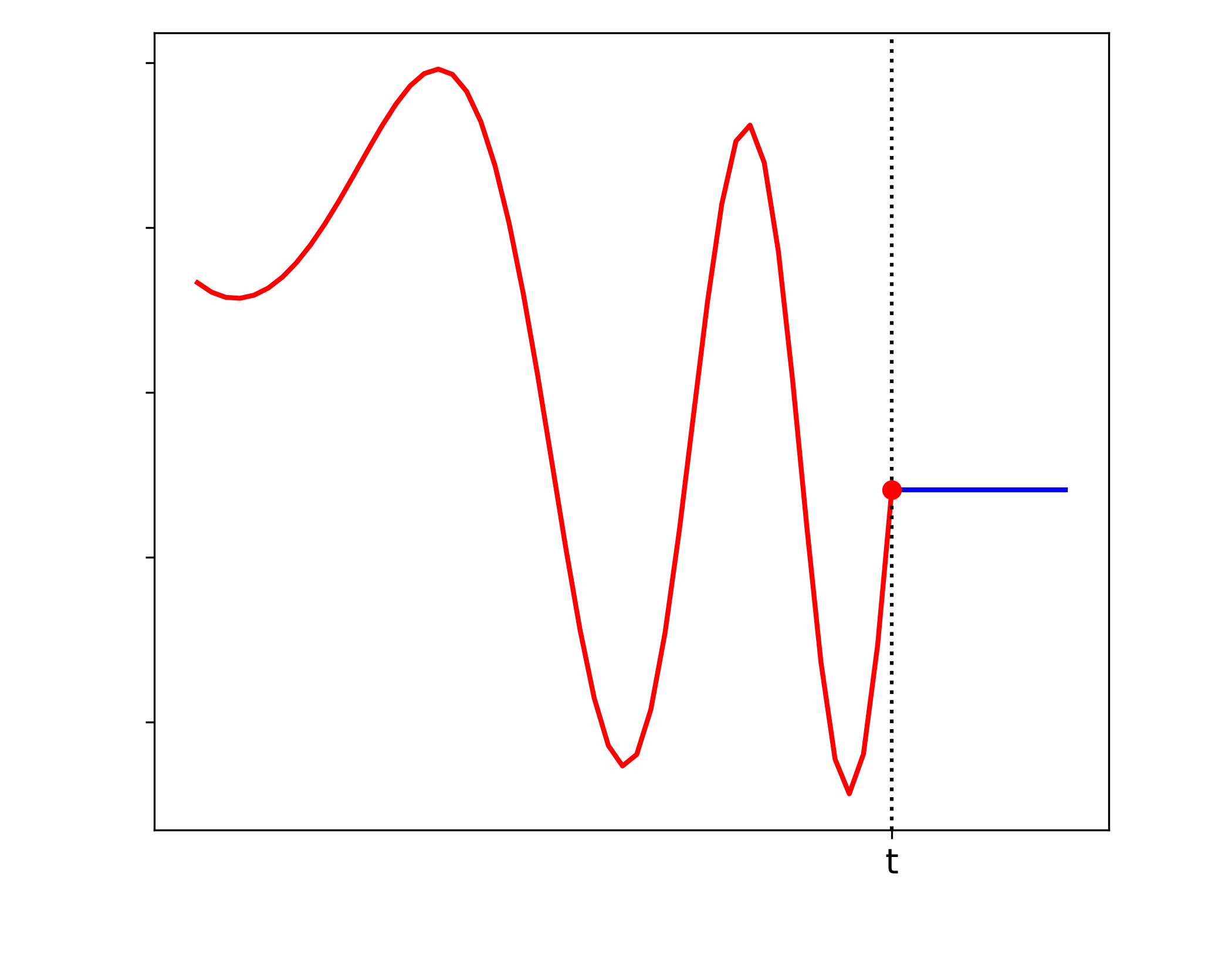} 
    \end{minipage}
\end{figure}

\begin{remark}\label{rem: time-stretched version}
\begin{enumerate}
\item Note that this time-stretching operation removes simply all constant parts of the path (except when $\XX$ is constant on the whole interval $[0,t]$).
    \item \label{rem: time-stretched version i}The definition of the time-stretched version path $\XX^{t,\rhd}$ of a non-constant path on $[0,t]$ depends on the specific bijections $\theta ^k$. Different bijections determine stretched versions that are time reparametrizations of one another, for some reparametrizations that map  $[0,t]$ into   $[0,t]$.
    \item If one component of $\XX$ is strictly increasing, for every $t\in(0,1]$ and $u\in [0,t)$, $\XX^{t,\rhd}_u=\XX_u$.
\end{enumerate}

\end{remark}

\subsection{Marcus transformation of weakly geometric \cadlag rough paths}\label{3sec: Marcus transformation}

In this section, we recall the notion of the Marcus transformation of weakly geometric~\cadlag $p$-rough paths for $p\in [1,3)$, introduced in~\cite{C:17}, (see also Section 2.3 in~\cite{CF:19}). This is an operation that associates every \cadlag path with continuous one obtained by introducing an additional time interval at each jump time and connecting the states before and after the jump via the so-called \textit{log-linear path-function} denoted by $\ell\ell$. Let us start by listing all the objects required for the precise definition.

\begin{enumerate}\phantomsection
    \item Fix $\XX\in D([0,1],G^{[p]}(\Rset^d))$, and let $(t_k)_{k\in\Nset}$ denote the sequence of its jumps times. 
     \item \label{3enumerate: Marcus transformation ii}Fix a sequence $R:=(r_k)_k$ such that for all $k\in \Nset $, $r_k>0$ and $\Sigma_R:=\sum_{k=1}^\infty r_k<\infty$, and define for all $t\in [0,1]$,
\begin{align}\label{3eq:tau X,R}
    \tau_{\XX,R}(t):=t+\sum_{k=1}^\infty r_k1_{\{t_k\leq t\}}.
\end{align}
Notice that $\tau_{\XX,R}$ is an increasing \cadlag function from $[0,1]$ with values in $[0,1+\Sigma_R]$ whose sequence of jumps times coincide with the sequence of jumps times of $\XX$. 
 \item Consider the \textit{log-linear path function} $\ell\ell$ defined as follows:
     \begin{align}\label{eq: ll path funcitonal}
	\ell\ell:G^{[p]}(\Rset^d)\times G^{[p]}(\Rset^d)&\rightarrow C([0,1],G^{[p]}(\Rset^d))\nonumber \\
	(\mathbf{x},\mathbf{y})&\mapsto\Big(s\mapsto\mathbf{x}\otimes \exp^{([p])}(s\log^{({[p]})}(\mathbf{x}^{-1}\otimes\mathbf{y} ))\Big).
	\end{align}
 Notice that for all $(\x,\y)\in G^{[p]}(\Rset^d)\times G^{[p]}(\Rset^d)$, $\ell\ell(\x,\y)(0)=\x$ and $\ell\ell(\x,\y)(1)=\y$. 
 
\item \label{3enumerate: Marcus transformation iv} Define $\YY\in C([0,1+\Sigma_R],G^{[p]}(\Rset^d))$ as follows: for all $s\in [0,1+\Sigma_R]$,
     \begin{equation*}
 \mathbf{Y}_s:=
\begin{cases}
    \XX_t,& \text{if } s=\tau_{\XX,R}(t) \text{ for some }t\in [0,1];\\
   \XX_{t_k^-}\otimes \exp^{([p])}\left(\frac{(s-\tau_{\XX,R}(t_k^-))}{r_k}\log^{([p])}(\Delta \XX_{t_k})\right) & \text{if }s\in [\tau_{\XX,R}(t_k^-),\tau_{\XX,R}(t_k)], \ k\in \Nset. 
\end{cases}
\end{equation*}
\end{enumerate}

\begin{definition}\label{3def: Marcus transformation}
    Let $\XX\in D([0,1],G^{[p]}(\Rset^d))$. Fix a sequence $R$ as in~\ref{3enumerate: Marcus transformation ii} and let $\tau_{\XX,R}$ be the increasing \cadlag function built through $R$ as in equation~\eqref{3eq:tau X,R}. Let $\YY$ be the continuous path defined via $\tau_{\XX,R}$ and the log-linear path function $\ell\ell$ as in~\ref{3enumerate: Marcus transformation iv}, and let $\psi_{R}$ denote an increasing bijection from $[0,1]$ to $[0,1+\Sigma_R]$. The \textit{Marcus-transformed path of $\XX$} associated with the pair $(R,\psi_R)$ is the continuous path $\widetilde \XX\in C([0,1],G^{[p]}(\Rset^d)) $ given by
    \begin{align*}
    \widetilde \XX_{\cdot}:=\YY_{\psi_R(\cdot)}.
\end{align*}

\end{definition}

\begin{figure}[h]
 \captionsetup{skip=1pt}
 \captionsetup{position=above}
    \centering
    \begin{minipage}{0.45\textwidth}
        \centering
         \caption*{\Cadlag $1$-dimensional path\\
         \vspace{0.1cm}$\XX$}\includegraphics[width=\textwidth]{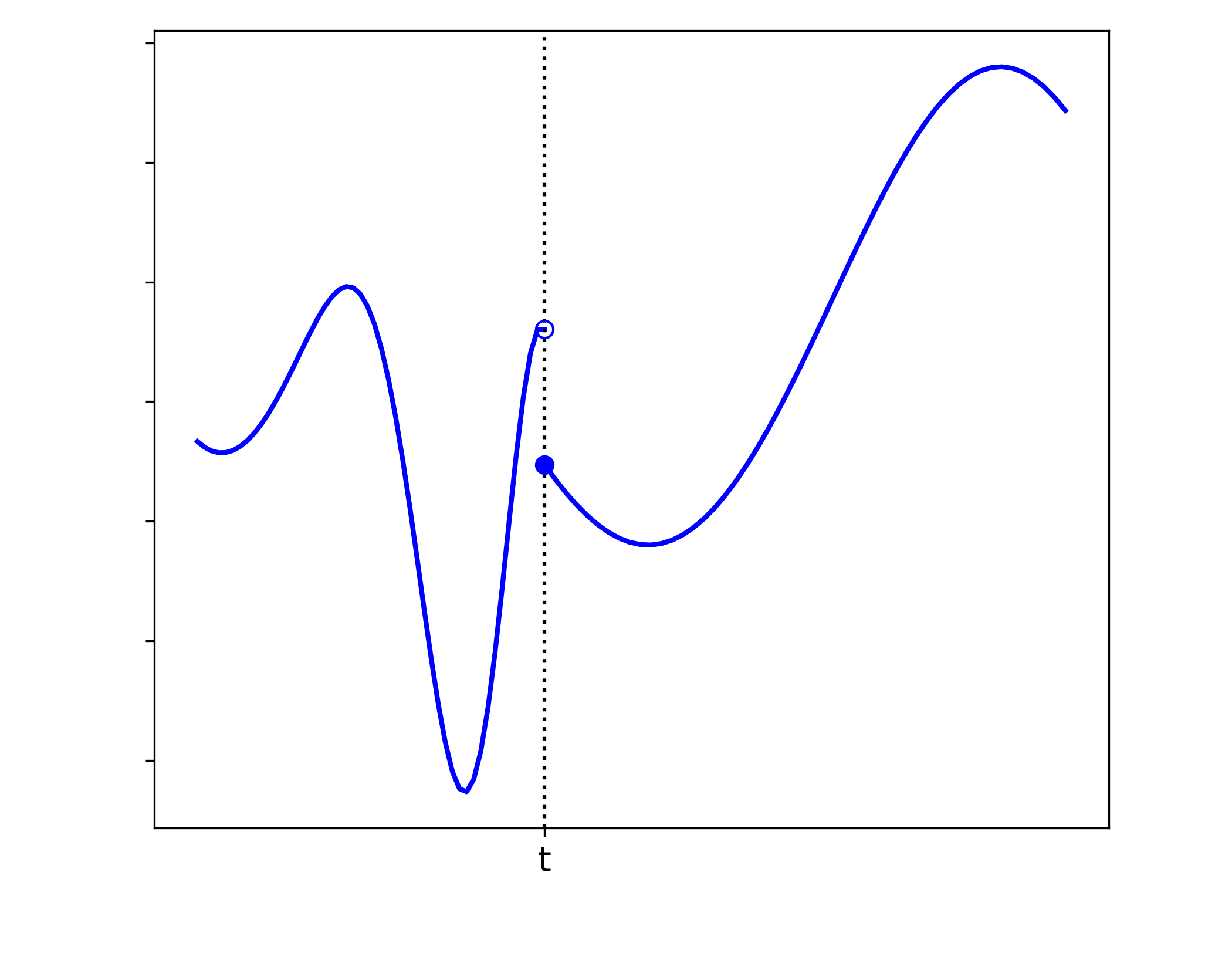}

    \end{minipage}
   \hspace{-3mm}
    \begin{minipage}{0.45\textwidth}  
        \centering
         \caption*{Marcus-transformed path\\ \vspace{0.1cm}$\widetilde{\XX}$}\includegraphics[width=\textwidth]{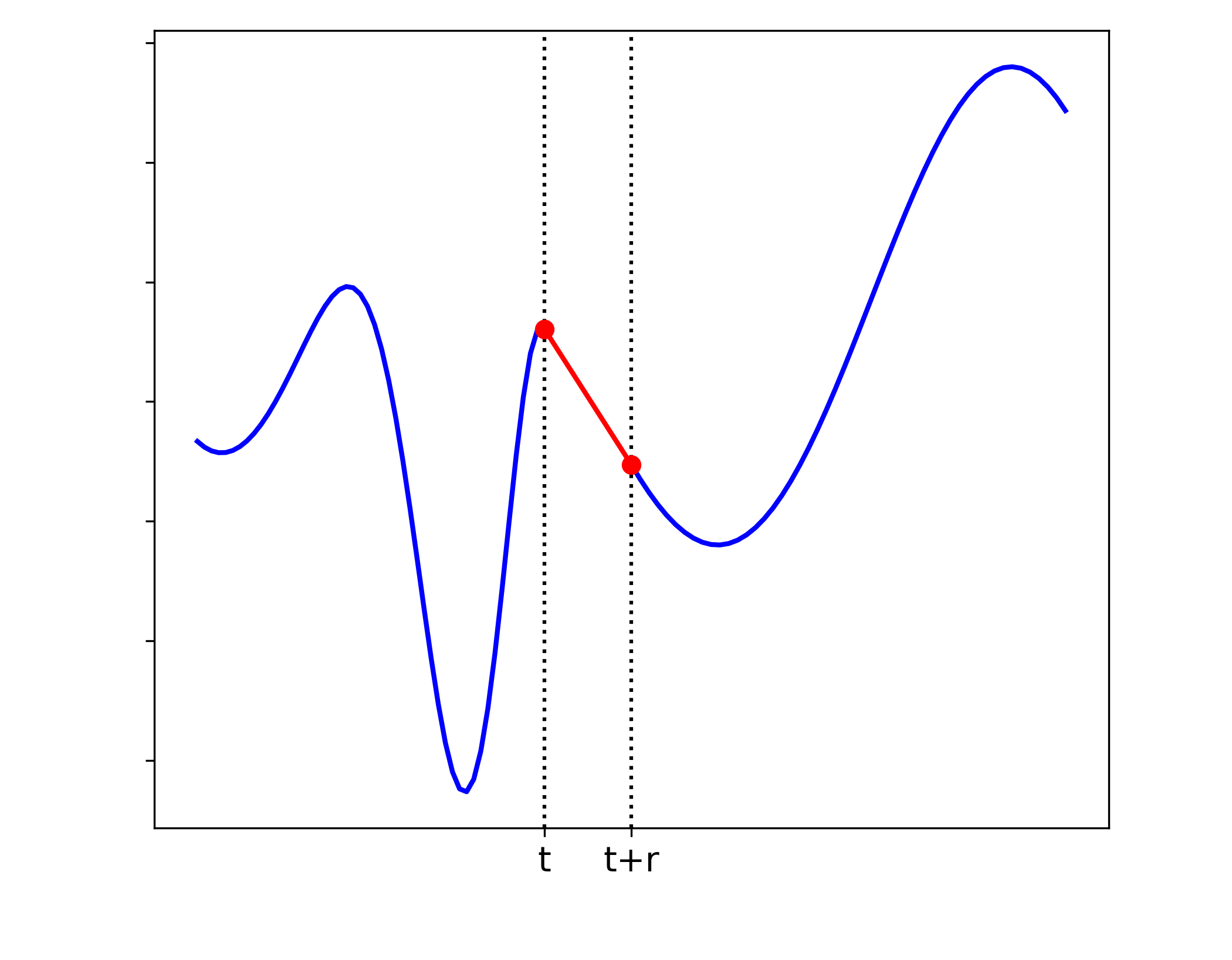}  
    \end{minipage}
\end{figure}

\begin{assumpt}
From now on, whenever we refer to the Marcus-transformed path of a \cadlag path $\XX\in D([0,1],G^{[p]}(\Rset^d))$ associated to a pair $(R,\psi_R)$, we implicitly assume that $R$ is a sequence satisfying condition~\ref{3enumerate: Marcus transformation ii} and $\psi_{R}$ is an increasing bijection from $[0,1]$ to $[0,1+\Sigma_R]$.
\end{assumpt}

\begin{remark}\phantomsection \label{3rem: Marcus transformed path} Fix $\XX\in D([0,1],G^{[p]}(\Rset^d))$.
    \begin{enumerate}
        \item If $\XX\in C([0,1],G^{[p]}(\Rset^d)),$ for any sequence $R$ as in 
       ~\ref{3enumerate: Marcus transformation ii}, the map $\tau_{\XX,R}$ in~\eqref{3eq:tau X,R} is the identity map, i.e., $\tau_{\XX,R}(t)=t$ for all $t\in [0,1]$, and $\YY=\XX$, for $\YY$ as in~\ref{3enumerate: Marcus transformation iv}. In this case, we establish the convention that $\psi_R$ is an increasing bijection from $[0,1]$ to $[0,1]$ and thus $\widetilde \XX$ is simply a time-reparametrization of $\XX$.
          \item Let $\widetilde \XX$ be the Marcus-transformed path of $\XX$ associated to some pair $(R,\psi_R)$. Observe that we can recover $\XX$ from $\widetilde \XX$ via
        \begin{align*}
            \XX_\cdot :=\widetilde \XX_{\psi_R^{-1}(\tau_{\XX,R}(\cdot))}.
        \end{align*}
        \item \label{3rem: Marcus transformed path iii}  The Marcus-transformed paths of $\XX$ associated with two different pairs $(R,\psi_R)$ and $(\widetilde R,\psi_{\widetilde R})$, are simply  time-reparametrizations of one another. More precisely, let $\ZZ$ and $\widetilde{\ZZ}$ denote the transformed paths associated to $(R,\psi_R)$ and $(\widetilde R,\psi_{\widetilde R})$, respectively. Then  ${\ZZ}_\phi=\widetilde{\ZZ}$ for some time-reparametrization $\phi$ such that for all $t\in [0,1]$, $$\phi(\psi_R^{-1}(\tau_{\XX,R}(t)))=\psi^{-1}_{\widetilde{R}}(\tau_{\XX,\widetilde{R}}(t)).$$
      
    \end{enumerate}
\end{remark}

\subsection{Signature of weakly geometric \cadlag rough paths}\label{sec: signature}
In this section, we recall the notion of the signature of weakly geometric \cadlag $p$-rough paths and the key idea behind its construction. More detailed discussions can be found in \cite{FS:17} and Chapter 1 in \cite{P:24}. 
 
The concept of the signature of a \cadlag rough path builds upon the established framework for continuous paths (see e.g.,  \cite{L:98}). More precisely, to compute the signature of a \cadlag rough path $\XX\in D^p([0,1],G^{[p]}(\Rset^d))$ with $\XX_0=\mathbf{1}$, the initial step involves transforming 
$\XX$ into a continuous path via the Marcus transformation (with respect to some pair $(R,\psi_R)$) detailed in Section~\ref{3sec: Marcus transformation}. This transformation results in a path $\widetilde\XX$  which is a continuous weakly geometric $p$-rough path by construction. Lyons's extension theorem guarantees the existence (and uniqueness) of the signature of $\widetilde{\XX}$ (see e.g., Theorem 9.5 in \cite{FV:10}). The signature of the original \cadlag rough path $\XX$ is defined then as the unique $G((\Rset^d))$-valued path $\X$ such that the projection paths over $G^N(\Rset^d)$, denoted by $\X^N$, for $\Nset\ni N>[p]$, are the \cadlag paths given by
 \begin{align}\label{eq: def sig cadlag}
\X^N_\cdot:=\widetilde{\X}^N_{\psi_R^{-1}(\tau_{\XX,R}(\cdot))}.
 \end{align}
Here $\widetilde{\X}^N$ denotes the extension path of $\widetilde{\XX}$ provided by Theorem 9.5 in \cite{FV:10} and $\tau_{\XX,R}$ the \cadlag map defined in equation~\eqref{3eq:tau X,R}. 

These are the key ideas underlying the proof of the analogous theorem in the \cadlag context of Lyons' extension theorem, on which the notion of signature relies.

\begin{theorem}\label{1minimal jump extension}

(\textit{Minimal jump extension theorem}, Theorem 20 in~\cite{FS:17})
   Let $\Nset \ni N>[p]$. Every $\mathbf{X}\in D^p([0,1],G^{[p]}(\Rset^d))$ with $\XX_0=\mathbf{1}\in G^{[p]}(\Rset^d))$ admits a unique extension to a \cadlag path $\mathbb{X}^{N}:[0,1]\rightarrow G^{N}(\Rset^d)$, such that $\mathbb{X}^{N}$ starts from $\mathbf{1}\in G^N(\Rset^d)$, is of finite $p$-variation with respect to 
   $d_{CC}$ on $G^N(\Rset^d)$, and satisfies the following  condition:
\begin{equation}\label{1eqn11}
	\log^{(N)}(\Delta \mathbb{X}_t^N)=\log^{([p])}(\Delta \mathbf{X}_t) \quad  \text{for all }\ t\in[0,1].
	\end{equation}
\end{theorem}

\begin{definition}\label{def: signature cadlag}
    Let  $\mathbf{X}\in D^p([0,1],G^{[p]}(\Rset^d))$ with $\XX_0=\mathbf{1}\in G^{[p]}(\Rset^d))$. The \textit{signature} of $\mathbf{X}$ is defined as the unique path
    \begin{align*}
        \mathbb{X}:[0,1]\rightarrow G((\Rset^d)),
    \end{align*}
    such that for all $\Nset \ni N>[p]$, $\pi_{\leq N}(\mathbb{X})=\mathbb{X}^N$, where $\mathbb{X}^N$ denotes the unique extension path of $\mathbf{X}$ in $G^N(\Rset^d)$ provided by Theorem~\ref{1minimal jump extension}.
\end{definition}
\begin{notation}
    From now on, given $\XX\in D^p([0,1],G^{[p]}(\Rset^d))$, we refer to $\mathbb{X}$ as \textit{signature} of $\mathbf{X}$ and to $\mathbb{X}^N$ as \textit{truncated signature} of order $N$ of $\mathbf{X}$.
\end{notation}

Finally, the truncated signature of a weakly geometric \cadlag rough path can be computed by solving a  (Marcus-type) RDE (see \cite{CF:19}).
\begin{corollary}[Corollary 39 in~\cite{FS:17}]\label{coro: equation signature}
     Let $\mathbf{X}\in D^p([0,1],G^{[p]}(\Rset^d))$ with $\XX_0=\mathbf{1}\in G^{[p]}(\Rset^d))$ and $\Nset \ni N>[p]$. The unique extension path $\mathbb{X}^{N}$ of $\XX$ with values in $G^N(\Rset^d)$ provided by Theorem~\ref{1minimal jump extension} satisfies the  following linear Marcus-type RDE
	\begin{align}\label{eq: MarcusRDE}
&d\mathbb{X}^{N}=\mathbb{X}^{N}\otimes  \diamond d\mathbf{X},\qquad \mathbb{X}^{N}_0=\mathbf{1}\in G^N(\Rset^d),
	\end{align}
 which admits a unique solution, whose explicit form can be written as 
\begin{align}\label{1Rdemarcusexplicit}
	\mathbb{X}^{N}_\cdot=1+\int_0^\cdot\mathbb{X}^{N}_{s^-}\otimes d\mathbf{X}_s+\sum_{0<s\leq \cdot}\mathbb{X}^{N}_{s^-}\otimes \big(\exp^{(N)}(\log^{([p])}(\Delta \mathbf{X}_s))-\Delta \mathbf{X}_s\big).
\end{align}
The integral in~\eqref{1Rdemarcusexplicit} is understood as a  Young (if $p\in [1,2)$) or level 2 rough  (if $p\in [2,3)$) integral  (see Section \ref{sec: rough=young}) and the summation term is well-defined as an absolutely summable series.
\end{corollary}

 Fix $\u \in T(\Rset^d)$ and recall the shifts introduced in equation~\eqref{3eq: shifts}. Let  $\mathbf{X}\in D^p([0,1],G^{[p]}(\Rset^d))$ with $\mathbf{X}_0=\mathbf{1}$ 
and denote by $\X$ its signature.  Then, a projection of equation~\eqref{1Rdemarcusexplicit} along $\u$ combined with Lemma 2.9 in~\cite{FZ:17} yields that 
\begin{enumerate}
      \item if $p\in [1,2)$,
      \begin{align}\label{eq: u,signature p12}
		\langle \u,\mathbb{X}_{\cdot}\rangle=&\langle \u,\mathbb{X}_{0}\rangle+ \int_0^\cdot \langle \u^{(1)},\mathbb{X}_{s^-}\rangle d\XX_{s}\\
  &+\sum_{0<s\leq \cdot}\langle \u,\mathbb{X}_{s}\rangle -\langle \u,\mathbb{X}_{s^-}\rangle -\langle  
 \u^{(1)},\mathbb{X}_{s^-}\rangle \Delta X_{s},  \nonumber 
		\end{align}
 where the integral 
 is a Young integral of $\langle \u^{(1)},\mathbb{X}_{\cdot}\rangle \in D^p([0,1],\Rset^{d+1})$ with respect to $\XX$, and for all $t\in (0,1]$, we set $(1,\Delta X_t):=\Delta \XX_t$. 
  \item if $p\in [2,3)$,
      \begin{align}\label{eq: u,signature p23}
		\langle \u,\mathbb{X}_{\cdot}\rangle =&  \langle \u,\mathbb{X}_{0}\rangle+\int_0^\cdot \langle \u^{(1)},\mathbb{X}_{s^-}\rangle d\XX_s\\
  &+\sum_{0<s\leq \cdot}\langle \u,\mathbb{X}_{s}\rangle -\langle \u,\mathbb{X}_{s^-}\rangle -\langle \u^{(1)},\mathbb{X}_{s^-}\rangle \Delta X_{s}-\langle \u^{(2)},\mathbb{X}_{s^-}\rangle \Delta  \X^{(2)}_{s},\nonumber 
		\end{align}
   where the integral is a  (level 2) rough integral of the controlled rough path 
   $\big(\langle \u^{(1)},\mathbb{X}_{\cdot}\rangle,\langle \u^{(2)},\mathbb{X}_{\cdot}\rangle\big) \in \mathcal{V}^{p,\frac{p}{2}}_{\XX}$  (Definition \ref{3def: controlled path qr}) with respect to $\XX$, and for all $t\in (0,1]$, we set $(1,\Delta X_t, \Delta  \X^{(2)}_{t}):=\Delta \XX_t$. 
  \end{enumerate}

   \begin{assumpt}\label{assumption}
     Throughout, we assume that all the weakly geometric~\cadlag rough paths start at $\mathbf{1}\in G^{[p]}(\Rset^d)$, i.e., $\mathbf{X}_0=\mathbf{1}$. 
 \end{assumpt}

\section{Functionals of weakly geometric~\cadlag rough paths}\label{sec: funcionals of RP}

\subsection{Non-anticipative Marcus canonical  path functionals}

Inspired by the notion of Marcus-type-RDEs (see~\cite{CF:19}),  we here introduce the class of so-called \textit{non-anticipative Marcus canonical path functionals}. To this end, we start with the notion of non-anticipative path functionals, which in turn relies on the definition of \textit{stopped-paths}. In the following, we set $\Nset_0:=\N\cup\{0\}.$

\begin{definition}
  Given ${{\mathbf{X}}}\in {D}^p([0,1],G^{[p]}(\Rset^{d})))$ and $t\in [0,1]$, we define the \textit{stopped weakly geometric} \textit{\cadlag rough path of ${\mathbf{X}}$} stopped at time $t$ as the \cadlag path ${\mathbf{X}}^t:[0,1]\rightarrow G^{[p]}(\Rset^{d})$  given by
$${\mathbf{X}}^t_u:=
         {\mathbf{X}}_u1_{\{u<t\}}+{\mathbf{X}}_t1_{\{u\geq t\}},$$
         for all $u\in [0,1]$.
 \end{definition}
 \begin{figure}[h]
 \captionsetup{skip=1pt}
 \captionsetup{position=above}
    \centering
    \begin{minipage}{0.45\textwidth}
        \centering
         \caption*{\Cadlag $1$-dimensional path\\
         \vspace{0.1cm}$\XX$}\includegraphics[width=\textwidth]{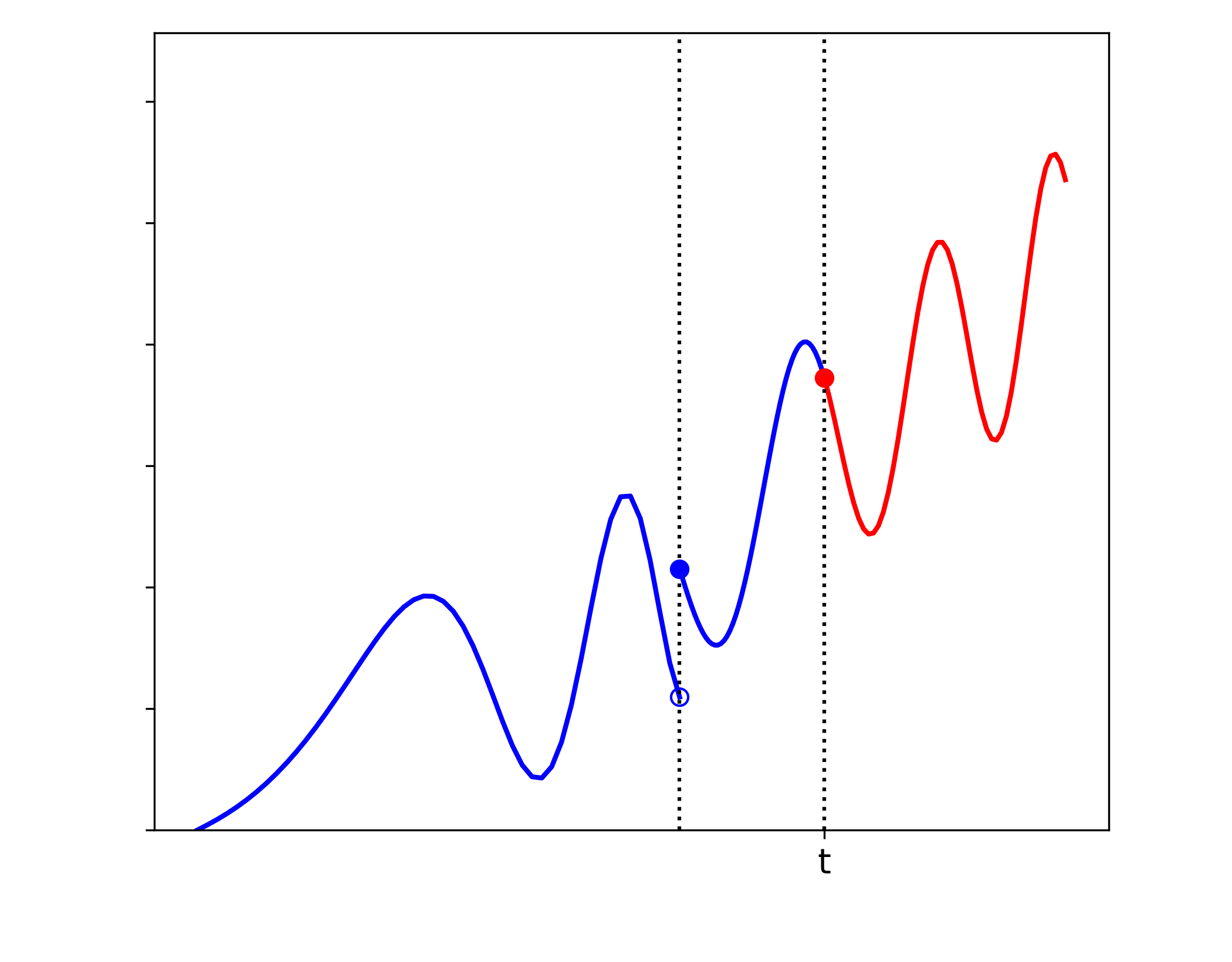}

    \end{minipage}
   \hspace{-3mm}
    \begin{minipage}{0.45\textwidth}  
        \centering
         \caption*{\Cadlag path stopped at time $t$\\
         \vspace{0.1cm}$\XX^t $}\includegraphics[width=\textwidth]{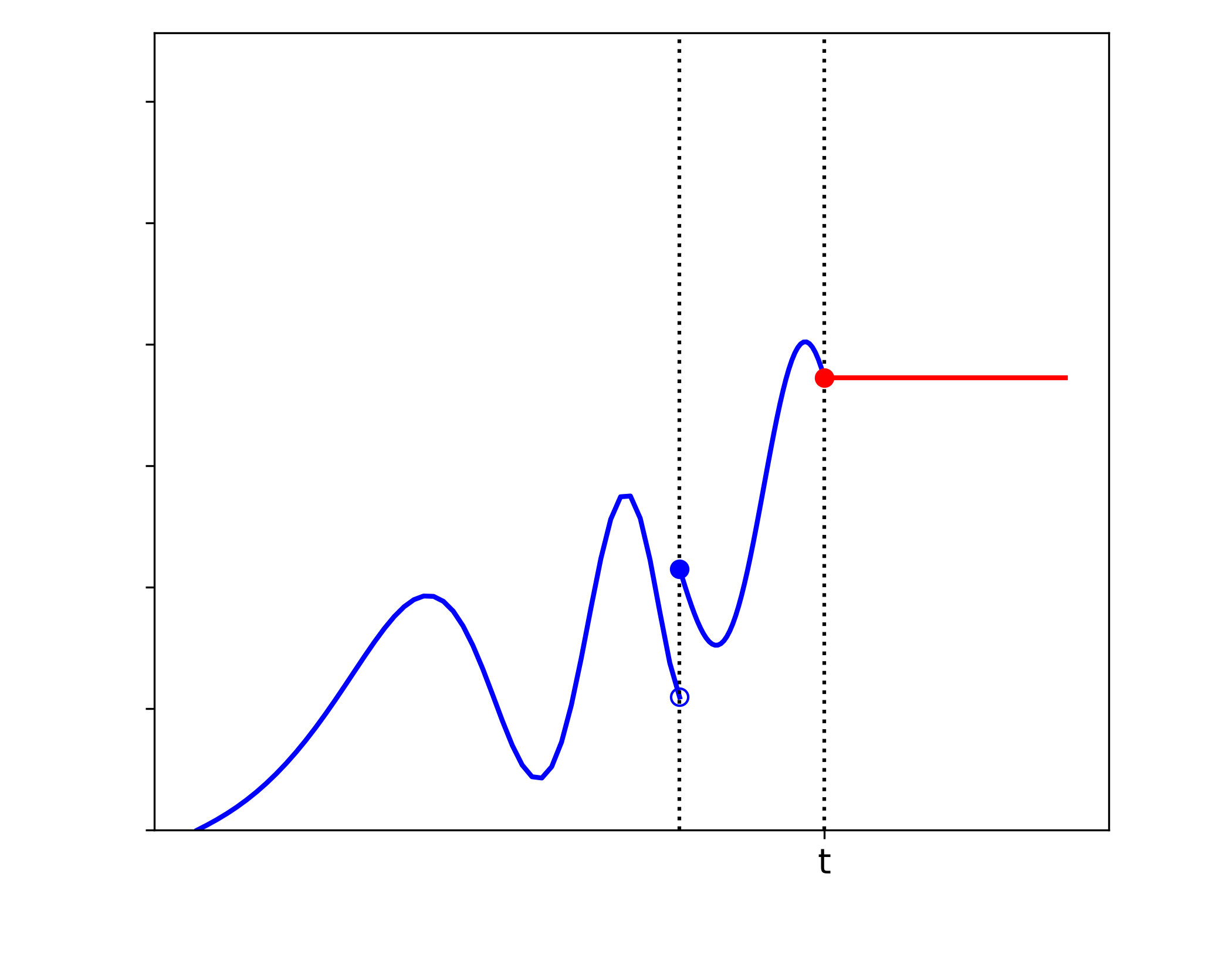} 
    \end{minipage}
\end{figure}

\begin{definition}\label{def: non anticipative path functional}
    Let $F:[0,1]\times D([0,1],G^{[p]}(\Rset^d))\rightarrow \Rset$.  We say that $F$ is a \textit{non-anticipative path functional} if for all $(t,\XX)\in [0,1]\times D([0,1],G^{[p]}(\Rset^d)) $, it holds that 
    \begin{align*}
        F(t,\XX)=F(t,\XX^t).
    \end{align*}
\end{definition}

\begin{remark}\label{rem: t variable} The variable $t$ has the role of a parameter and not of a component of the path $\XX$. In fact, it is the parameter needed to specify that the path functional is non-anticipative.

\end{remark}

Let us now introduce the class of the non-anticipative Marcus canonical  path functionals. We refer to \cite{CF:19} and also ~\cite{FK:85}, \cite{AT:92}, \cite{KP:95} for some related definitions from the literature. 

For $\XX\in D([0,1],G^{[p]}(\Rset^{d}))$, recall the notion of a Marcus-transformed path of $\XX$ associated with some pair $(R,\psi_R)$ given in Definition~\ref{3def: Marcus transformation}  and the one of time-stretched path given in Definition \ref{def: time-stretched version}.
\begin{definition}\label{3def: Marcus canonical path functional}
    Let $F:[0,1]\times D([0,1],G^{[p]}(\Rset^{d}))\rightarrow \Rset$ be a non-anticipative path functional. We say that $F$ is a \textit{non-anticipative Marcus canonical path functional} if
    \begin{enumerate}
        \item \label{item def Mcpf i}for all $\XX\in C([0,1],G^{[p]}(\Rset^d))$, $t\in [0,1]$, $\phi$ time-reparametrization,$$F(t,\XX_\phi)=F(\phi(t),\XX);$$
        \item  \label{item def Mcpf ii}for all $\XX\in C([0,1],G^{[p]}(\Rset^d)) $ and all $t\in [0,1]$,
    $$F(t,\XX)=F(t,\XX^{t,\rhd}),$$
where $\XX^{t,\rhd}$ denotes some time-stretched version of $\XX$ on $[0,t]$;

        \item \label{item def Mcpf iii}for all $(t,\XX)\in[0,1]\times  D([0,1],G^{[p]}(\Rset^d)),$
    \begin{align*}
F(t,\XX)=F(\psi_R^{-1}(\tau_{\XX,R}(t)),\widetilde\XX),\end{align*}
    where $\widetilde \XX$ denotes the Marcus-transformed path of $\XX$ with respect to some pair $(R,\psi_R)$. 
    \end{enumerate}
     We denote the space of such functionals by $\mathcal{M}_{[p]}^0$.
\end{definition}

\begin{notation}
     We  set $((\Mcalp^0)^d)^{\otimes 0}:=\Mcalp^0$ and for $m\in \Nset_0$,  we write $F\in ((\Mcalp^0)^d)^{\otimes m}$ if 
    \begin{align*}
        F:[0,1]\times D([0,1],G^{[p]}(\Rset^{d}))\rightarrow (\Rset^d)^{\otimes m}
    \end{align*}
    is a path functional whose components
    are non-anticipative Marcus canonical path functionals in the sense of Definition~\ref{3def: Marcus canonical path functional}.
\end{notation}

\begin{remark} \label{3rem: Marcus functional independent on the pair}
\begin{enumerate}

\item  The symbol $\Mcalp^0$ does not include the dimension $d$, as it will
always be clear from the context.

 \item \label{item I Marcus canonical does not depend para }Property \ref{item def Mcpf i}  in Definition \ref{3def: Marcus canonical path functional} guarantees that the subsequent conditions \ref{item def Mcpf ii}, \ref{item def Mcpf iii} are independent
 of the specific stretched version of $\XX$ on $[0,t]$ and independent of the specific Marcus-transformed path, respectively. Indeed, let ${\ZZ}$, ${\widetilde\ZZ}$ be two Marcus transformations of $\XX$ associated with  $(R,\psi_R)$ and  $(\widetilde{R},\psi_{\widetilde{R}})$, respectively. By Remark~\ref{3rem: Marcus transformed path}~\ref{3rem: Marcus transformed path iii},  ${\ZZ}_\phi=\widetilde{\ZZ}$, for some  time-reparametrization $\phi$, and for all $t\in [0,1]$, $$\phi(\psi^{-1}_{\widetilde{R}}(\tau_{\XX,\widetilde{R}}(t)))=\psi_R^{-1}(\tau_{\XX,R}(t)).$$ Therefore, by condition \ref{item def Mcpf i},
  \begin{align}\label{3eq: Marcus functional independent on the pair}
      F(\psi^{-1}_{\widetilde{R}}(\tau_{\XX,\widetilde{R}}(t)),\widetilde{\ZZ})=F(\psi^{-1}_{{R}}(\tau_{\XX,{R}}(t)),\ZZ).
  \end{align} 
  Recalling Remark \ref{rem: time-stretched version}\ref{rem: time-stretched version i}, a similar argument holds for condition \ref{item def Mcpf ii}. 

  \item \label{iii2}  A non-anticipative path functional defined only on the set of continuous paths that satisfies conditions~\ref{item def Mcpf i},\ref{item def Mcpf ii} of Definition~\ref{3def: Marcus canonical path functional} can always be extended to a path functional on the entire set of \cadlag path via condition~\ref{item def Mcpf iii}  of Definition~\ref{3def: Marcus canonical path functional}. The resulting functional is a well-defined Marcus canonical path functional. Moreover, this extension is also unique. This is in fact the approach proposed in~\cite{CF:19}.
\end{enumerate}

\end{remark}

\begin{notation}\label{3notation: mut}
      Given the independence on the specific Marcus-transformed path discussed in Remark~\ref{3rem: Marcus functional independent on the pair} and to ease notation, for all $t\in [0,1]$, $\XX\in D([0,1],G^{[p]}(\Rset^d))$ and a pair $(R,\psi_R)$,  we set  $\mu_{t}:=\psi_R^{-1}(\tau_{\XX,R}(t))$ whenever there is no ambiguity.
\end{notation}

 We present some examples of non-anticipative Marcus canonical  path functionals.
\begin{example}\phantomsection\label{example: Marcus canonical}
   \begin{enumerate}
        \item  Let $F(t,\XX):=Y_t$, where $Y$ denotes the solution of a Marcus-type RDE in the sense of Definition 3.1 in \cite{CF:19}, driven by some \cadlag path $\XX\in D^p([0,1],G^{[p]}(\Rset^d))$, at time $t\in (0,1]$. By the solution concept and the property of the rough integral, $F$ is a non-anticipative Marcus canonical  path functional. 

        \item  Let $F$ be the functional such that  for all $(t,\XX)\in [0,T]\times D^p([0,1],G^{[p]}(\Rset^d))$, $$ F(t,\XX):=\sup_{s\in [0,t]}\|\XX_s\|.$$
 Then F is a non-anticipative Marcus canonical  path functional.
 \item \label{item examples Marcus iii}For $\XX\in C^1([0,1],G^1(\Rset^{2}))$, set $\XX^1_\cdot:=\langle  \e_1,\XX_\cdot\rangle=$  and $\XX^2_\cdot:=\langle  \e_2,\XX_\cdot\rangle$, and consider the  path functional (defined only on the set of continuous paths) via $F(t,\XX):=\int_0^t \XX^1_sd\XX^2_s $. Since $F$ verifies~\ref{item def Mcpf i},\ref{item def Mcpf ii} of Definition~\ref{3def: Marcus canonical path functional}, following the discussion in Remark~\ref{3rem: Marcus functional independent on the pair}\ref{iii2}, we can extend it to the set $D^1([0,1],G^1(\Rset^{2}))$ via condition~\ref{item def Mcpf iii}. A direct computation shows that the resulting functional, which is non-anticipative Marcus canonical and we still denote by $F$, explicitly reads as
 \begin{align}\label{eq: F example Marcus}
     F(t,\XX):=\int_0^t \XX^1_{s^-}d\XX^2_s +\frac{1}{2}\sum_{0<s\leq t}\Delta \XX^1_s\Delta \XX^2_s,
 \end{align}
 for all $(t,\XX)\in D^1([0,1],G^1(\Rset^{2}))$.
    \end{enumerate}
\end{example}

\begin{remark}\label{remark: non marucs like}\phantomsection
\begin{enumerate}
    \item  \label{remark: non marucs like i}It is important to note that not every non-anticipative path functional that is well defined on the space of \cadlag paths is a Marcus canonical  path functional. Consider for instance the functional defined via $F(t,\XX):=\int_0^t \XX^1_{s^-}d\XX^2_s $, for all $(t,\XX)\in D^1([0,1],G^1(\Rset^{2}))$. Then, $F$ is not Marcus canonical . Indeed, for  $\XX\in D^1([0,1],G^1(\Rset^{2}))\setminus C^1([0,1],G^1(\Rset^{2}))$, let $\widetilde \XX$ be its Marcus-transformed path with respect to some pair $(R,\psi_R)$. Then, 
    \begin{align*}
        F(\psi_R^{-1}(\tau_{\XX,R}(t)),\widetilde \XX)=\int_0^t \XX^1_{s^-}d\XX^2_s +\frac{1}{2}\sum_{0<s\leq t}\Delta \XX^1_s\Delta \XX^2_s\neq F(t,\XX).
    \end{align*}
    Similarly, the path functional given by $F(t,\XX):=\sum_{0<s\leq t}\Delta \XX^1_s\Delta \XX^2_s$ is not Marcus canonical  as  $0=F(\psi_R^{-1}(\tau_{\XX,R}(t)),\widetilde \XX)\neq F(t,\XX)$ for pure jump paths.

    \item \label{remark: non marucs like ii}
Note however that 
     some non-anticipative path functionals that  do not appear to be Marcus canonical at first glance can be easily turned into Marcus canonical ones. This is the case for instance for the functional $F$ defined via  $F(t,\XX):=\int_0^t \XX_{s}ds $ for all $(t,\XX)\in[0,1]\times  C^1([0,1],G^1(\Rset))$, which does not satisfy condition \ref{item def Mcpf i} in Definition \ref{3def: Marcus canonical path functional}. However, considering the path functional $\bar{F}(t,\XX):=\int_0^t\XX^2_sd\XX_s^1$ defined on $[0,1]\times  C^1([0,1],G^1(\Rset^2))$, we get that it satisfies property \ref{item def Mcpf i} in Definition \ref{3def: Marcus canonical path functional}, and for all $(t,\wXX)\in [0,1]\times C^1([0,1],G^1(\Rset^2))$ where $\wXX:=(Id,\XX)$, with $Id_u:=u$ for all $u\in [0,1]$ and $\XX\in C^1([0,1],G^1(\Rset))$, it holds that 
    $F(t,\XX)=\bar F(t,\wXX)$. Therefore, a time-extension of the original path can be crucial to satisfy the conditions specified in Definition \ref{3def: Marcus canonical path functional}.

    The above example also illustrates that a functional may have multiple representations. To apply the theory that follows, it is necessary to select the representation that satisfies the conditions for being a non-anticipative Marcus canonical path functional.

\end{enumerate}  
\end{remark}

Next, we introduce the notion of path functionals that are invariant under reparametrization, which is 
the same as condition (i) in Definition~\ref{3def: Marcus canonical path functional}, however on the whole space of \cadlag paths $D([0,1],G^{[p]}(\Rset^d))$ and not only on $C([0,1],G^{[p]}(\Rset^d))$.

\begin{definition}\label{3def: invariant under reparametri}
     Let $F:[0,1]\times D([0,1],G^{[p]}(\Rset^d))\rightarrow \Rset$ be a non-anticipative path functional. We say that $F$ is \textit{invariant under reparametrization} if for all $\XX\in D([0,1],G^{[p]}(\Rset^d))$, $t\in [0,1]$ and $\phi$ time-reparametrization,  
\begin{align*}
F(t,\XX_\phi)=F(\phi(t),\XX).
\end{align*}
\end{definition}

We now show that every $F\in \mathcal{M}^0_{[p]}$ satisfies this property. The proof of the following proposition is given in Appendix \ref{proof 3prop: F Marcus canonical then F invariant}.

\begin{proposition}\label{3prop: F Marcus canonical then F invariant}
    Let $F\in \mathcal{M}^0_{[p]}$. Then $F$ is invariant under reparametrization.
\end{proposition}

To conclude, we introduce the notion of the remainder path functional related to some functionals of $G^{2}(\Rset^{d})$-valued paths.

\begin{definition}\label{3def: 2 parameter functional}
    Let $F\in ((\Mcalp^0)^d)^{\otimes m}$  and $F'\in ((\Mcalp^0)^d)^{\otimes m+1}$, for $m\in \Nset_0$.
We define the \textit{remainder path functional} (related to $F$ and $F'$) as follows:
$$R^{F,F'}:\Delta_1\times D([0,1],G^{2}(\Rset^{d}))\rightarrow (\Rset^d)^{\otimes m},$$  
given by 
     \begin{align*}
       R^{F,F'}((s,t),\XX):=F(t,\XX)-F(s,\XX)-F'(s,\XX)\pi_1(\XX_{s,t}),
     \end{align*}
     for all $((s,t),\XX)\in \Delta_1\times D([0,1],G^{2}(\Rset^{d}))$.
\end{definition}

\subsection{Vertically differentiable path functionals}

\begin{definition}\label{def: vertical differentiable path functional}
     Let $F\in \mathcal{M}^0_{[p]}$. We say that $F$ is vertically differentiable at $(t,\XX)\in [0,1]\times D([0,1],G^{[p]}(\Rset^{d}))$ in the direction $i=1,\dots,d$ if the map
\begin{align*}
   \Rset\ni h\mapsto F(\mu_t,\widetilde{\XX}\otimes \exp^{([p])}(h\e_i)1_{\{\cdot \geq \mu_t\}})
\end{align*}
is differentiable at $0$, for some Marcus-transformed path $\widetilde{\XX}$ and $\mu_t\in [0,1]$ given in Notation~\ref{3notation: mut}. In this case we call 
\begin{align}\label{eq: vertical derivativ}
    \frac{d}{dh}F(\mu_t,\widetilde{\XX}\otimes \exp^{([p])}(h\e_i)1_{\{\cdot \geq \mu_t\}})|_{h=0}.
\end{align}
the vertical derivative of $F$ at $(t,\XX)$ in the direction $i$. If $F$ is vertically differentiable in all directions $i=1,\dots,d$ at all $(t,\XX)$, we say that $F$ is vertically differentiable. We denote the space of such functionals by $\mathcal{M}_{[p]}^1$.
\end{definition}

\begin{remark}\phantomsection\label{3rem: vertical derivative indipendent on the pair}
\begin{enumerate}
\item The value of the vertical derivatives of $F\in \Mcalp^1$ at some $(t,\XX)$ is independent of the specific Marcus-transformed path. This follows from Proposition~\ref{3prop: F Marcus canonical then F invariant} and a similar
reasoning as in Remark~\ref{3rem: Marcus functional independent on the pair}.

\item \label{3rem: vertical derivative indipendent on the pair item ii} In~\eqref{eq: vertical derivativ} the Marcus transformation of the path is computed before the vertical perturbation. This is, in fact, a key aspect of this definition. On one hand, it allows the preservation of the Marcus property of the functionals also at the level of the (functional) derivative (see Proposition~\ref{3prop: F Marcus canonical UiF Marcus canonical}). On the other hand, if the original path $\XX$ admits a jump a time $t$, the Marcus property of the functional $F$ allows to interpret the derivative in~\eqref{eq: vertical derivativ} as a derivative computed via a delayed perturbation of the original path:
\begin{align}\label{eq: delayed pert 1}
      \XX^t+h\e_i1_{\{\cdot\geq t+\varepsilon\}},
\end{align}
for some $\varepsilon>0$ independent on the specific pair employed for computing $\tXX$. 

To clarify this aspect, suppose that $\XX$ is a weakly geometric \cadlag $p$-rough path, for $p\in [1,2)$,  that admits a jump only at time $t$. For simplicity, we identify $G^1(\Rset^d)$ with $\Rset^d$. 
Then, applying the Marcus transformation first to the path $\XX$ and then to the perturbed path $\tXX+ h\e_i1_{\{\cdot\geq\mu_t\}}$ yields that~\eqref{eq: vertical derivativ} explicitly reads as
\begin{align*}
    \frac{d}{dh}F(t+\delta_2+\delta_1,\YY^{[i]}(h))|_{h=0},
\end{align*}
for $\YY^{[i]}(h)$ defined  up to time $t+\delta_2+\delta_1$ via 
\begin{equation*}
  \YY^{[i]}(h)_s:=
\begin{cases}
   \XX_{s^-}& \text{if }s\in [0,t],\\
   \XX_{t^-}+\frac{s-t}{\delta_2} \Delta \XX_{t}&  \text{if }s\in [t,t+\delta_2],\\
    \XX_{t^-}+\Delta \XX_t+\frac{s-t-\delta_2}{\delta_1} h\e_{i} &  \text{if }s\in [t+\delta_2, t+\delta_2+\delta_1],
\end{cases}
\end{equation*}
for some $\delta_2,\delta_1>0$ such that $t+\delta_2+\delta_1<1$, $\mu_t=t+\delta_2$.  Since the path $\YY^{[i]}(h)$ stopped at time $t+\delta_2+\delta_1$  is a time-stretched version of some Marcus-transformed path of \eqref{eq: delayed pert 1} on $[0,\mu_{t+\varepsilon}]$, for $\mu_{t+\varepsilon}=t+\delta_2+\delta_1$, by conditions \ref{item def Mcpf ii} and \ref{item def Mcpf iii}, we get
\begin{align*}
   \frac{d}{dh}F(t+\varepsilon,\XX^t+h\e_i1_{\{\cdot \geq t+\varepsilon\}})|_{h=0}= \frac{d}{dh}F(t+\delta_2+\delta_1,\YY^{[i]}(h))|_{h=0}.
\end{align*}

\begin{figure}[h!]
 \captionsetup{skip=1pt}
 \captionsetup{position=above}
    \centering
    \begin{minipage}{0.45\textwidth}
        \centering
         \caption*{Delayed vertical perturbation \\  \vspace{0.1cm}$\XX^t+h\e_11_{\{\cdot\geq t+\varepsilon\}}$}\includegraphics[width=\textwidth]{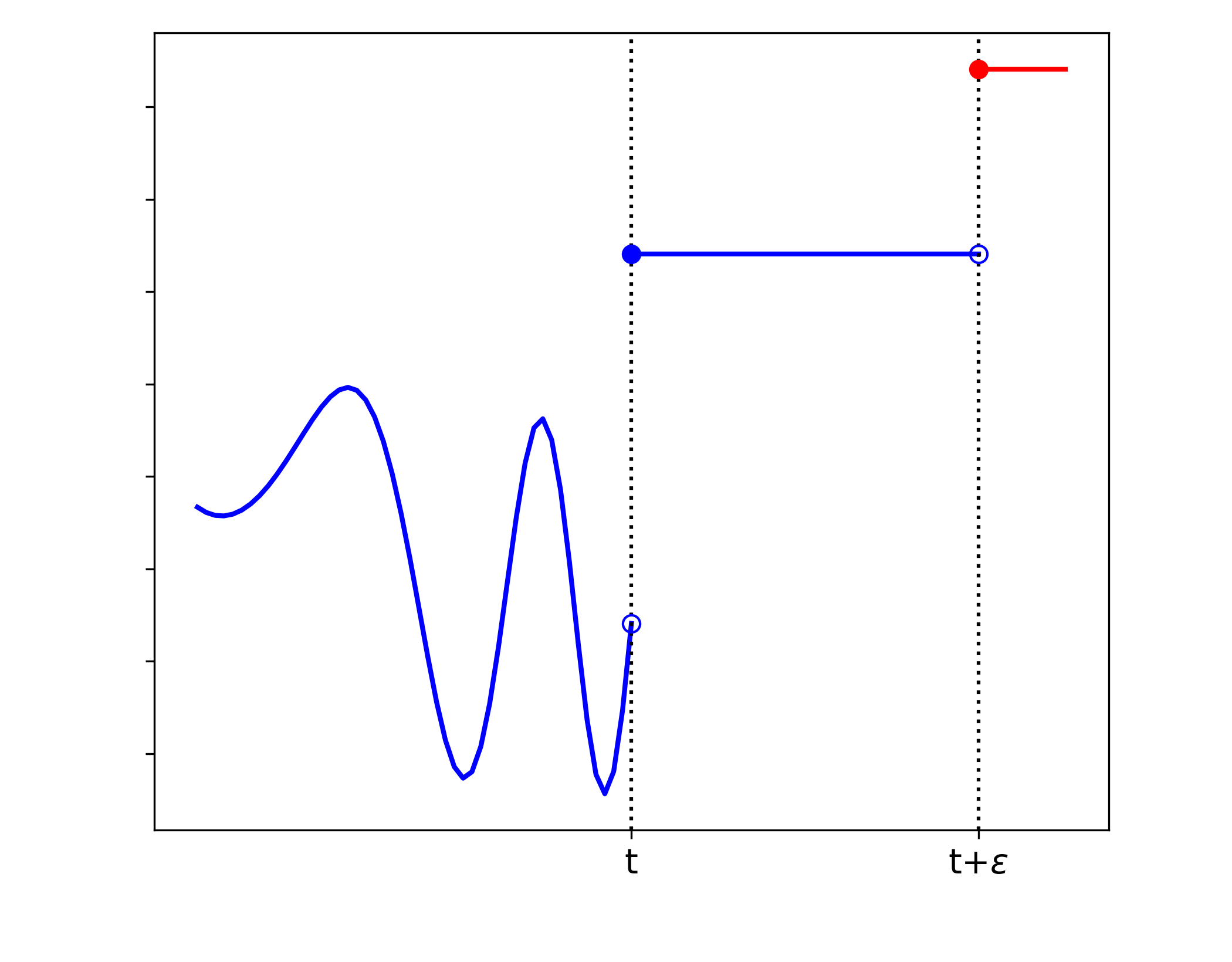}

    \end{minipage}
   \hspace{-3mm}
    \begin{minipage}{0.45\textwidth}  
        \centering
         \caption*{Stretched Marcus tranformation\\ \vspace{0.1cm} $ \YY^{[1]}(h)$}\includegraphics[width=\textwidth]{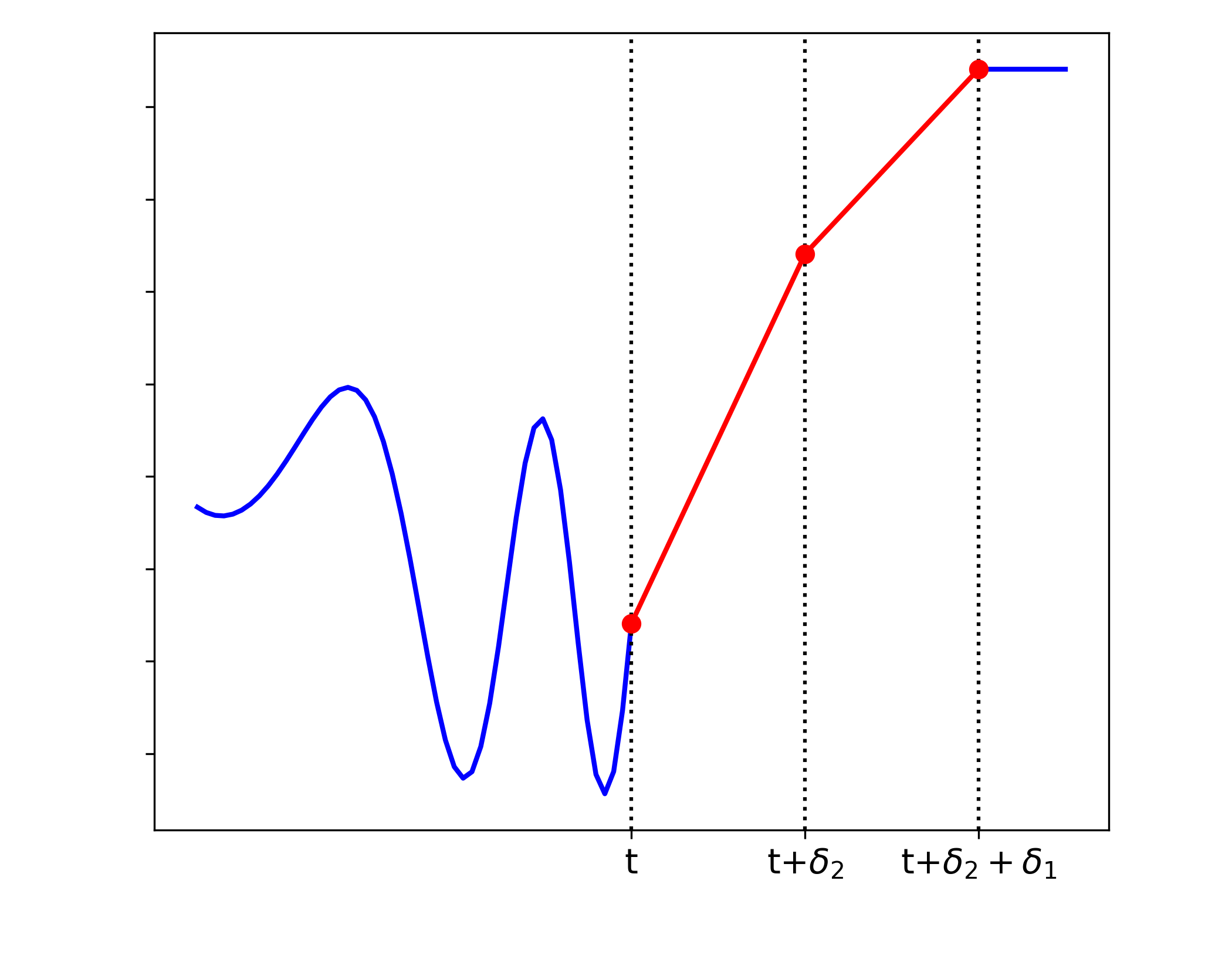}

    \end{minipage}
\end{figure}

This confirms that our notion of vertical derivative involving Marcus transformations of the original path allows to view the derivative in~\eqref{eq: vertical derivativ} as a derivative computed via a delayed perturbation of the original path. This aspect affects in particular the higher order vertical derivatives where a jump at the current time always occurs (see Remark~\ref{rem: non commutativity derivate} and Example~\ref{example non commutative}).  Notice furthermore that the independence on the specific $\varepsilon>0$ follows by the invariance of the functional with respect to time-reparametrization.

   \item   If $p\in [1,2)$, the notion of vertical derivatives given in Definition~\ref{def: vertical differentiable path functional} corresponds to the notion of vertical derivatives introduced in~\cite{D:09} (see also Definition 8 in~\cite{CRF:10}) when evaluated at continuous paths.
   
   Similarly, if for $p\in [2,3)$ a non-anticipative path functional depends only on $\pi_1(\XX)$ for all $(t,\XX)$, then~\eqref{eq: vertical derivativ} matches again the notion of vertical derivative introduced in~\cite{D:09} evaluated at continuous paths.

\end{enumerate}
\end{remark}

Next, we introduce the notion of higher-order vertical derivatives, which are given via an iterative application of the computation in \eqref{eq: vertical derivativ}. To make the argument precise, we formally introduce the differential that associates to every path functional its derivative functional. The well-posedness of this concept relies on the following proposition, whose proof is given in Appendix \ref{proof 3prop: F Marcus canonical UiF Marcus canonical}.

\begin{proposition}\label{3prop: F Marcus canonical UiF Marcus canonical}
    Let $F\in \mathcal{M}^1_{[p]}$. The path functional given by 
    \begin{align}\label{3eq: pre UiF path functional}
  [0,1]&\times D([0,1],G^{[p]}(\Rset^{d}))\rightarrow\Rset \\
    &(t,\XX)\longmapsto \frac{d}{dh}F(\mu_t,\widetilde{\XX}\otimes \exp^{([p])}(h\e_i)1_{\{\cdot \geq \mu_t\}})|_{h=0},\nonumber
\end{align}
is a non-anticipative Marcus canonical  path functional. 
\end{proposition}

\begin{definition}\label{3def: operator Ui}
    For all $i=1,\dots,d$, we define the differential operators
    \begin{align*}
U^i: \ &\mathcal{M}^1_{[p]}\rightarrow\mathcal{M}^0_{[p]}\\
        & \ F \mapsto U^i(F),
    \end{align*}
   where for all $F\in \mathcal{M}^1_{[p]}$ and all $(t,\XX)\in [0,1]\times D([0,1],G^{[p]}(\Rset^{d}))$,
    \begin{align*}
        U^i(F)(t,\XX):= \frac{d}{dh}F(\mu_t,\widetilde{\XX}\otimes \exp^{([p])}(h\e_i)1_{\{\cdot \geq \mu_t\}})|_{h=0},
    \end{align*}
    for some Marcus-transformed path $\widetilde{\XX}$ and $\mu_t\in [0,1]$ given in Notation~\ref{3notation: mut}.

\end{definition}

Finally, we introduce the notion of higher-order vertical derivatives.
\begin{definition}\label{def: higher order derivatives}
    Let $F\in \mathcal{M}^0_{[p]}$ and $K\in \Nset$. We say that $F$ is $K$ times vertically differentiable if for all $l=1,\dots,K$ the path functionals defined recursively by 
    \begin{align*}
    &U^{i_0}F:=F,\qquad \qquad \qquad \qquad \qquad \qquad \quad \qquad \quad \text{for }l=1,\\
    &U^{i_{l-1}}\dots U^{i_1}U^{i_0}F:=U^{i_{l-1}}(U^{i_{l-2}}\dots U^{i_1}U^{i_0}F), \quad \text{for }l=2,\dots,K, \ (i_1,\dots,i_{l-1})\in \{1,\dots,d\}^{l-1},\nonumber 
\end{align*} 
are vertically differentiable at all $(t,\XX)\in  [0,1]\times D([0,1],G^{[p]}(\Rset^{d}))$.  In this case, we call
\begin{align}\label{3eq: value vertical derivatives higher order}
    U^{i_K}\dots U^{i_1}  F(t,\XX):= U^{i_K}(U^{i_{K-1}}\dots U^{i_1}U^{i_0} F)(t,\XX)
\end{align}
the vertical derivative of order $K$ of $F$ at $(t,\XX)$ in the directions $(i_1,\dots,i_K)\in \{1,\dots,d\}^{K}$.
We denote the space of such functionals by $\mathcal{M}_{[p]}^K$.
\end{definition}

\begin{notation}\label{3notation higher order derivatives}
    In the following, for $K\in \Nset$, $F\in \Mcalp^K$, $l=1,\dots,K$, 
   we let $\nabla^l F$ denote the  $(\Rset^d)^{\otimes l}$-valued 
      path functional such that for all $I=(i_1,\dots,i_l)$,
    \begin{align*}
    \nabla ^lF(t,\XX)_I:=U^{i_l}\dots U^{i_1}F(t,\XX),
\end{align*}
for all $(t,\XX)\in [0,1]\times D([0,1],G^{[p]}(\Rset^{d}))$.
Notice that $\nabla^l F$ is a $(\Rset^d)^{\otimes l}$-valued path functional whose components are in $\Mcalp^{K-l}$.  For notational convenience, we also write $\nabla^0F:=F$ and $\nabla F:=\nabla^1 F$.

\end{notation}

The computation of the higher-order vertical derivatives can be done by considering the iterative procedure described below.  For notational simplicity, we here consider only  \cadlag paths with values in $G^{1}(\Rset^{d})$, identify $G^{1}(\Rset^{d})$ as $\Rset^d$, and explicitly write the procedure for computing the derivatives up to the second order.

\noindent Let $F: [0,1]\times D([0,1],G^{1}(\Rset^{d}))\rightarrow\Rset$ and assume $F\in \mathcal{M}^2_1$. Fix $\XX\in D([0,1],G^1(\Rset^d))$, $t\in[0,1]$.

\begin{equation*}
    \begin{array}{l l}
   \text{Step 0:} &\XX^{[0]}:=\XX^t, \\
    &\YY^{[0]}:=\reallywidetilde{\XX^{[0]}}, \\
    &\mu_t^{[0]}:=\psi_{R_0}^{-1}(\tau_{\XX^{[0]},R_0}(t)), \text{ for some pair } (R_0,\psi_{R_0}).
    \end{array}
\end{equation*}

\begin{figure}[H]
 \captionsetup{skip=1pt}
 \captionsetup{position=above}
    \centering
    \caption{Step 0}\label{step0}
    \begin{minipage}{0.45\textwidth}
        \centering
         \caption*{$\XX^{[0]}$}\includegraphics[width=\textwidth]{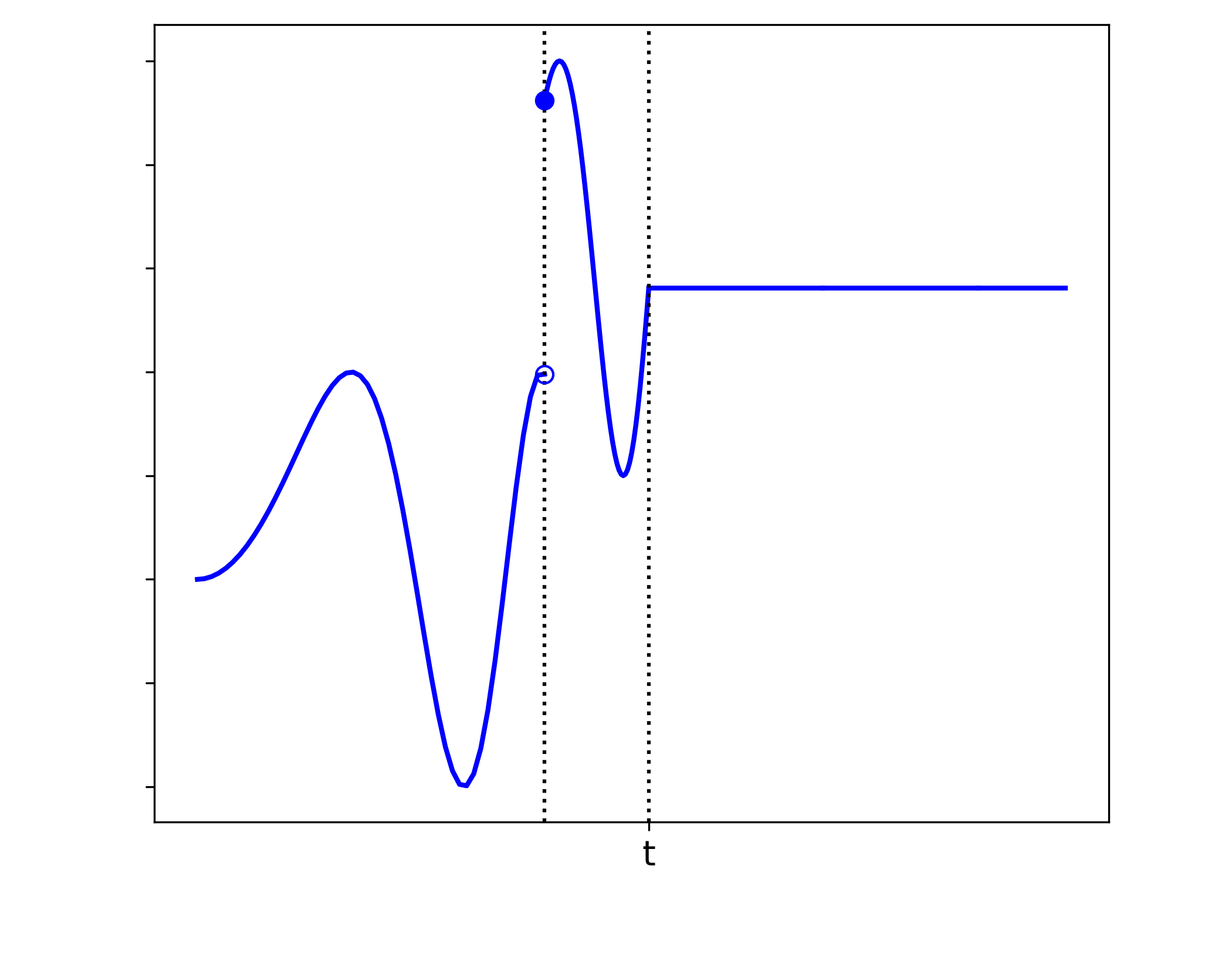}

    \end{minipage}
   \hspace{-3mm}
    \begin{minipage}{0.45\textwidth}  
        \centering
         \caption*{$\YY^{[0]}$}\includegraphics[width=\textwidth]{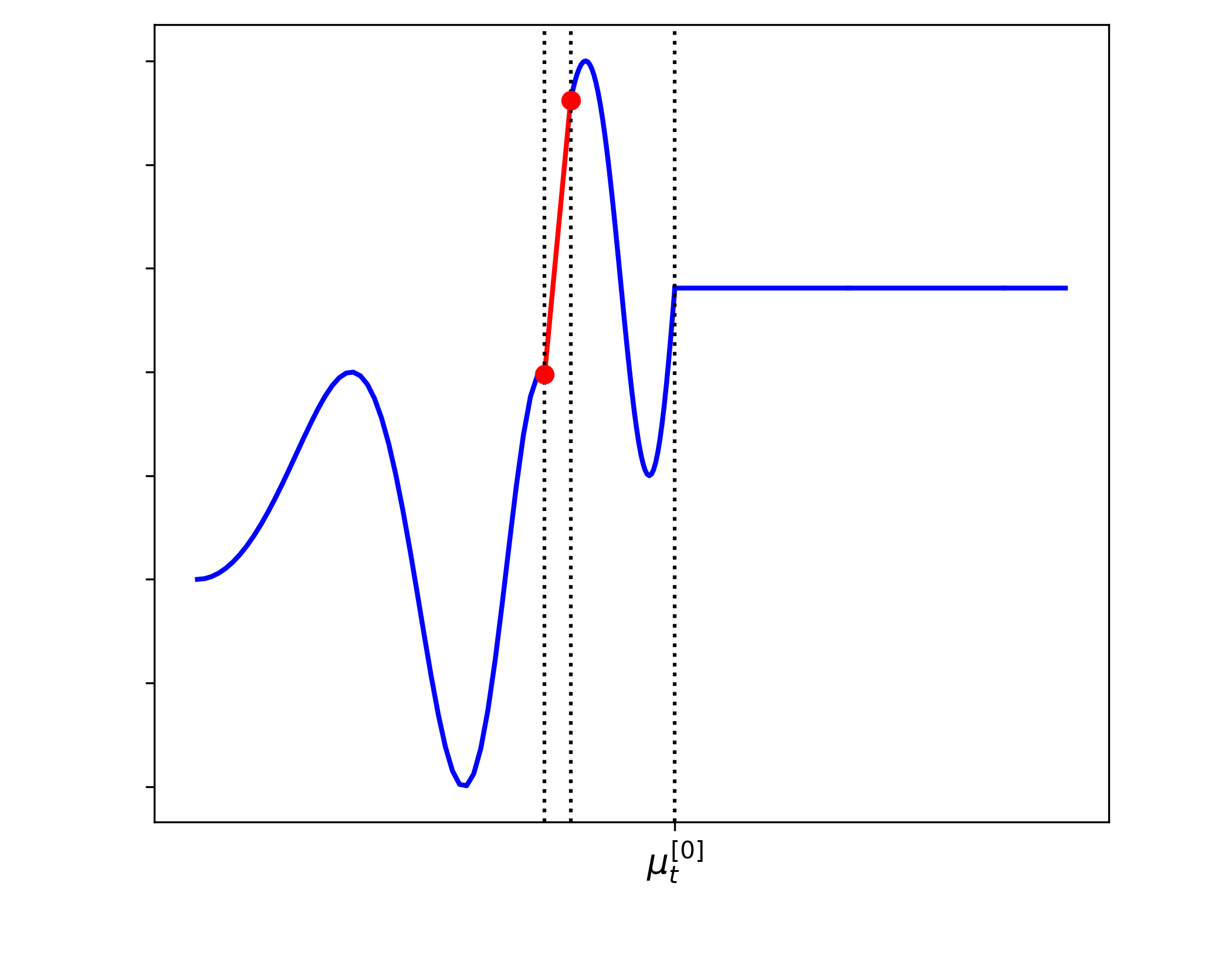} 
    \end{minipage}
\end{figure}

\begin{equation*}
    \begin{array}{l l}
      \text{Step I:} &\XX^{[i_2]}(h_2):=\YY^{[0]}+h_{2} \e_{i_2}1_{\{\cdot \geq \mu_t^{[0]}\}},\\
    & \YY^{[i_2]}(h_2):=\reallywidetilde{\XX^{[i_2]}(h_2)},\\
    & \mu_t^{[i_2]}(h_2):=\psi_{R_2}^{-1}(\tau_{\XX^{[i_2]}(h_2),R_2}(\mu_t^{[0]})), \text{ for some pair } (R_2,\psi_{R_2}).
    \end{array}
\end{equation*}

\begin{figure}[H]
 \captionsetup{skip=1pt}
 \captionsetup{position=above}
    \centering
    \caption{Step I}\label{step1}
    \begin{minipage}{0.45\textwidth}
        \centering
         \caption*{$\XX^{[i_2]}(h_2)$}\includegraphics[width=\textwidth]{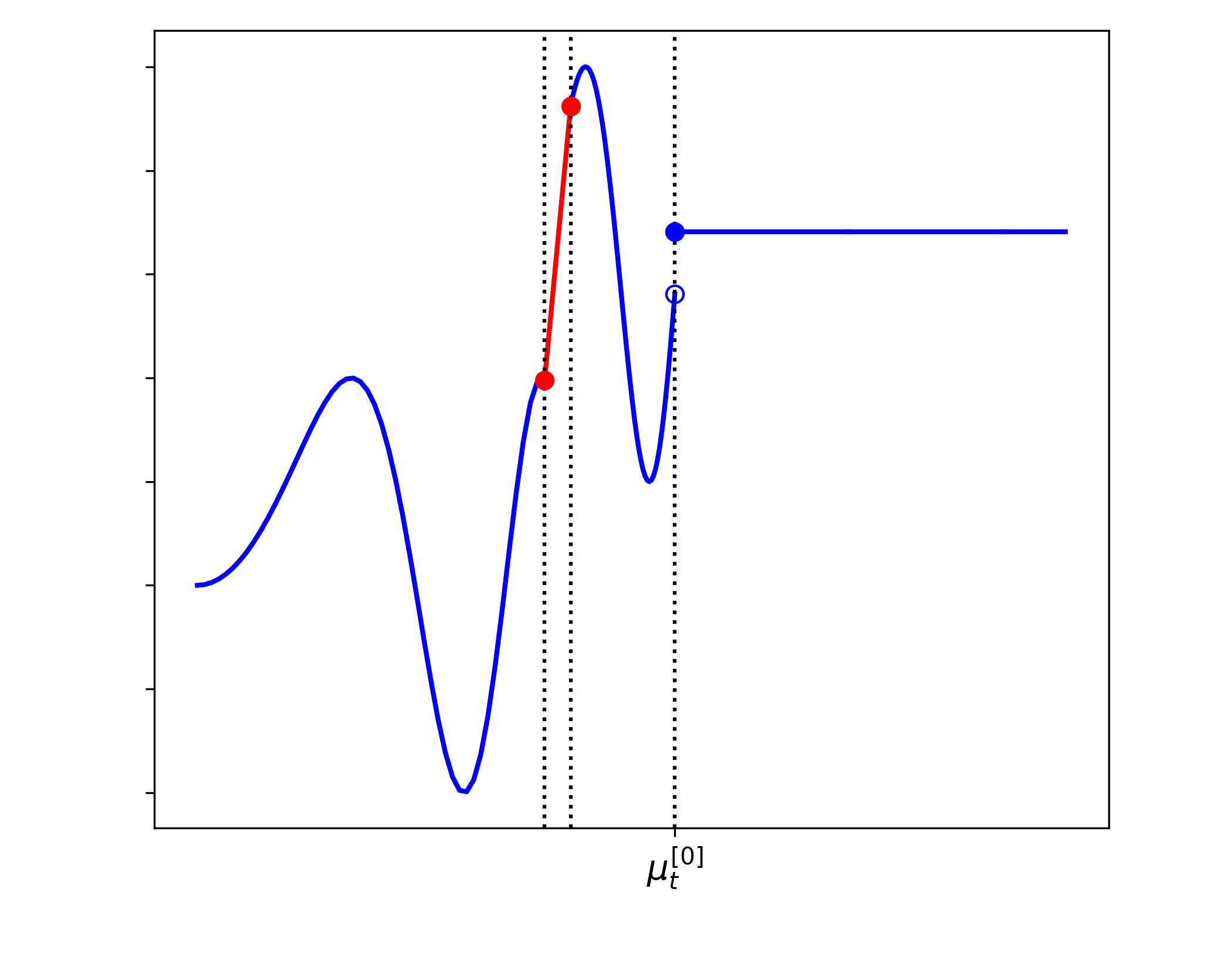}

    \end{minipage}
   \hspace{-3mm}
    \begin{minipage}{0.45\textwidth}  
        \centering
         \caption*{$\YY^{[i_2]}(h_2)$}\includegraphics[width=\textwidth]{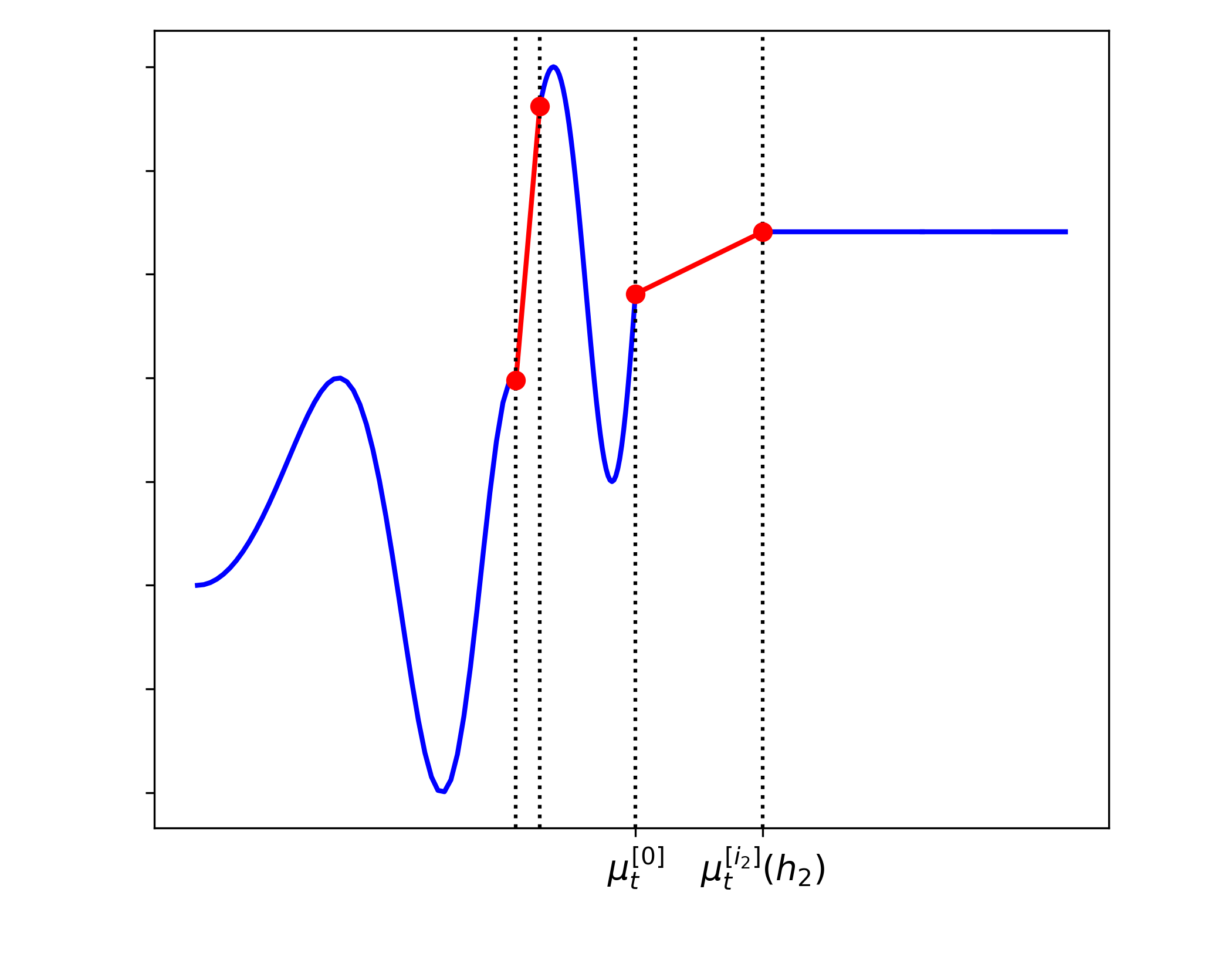}  
    \end{minipage}
\end{figure}

\begin{equation*}
    \begin{array}{l l}
   \text{Step II:}&  \XX^{[i_2,i_1]}(h_2,h_1):=\YY^{[i_2]}(h_{2})+h_{1} \e_{i_1}1_{\{\cdot \geq \mu_t^{[i_2]}(h_2)\}}, \\
    & \YY^{[i_2,i_1]}(h_2,h_1):=\reallywidetilde{\XX^{[i_2,i_1]}(h_2,h_1)},\\
   &  \mu_t^{[i_2,i_1]}(h_2,h_1):=\psi_{R_1}^{-1}(\tau_{\XX^{[i_2,i_1]}(h_2,h_1),R_1}(\mu_t^{[i_2]}(h_2)), \text{ for some pair } (R_1,\psi_{R_1}).
    \end{array}
\end{equation*}

\begin{figure}[H]
 \captionsetup{skip=1pt}
 \captionsetup{position=above}
 \caption{Step II}\label{step 2}
    \centering
    \begin{minipage}{0.45\textwidth}
        \centering
         \caption*{$\XX^{[i_2,i_1]}(h_2,h_1)$}\includegraphics[width=\textwidth]{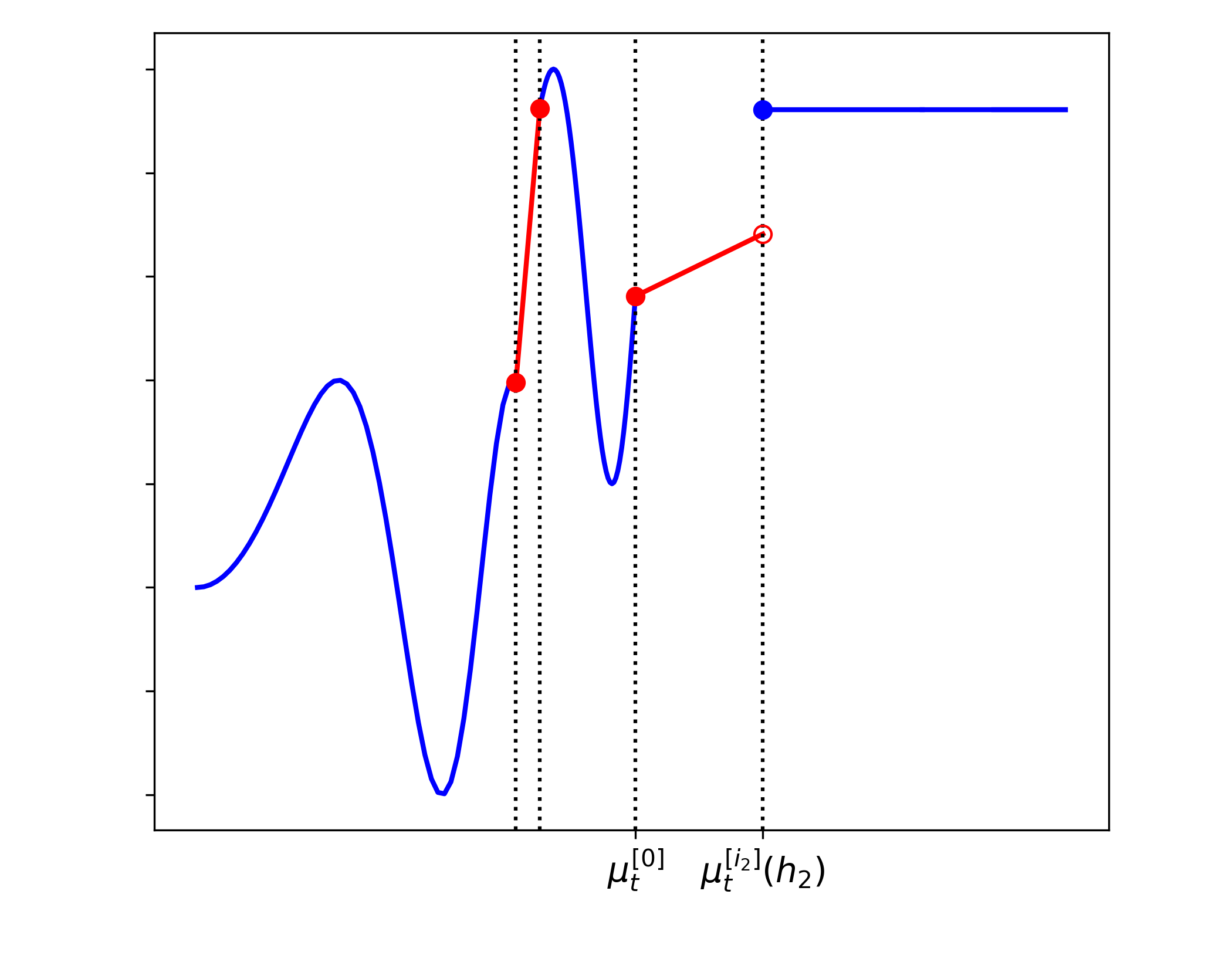}

    \end{minipage}
   \hspace{-3mm}
    \begin{minipage}{0.45\textwidth}  
        \centering
         \caption*{$\YY^{[i_2,i_1]}(h_2,h_1)$}\includegraphics[width=\textwidth]{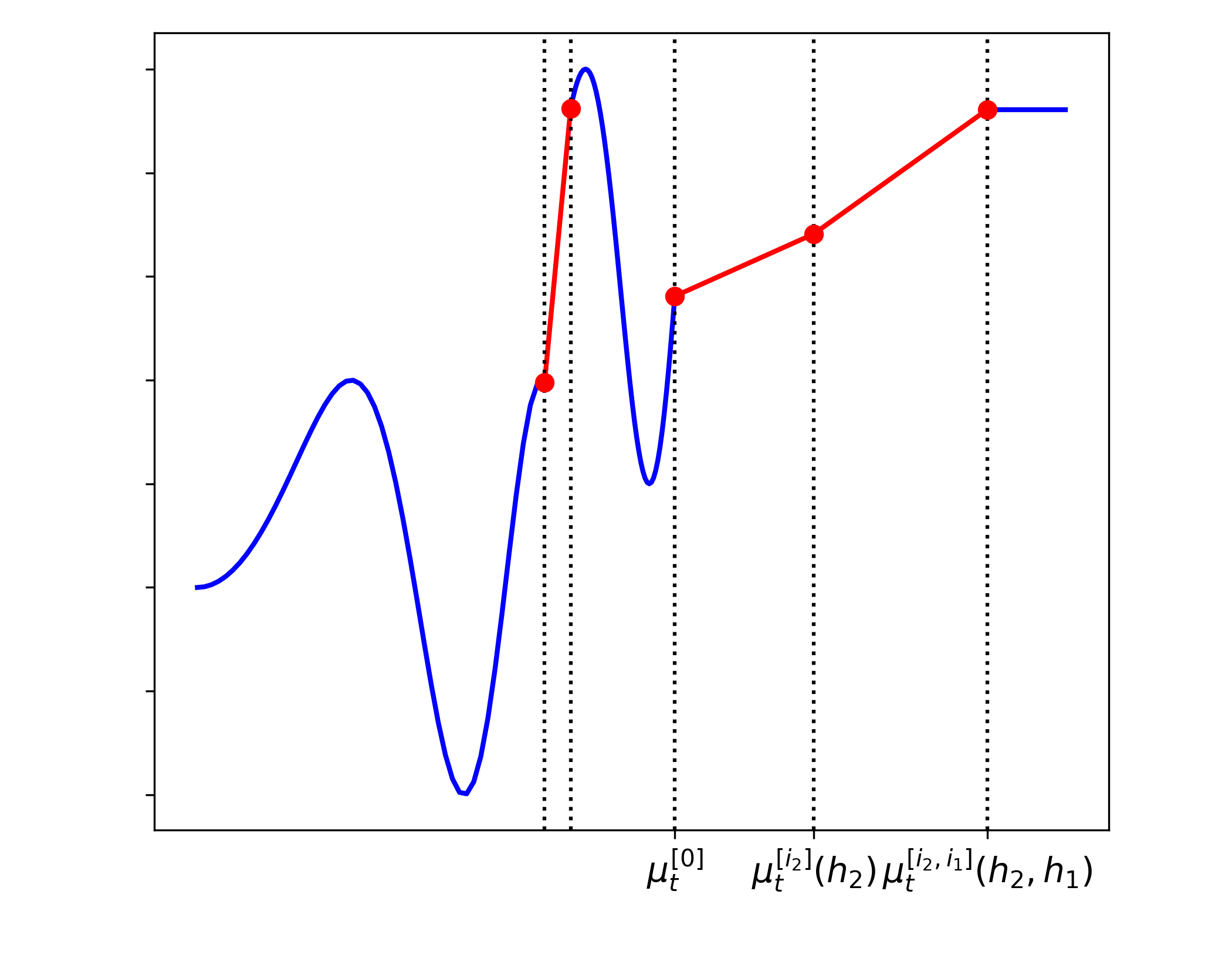}  
    \end{minipage}
\end{figure}
Notice that due to the one-dimensionality of the graphs, in Figures \ref{step0}, \ref{step1}, and \ref{step 2} we consider a one-dimensional path and thus $i_1=i_2$.

\begin{figure}[H]
 \captionsetup{skip=2pt}
 \setlength{\belowcaptionskip}{10pt}
  \caption{Vertically perturbed $2$-dimensional path  }\label{figure1}
 \captionsetup{position=above}
    \centering
    \begin{minipage}{0.45\textwidth}
        \centering
         \caption*{$\XX^t+h_2\e_{2}1_{\{\cdot \geq t+\varepsilon_2\}}+h_1\e_{1}1_{\{\cdot \geq t+\varepsilon_2+\varepsilon_1\}}$}
        \includegraphics[width=\textwidth]{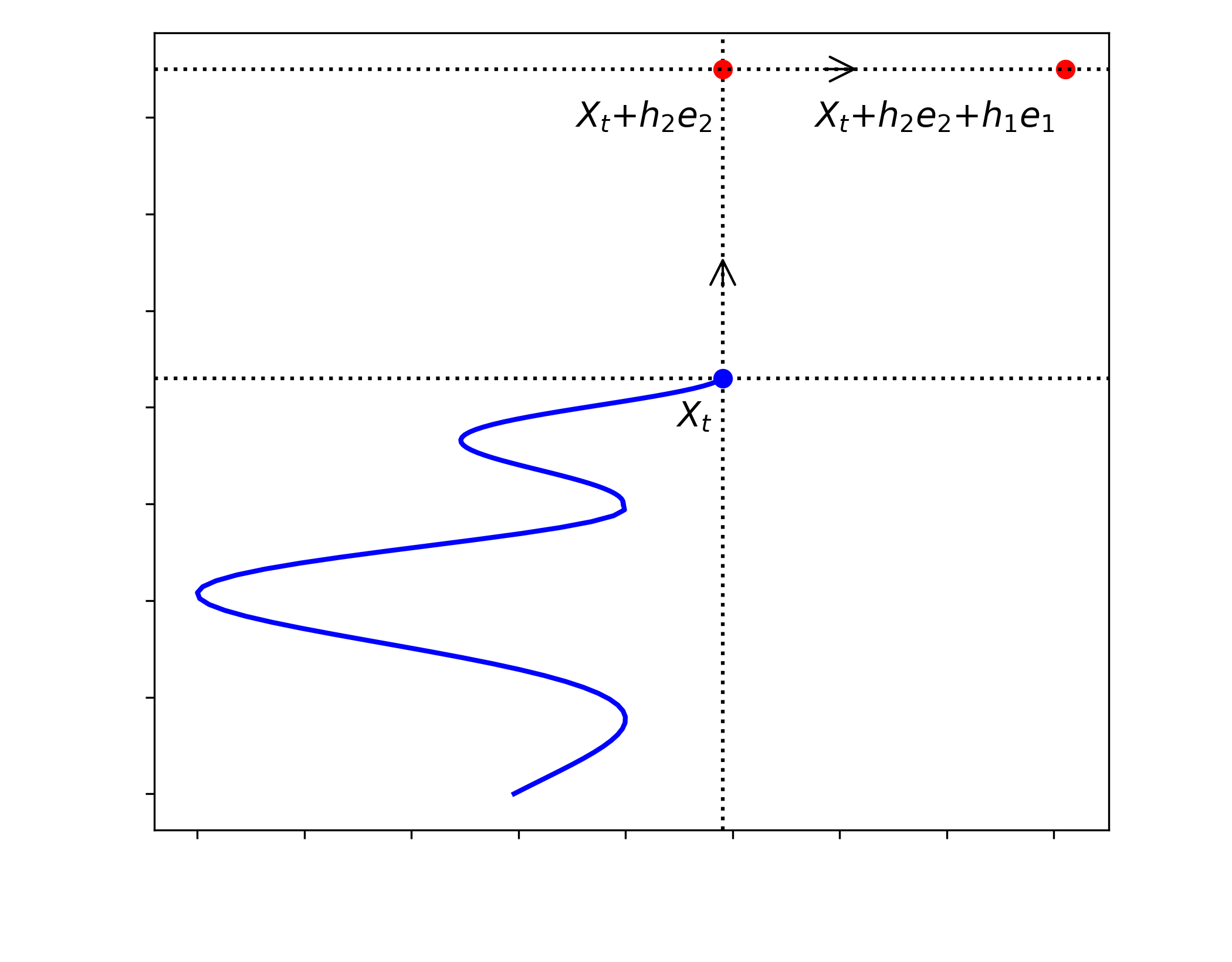}

    \end{minipage}
   \hspace{-3mm}
     \begin{minipage}{0.45\textwidth}
        \centering
         \caption*{$\XX^t+h_1\e_{1}1_{\{\cdot \geq t+\varepsilon_1\}}+h_2\e_{2}1_{\{\cdot \geq t+\varepsilon_1+\varepsilon_2\}}$}
        \includegraphics[width=\textwidth]{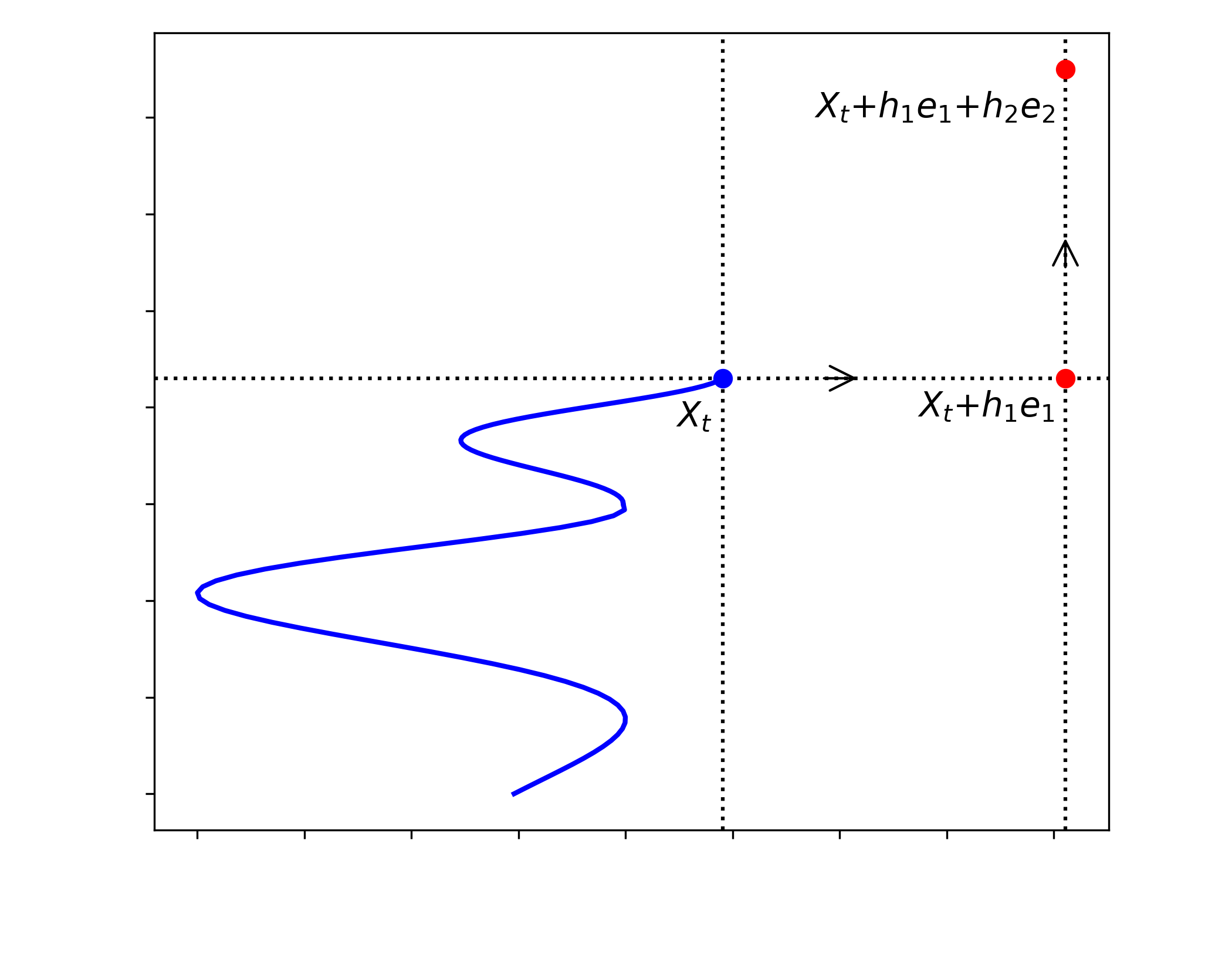}

    \end{minipage}
\end{figure}

Then, for $i_1,i_2=1,\dots,d$,
\begin{align*}
   U^{{i_1}}F(t,\XX)=& \frac{d}{ dh_1}F(\mu_t^{[i_2,i_1]}(0,h_1),\YY^{[i_2,i_1]}(0,h_1))|_{h_1=0},\\
  U^{i_{2}} U^{i_{1}}F(t,\XX)=&\frac{d^2}{dh_{2}dh_1}F(\mu_t^{[i_2,i_1]}(h_2,h_1),\YY^{[i_2,i_1]}(h_2,h_1))|_{h_{1}=h_2=0}.
\end{align*}
In particular, if $i_2\neq i_1$, $ U^{i_{2}} U^{i_{1}}F(t,\XX)$ and $ U^{i_{1}} U^{i_{2}}F(t,\XX)$ are computed by evaluating the functional at different paths, that is, 
\begin{align}\label{eq: second derivative 1}
    U^{i_{2}} U^{i_{1}}F(t,\XX)=&\frac{d^2}{dh_{2}dh_1}F(\mu_t^{[i_2,i_1]}(h_2,h_1),\YY^{[i_2,i_1]}(h_2,h_1))|_{h_{1}=h_2=0}, \\
     U^{i_{1}} U^{i_{2}}F(t,\XX)=&\frac{d^2}{dh_{1}dh_2}F(\mu_t^{[i_1,i_2]}(h_1,h_2),\YY^{[i_1,i_2]}(h_1,h_2))|_{h_{1}=h_2=0}, \nonumber 
\end{align}

with  $\YY^{[i_2,i_1]}(h_2,h_1))$ and $\YY^{[i_1,i_2]}(h_1,h_2))$ being different.  This is the crucial point as it is precisely for this reason that $ U^{i_{2}} U^{i_{1}}F(t,\XX)$ and $ U^{i_{1}} U^{i_{2}}F(t,\XX)$ are not necessarily equal, implying that the mixed vertical derivatives do not commute. One may notice that the derivatives with respect to $h_1$ and $h_2$ commute in each of the equations in  \eqref{eq: second derivative 1} if the second order partial derivatives are continuous. However, this is not relevant.
Indeed, the computation of different mixed vertical derivatives is not about reversing the order of differentiation with respect to $h_1$ and $h_2$, but, instead, requires evaluating the functional at different paths. This difference becomes  more evident when recognizing that $\YY^{[i_2,i_1]}(h_2,h_1)$ and $\YY^{[i_1,i_2]}(h_1,h_2)$  are time-stretched versions of some Marcus transformation of the paths
\begin{align}\label{eq: delayed perturbed path}
    &\XX^t+h_2\e_{i_2}1_{\{\cdot \geq t+\varepsilon_1\}}+h_1\e_{i_1}1_{\{\cdot \geq t+\varepsilon_1+\varepsilon_2\}},\\
    &\XX^t+h_1\e_{i_1}1_{\{\cdot \geq t+\varepsilon_1\}}+h_2\e_{i_2}1_{\{\cdot \geq t+\varepsilon_1+\varepsilon_2\}}, \nonumber 
\end{align}
respectively, and that
by conditions \ref{item def Mcpf ii},\ref{item def Mcpf iii} in Definition \ref{3def: Marcus canonical path functional},
\begin{align*}
    U^{i_{2}} U^{i_{1}}F(t,\XX)=&\frac{d^2}{dh_{2}dh_1}F(t+\varepsilon_1+\varepsilon_2,\XX^t+h_2\e_{i_2}1_{\{\cdot \geq t+\varepsilon_1\}}+h_1\e_{i_1}1_{\{\cdot \geq t+\varepsilon_1+\varepsilon_2\}})|_{h_{1}=h_2=0}.\\
U^{i_{1}} U^{i_{2}}F(t,\XX)=&\frac{d^2}{dh_{1}dh_2}F(t+\varepsilon_1+\varepsilon_2,\XX^t+h_1\e_{i_1}1_{\{\cdot \geq t+\varepsilon_1\}}+h_2\e_{i_2}1_{\{\cdot \geq t+\varepsilon_1+\varepsilon_2\}})|_{h_{1}=h_2=0}.
\end{align*}
Therefore, computing the mixed vertical derivatives reduces in fact to considering different delayed perturbations. In particular, in the calculation of $U^{i_{2}} U^{i_{1}}F(t,\XX)$ the  path is vertically perturbed first in the direction of the canonical vector $\e_{i_2}$, followed by a perturbation in the direction of $\e_{i_1}$, which results in
\begin{align*}
  \frac{d}{dh_1}F(t+\varepsilon_1+\varepsilon_2,\XX^t+h_2\e_{i_2}1_{\{\cdot \geq t+\varepsilon_1\}}+h_1\e_{i_1}1_{\{\cdot \geq t+\varepsilon_1+\varepsilon_2\}})|_{h_{1}=0}&=U^{i_{1}}F(t+\varepsilon_1,\XX^t+h_2\e_{i_2}1_{\{\cdot \geq t+\varepsilon_1\}})\\
    \frac{d}{dh_2}U^{i_{1}}F(t+\varepsilon_1,\XX^t+h_2\e_{i_2}1_{\{\cdot \geq t+\varepsilon_1\}})|_{h_2=0}&= U^{i_{2}} U^{i_{1}}F(t,\XX).
\end{align*}
In the calculation of $U^{i_{1}} U^{i_{2}}F(t,\XX)$ instead, the order of the vertical perturbation of the path is reversed. To clarify this argument further, we illustrate in Figure \ref{figure1} the vertical perturbations of a 2-dimensional path. The image on the left shows the perturbation required for computing $U^{2} U^{1}F(t,\XX)$, while the image on the right the perturbation for computing $U^1 U^{2}F(t,\XX)$.

\begin{example}\label{example non commutative}
Let us consider the path functional introduced in Example~\ref{example: Marcus canonical}\ref{item examples Marcus iii}, that is
\begin{align*}
     F(t,\XX):=\int_0^t \XX^1_{s^-}d\XX^2_s +\frac{1}{2}\sum_{0<s\leq t}\Delta \XX^1_s\Delta \XX^2_s,
 \end{align*}
 for all $(t,\XX)\in D^1([0,1],G^1(\Rset^{2}))$.

A direct computation of the iterative procedure described above shows that for all $(t,\XX)\in [0,1]\times D([0,1],G^{1}(\Rset^{2}))$, one gets $$U^2U^2F(t,\XX)=U^1U^1F(t,\XX)=U^2U^1F(t,\XX)=U^1F(t,\XX)=0,$$ whereas, $$U^2F(t,\XX)=\XX^1_t \quad \text{and}\quad U^1U^2F(t,\XX)=1.$$

\end{example}

\begin{remark}\label{rem: non commutativity derivate}
The above considerations, in particular Example \ref{example non commutative}, show that the notion of vertical derivatives of higher order introduced in Definition~\ref{def: higher order derivatives} establishes a differential calculus for path functionals that allows for non-commutative derivation orders. This is in contrast to the setup in \cite{D:09},~\cite{CRF:10},~\cite{CR:12} (see e.g., p. 133 in~\cite{CR:12}). 

\noindent This different behavior can be explained as follows. For a path functional that is non-anticipative, the second-order vertical derivatives at some $(t,\XX)$ as computed in \cite{D:09} (see also Definition 9 in \cite{CRF:13}) explicitly read as   
\begin{align}\label{eq: derivative Dupire}
\partial_{i_2}\partial_{i_1}F(t,\XX)=&\frac{d^2}{dh_{2}dh_1}F(t,\XX^t+(h_2\e_{i_2}+h_1\e_{i_1})1_{\{\cdot \geq t\}})|_{h_{1}=h_2=0},
\end{align}
for $i_1,i_2=1,\dots,d$. In contrast, in the present framework,
\begin{align*}
    U^{i_{2}} U^{i_{1}}F(t,\XX)=\frac{d^2}{dh_{2}dh_1}F(t+\varepsilon_2+\varepsilon_1,\XX^t+h_2\e_{i_2}1_{\{\cdot \geq t+\varepsilon_2\}}+h_1\e_{i_1}1_{\{\cdot \geq t+\varepsilon_2+\varepsilon_1\}})|_{h_{1}=h_2=0}.
\end{align*}
As already argued above this means that, if $i_1\neq i_2$, $ U^{i_{2}} U^{i_{1}}F(t,\XX)$ and $ U^{i_{1}} U^{i_{2}}F(t,\XX)$ are computed by evaluating the functional at different paths. This differs from the computation in \eqref{eq: derivative Dupire} where, since the perturbation occurs at the same time, the functional is evaluated always at the same path. Observe furthermore that $$\partial_{i_2}\partial_{i_1}F(t,\XX)=\partial^2_{\xi_{i_{1}}\xi_{i_2}}G(\xi)|_{\xi=0},$$ for $G(\xi):=F(t,\XX+\xi1_{\{\cdot \geq t\}})$, $\xi \in \Rset^d$. Therefore, by Schwarz's theorem the mixed vertical derivatives as computed in \cite{D:09} are equal as long as the maps $\xi\mapsto \partial^2_{\xi_{i_{1}}\xi_{i_2}}G(\xi)$ are continuous in a neighborhood of $0$. Other considerations for accommodating the non-communicative nature of the higher-order derivatives have been made by \cite{DS19}.

The non-commutative behavior that we obtain here is particularly relevant when dealing with weakly geometric $p$-rough paths, for $p\in[2,3)$,  where the second order vertical derivatives appear in the \Ito-formula (see Sections~\ref{sec: into formula p23}). 
It is, however, also crucial for finite variation paths in view of the Taylor expansions (see Theorem~\ref{th: taylor 12}).

\end{remark}

\subsection{Var-continuous path functionals}

In this section, we introduce the necessary regularity conditions for the path functionals to derive the results presented in the subsequent sections. Unless otherwise specified,  we let $p\in [1,3)$ and for $q\geq p$, we denote by $d_q$ the $q$-variation pseudometric on the space of \cadlag paths as defined in Section~\ref{sec: rough paths and rough integration}.

\begin{definition}\label{3def: var continuous funcitonals}
    Fix  $m\in \Nset_0$. We say that a path functional $F\in ((\Mcalp^0)^d)^{\otimes m}$ is \textit{$[p,[p]+1)$\text{-}var continuous} if for every $(\XX^n, \XX)_{n\in\Nset}\subset C^p([0,1],G^{[p]}(\Rset^{d}))$ such that for some  $q\in [p,[p]+1)$
    \begin{align*}
        \lim_{n\rightarrow\infty}d_q(\XX^n,\XX)=0,
    \end{align*}
 it holds that 
    \begin{align*}
\lim_{n\rightarrow\infty}d_{q'}(F(\cdot,\XX^n),F(\cdot,\XX))=0,
    \end{align*}
    for all $q'>q$.
    \end{definition}

    \begin{definition}
 Fix $m\in \Nset_0$. For $p\in [2,3)$ and $F\in ((\Mcal_2^0)^d)^{\otimes m}$, $F'\in ((\Mcal_2^0)^d)^{\otimes m+1}$, we say that a remainder path functionals $R^{F,F'}$, defined in Definition~\ref{3def: 2 parameter functional}, is \textit{$[p,3)$\text{-}var continuous} if for every $(\XX^n, \XX)_{n\in\Nset}\subset C^p([0,1],G^{2}(\Rset^{d}))$ such that for some  $q\in [p,3)$
    \begin{align*}
        \lim_{n\rightarrow\infty}d_q(\XX^n,\XX)=0,
    \end{align*}
 it holds that
    \begin{align*}
\lim_{n\rightarrow\infty}d_{q'/2}(R((\cdot,\cdot),\XX^n),R((\cdot,\cdot),\XX))=0,
    \end{align*}
    for all $q'>q$.

\end{definition}
\subsection{Linear functionals of the signature}\label{3sec: linear functions as Marcus}
This section is devoted to the study of linear functionals of the signature.
\begin{definition}
    Fix $\u\in T(\Rset^d)$. We call the path functional  
    \begin{align}\label{3eq: linear functional as Marcus}
    F^\u:&[0,1]\times D^p([0,1],G^{[p]}(\Rset^d))\rightarrow\Rset \\ \nonumber 
    & \quad (t,\XX)\longmapsto F^\u(t,\XX):=\langle \u,\X_t\rangle
\end{align} 
linear functional of the signature.
\end{definition}
\noindent The proof of the following proposition is given in Appendix \ref{proof prop: derivative linear func}.
\begin{proposition}\label{prop: derivative linear func}
    Let $F^\u :[0,1]\times D^p([0,1],G^{[p]}(\Rset^d))\rightarrow \Rset$  be a linear functionals of the signature for some $\u\in T(\Rset^d)$. Then, 
    \begin{enumerate}
        \item  \label{item 1 Fu}$F^\u$ is a non-anticipative Marcus canonical  path functional. 
        \item  \label{item 2 Fu}$F^\u$ is infinitely many times vertically differentiable and for all $k\in \Nset$, 
        \begin{align}\label{3eq: derivatives linear functionals}
    \nabla^k F^\u(t,\XX)=\langle \u^{(k)},\X_t\rangle,
\end{align}
for all $(t,\XX)\in [0,1]\times D([0,1],G^{[p]}(\Rset^d))$.
\item \label{item 3 Fu}$F^\u$ is $[
p,[p]+1)$-var continuous.
    \end{enumerate}
\end{proposition}

Now, a direct combination of Proposition~\ref{prop: derivative linear func} and equations~\eqref{eq: u,signature p12} and~\eqref{eq: u,signature p23} yields the functional \Ito-formula for linear functionals of the signature. 
\begin{lemma}\label{lemma: into linear func}
    Let $F^\u :[0,1]\times D^p([0,1],G^{[p]}(\Rset^d))\rightarrow \Rset$  be a linear functional of the signature for some $\u\in T(\Rset^d)$. Then, the following functional \Ito-formulas hold.
    \begin{enumerate}
      \item If $p\in [1,2)$,
      \begin{align*}
	F^\u(t,\XX)=&F^\u(0,\XX)+ \int_0^t \nabla F^\u(s^-,\XX) d\XX_{s}\\
  &+\sum_{0<s\leq t}F^\u(s,\XX) -F^\u(s^-,\XX)-\Delta F^\u(s^-,\XX) \Delta X_{s}.  \nonumber 
		\end{align*}
The integral term 
 is a Young integral of $\nabla F^\u(\cdot,\XX) \in D^p([0,1],\Rset^{d+1})$ with respect to $\XX$.
  \item If $p\in [2,3)$,
      \begin{align*}
		F^\u(t,\XX) =&  F^\u(0,\XX)+\int_0^t F^\u(s^-,\XX)d\XX_s\\
  &+\sum_{0<s\leq t}F^\u(s,\XX) -F^\u(s^-,\XX) -\nabla F^\u(s^-,\XX) \Delta X_{s}-\nabla^2F^\u(s^-,\XX) \Delta \X^{(2)}_{s}.\nonumber 
		\end{align*}
The integral term is a rough integral of  $\big (\nabla F^\u(\cdot,\XX),\nabla^2 F^\u(\cdot,\XX)\big) \in \mathcal{V}^{p,\frac{p}{2}}_{\XX}$  with respect to $\XX$.
  \end{enumerate}

\end{lemma}

\section{Universal approximation theorem for vertically differentiable path functionals}\label{sec: uat differential}
In this section, we exploit the approximation properties of the signature of weakly geometric \cadlag $p$ rough paths to derive a universal approximation theorem (UAT) for non-anticipative Marcus canonical path functionals which are  vertically differentiable (Theorem \ref{3th: Nachbin for path functionals}).

This result for path functionals is analogous to the classic Nachbin theorem, which provides an analog to the Stone-Weierstrass theorem for algebras of $C^K$ functions on  $\Rset^d$, for $K\in \Nset$ (see~\cite{N:49}).

\subsection{The subspace of the tracking-jumps-extended paths}

First, we introduce the subset of paths considered for deriving the UAT. We call them \textit{tracking jumps-extended paths}. The distinctive feature of these paths is that they always admit a Marcus transformed path that is time-extended (see Remark~\ref{3rem: Marcus tracking-jumps paths}\ref{ii}). This is a key aspect both for the development of the proof of Theorem~\ref{3th: Nachbin for path functionals} and the main results in Sections~\ref{3sec: functional Ito formula} and~\ref{sec: functional taylor}.

The notion of tracking jumps-extended paths relies on the more general one of 
extended weakly geometric \cadlag rough path, which we give below. Roughly speaking this space consists of all the weakly geometric  \cadlag  $p$-rough paths with values in $G^{[p]}(\Rset^{d+1})$ obtained through the addition of an auxiliary path component of finite variation to some $\XX\in D^p([0,1],G^{[p]}(\Rset^{d}))$. When $p\in [1,2)$, this extension is quite straightforward. When $p\in [2,3)$ instead, 
one needs to ensure consistency between the Young and the (level 2) rough integration (see Section \ref{sec: rough=young}), while preserving the group-valued constraint. For the explanation of the index $0$ used here we refer to Notation~\ref{notation: index 0} below.
\begin{definition}\label{def: Marcus extended rough paths}
    Fix $Z\in D^1([0,1],\Rset)$ and $\widehat{\XX}\in D^p([0,1],G^{[p]}(\Rset^{d+1})).$
    \begin{enumerate}
        \item If $p\in [1,2)$, we say that $\wXX$ is a $Z$-extended weakly geometric \cadlag $p$-rough path if for all $i=1,\dots,d$, 
    \begin{align*}
        \langle \e_i,\wXX\rangle = \langle \e_i,\XX\rangle, \quad \text{and }\quad   \langle \e_{0},\wXX\rangle = Z,
    \end{align*}
     for some $\XX\in D^p([0,1],G^1(\Rset^d))$. 
     \item If $p\in [2,3)$,    
    we say that $\wXX$ is a $Z$-extended weakly geometric \cadlag $p$-rough path if for all $i,j=1,\dots,d$, 
    \begin{align*}
        \langle \e_i,\wXX\rangle = \langle \e_i,\XX\rangle, \qquad   \langle \e_{(ji)},\wXX\rangle = \langle \e_{(ji)},\XX\rangle, \qquad  \langle \e_0,\wXX\rangle = Z, 
    \end{align*}
    and for all $i=0,\dots,d$,
     \begin{align}\label{eq: components extended paths}
   \langle \e_{(i,0)},\wXX\rangle =\int_0^\cdot  \langle \e_i,\wXX_{0,s^-}\rangle dZ_s+\frac{1}{2} \sum_{0<s\leq \cdot}\Delta  \langle \e_i,\wXX_s\rangle \Delta Z_s.
    \end{align}
    for some $\XX\in D^p([0,1],G^2(\Rset^d))$, and the integral in~\eqref{eq: components extended paths} is of Young type.
    \end{enumerate}

\end{definition}

\begin{notation}\label{notation: index 0}
Unless explicitly mentioned, we use throughout the paper the index $0$ to denote the additional auxiliary component $Z$ of $\wXX$.
\end{notation}

\begin{remark}\phantomsection\label{rem: extended paths}
\begin{enumerate}

\item Let $p\in [2,3)$. Notice that by definition any $Z$-extended weakly geometric \cadlag $p$-rough path $\wXX$ lies in $G^2(\Rset^{d+1})$. Therefore, by the 
shuffle property~\eqref{3eq:shuffle}, for all $i=0,1,\dots,d$, $t\in [0,1]$,  $$\langle \e_{(0,i)},\widehat{\mathbf{X}}_t\rangle=Z_t \langle \e_{i},\widehat{\mathbf{X}}_t\rangle-\langle \e_{(i,0)},\widehat{\mathbf{X}}_t\rangle.$$
Moreover, since $Z\in D^1([0,1],\R)$ and $\langle \e_i,\XX\rangle \in D^p([0,1],\R)$, the Young integral in~\eqref{eq: components extended paths} is well defined and for all $t\in [0,1]$, 
\begin{align*}
  \sum_{0<s\leq t}|\Delta \langle \e_i,\wXX_s\rangle \Delta Z_s|\leq \left(\sum_{0<s\leq t}|\Delta \langle \e_i,\wXX_s\rangle|^{3}\right)^{1/3}\left(\sum_{0<s\leq t}| \Delta Z_s|^{\frac{3}{2}}\right)^{2/3}<\infty.
  \end{align*}
This implies that the path defined in equation~\eqref{eq: components extended paths} is of finite variation and the definition is well-posed. 

\item \label{item 4 remark extended path} One can see that given $Z\in D^1([0,1],\R)$ and $\XX\in D^p([0,1],G^{[p]}(\Rset^d))$, it is always possible to build a $Z$-extended weakly geometric $p$-rough path $\wXX$ as in Definition~\ref{def: Marcus extended rough paths}. The reverse is trivially true.

\end{enumerate}
    
\end{remark}

\noindent Now, we are ready to introduce the set of tracking jumps-extended paths.

 \begin{definition}\label{3def: tracking-jumps paths}
    We say that $\wXX\in {D}^p([0,1],G^{[p]}(\Rset^{d+1}))$  is a \textit{tracking-jumps-extended weakly geometric \cadlag $p$-rough path} if it is the $Z$-extension of some $\XX\in {D}^p([0,1],G^{[p]}(\Rset^{d}))$, for a \cadlag strictly increasing  piecewise linear function
    $ Z:[0,1]\rightarrow[0,1]$ such that $Z_0=0$, $Z_1=1$, and
 \begin{align*}
     \{s\in (0,1] \ : \ \Delta Z_s\neq 0\}= \{s\in (0,1] \ : \ \Delta \XX_s\neq 0\}.
 \end{align*}
 We call such a path $Z$ 
 \textit{tracking-jumps-path} of $\XX$ and denote the set of such $Z$-extended \cadlag 
 
 (continuous) paths by $\widehat{D}^p([0,1],G^{[p]}(\Rset^{d+1}))$ ($\widehat{C}^p([0,1],G^{[p]}(\Rset^{d+1}))$). 
 \end{definition}
 \begin{remark}\label{3rem: Marcus tracking-jumps paths}
 \begin{enumerate}
   \item \label{3item rem tracking-jumpsp ii} Given $\XX\in {D}^p([0,1],G^{[p]}(\Rset^{d}))$ it is always possible to build (non-uniquely) a tracking-jumps path associated with it. If $\XX$ is continuous, the tracking-jumps path is unique and given by the identity path $\text{Id}_u:=u$ for all $u\in [0,1]$. For $\XX$ \cadlag, there exists more than one tracking-jumps-extension $\wXX$.  We refer to every Id-extended continuous path as \textit{time-extended path}.
   
    Notice in particular that the existence of a tracking-jumps-extension is related only to the \cadlag property of the original path. 
   
     \item \label{ii}  Observe  that for all $\wXX\in \widehat{D}^p([0,1],G^{[p]}(\Rset^{d+1}))$  there always exists a Marcus-transformed path of $\wXX$ with respect to some pair $(R,\psi_R)$, which is a time-extended weakly geometric continuous $p$-rough path.
Indeed, since a tracking-jumps path is simply a strictly increasing piecewise linear path $Z$ such that $Z(0)=0$ and $Z(1)=1$ and jumps whenever the other components of the path do,  a pair $(R,\psi_R)$ can be chosen to ensure a Marcus-transformed path for which the transformed component $\widetilde{Z}$ satisfies $\widetilde{Z}_u=u$ for all $u\in [0,1]$.

 \end{enumerate}
 
 \end{remark}

\subsection{UAT for vertically differentiable path functionals} 
To present the main result of the section (Theorem~\ref{3th: Nachbin for path functionals}), we shall introduce some notation and new definitions. Since we will be mostly concerned with paths with values in $G^{[p]}(\Rset^{d+1})$, we introduce the relevant notions by considering directly the space $\Rset^{d+1}$.

\begin{notation}\label{3notation derivatives lie algebra}
     Let $K\in \Nset$ and $\mathfrak{g}^K(\Rset^{d+1})$ be the free step-$K$ nilpotent Lie algebra over $\Rset^{d+1}$ (see Section~\ref{3sec: algebra}). Recall that for every $  \bxi:=(0,\bxi^{(1)},\dots,\bxi^{(K)}) \in \mathfrak{g}^K(\Rset^{d+1})$, $$\bxi^{(j)}:=\pi_{j}(\bxi)\in \underbrace{[\Rset^{d+1},[\Rset^{d+1},\dots ,[\Rset^{d+1},\Rset^{d+1}]]]}_{(j-1)\text{ brackets}},$$  for $j=1,\dots, K$. Set $M:=\text{dim}(\mathfrak{g}^K(\Rset^{d+1}))=\sum_{j=1}^K M_j$ and 
\begin{align*}
    M_j:=\text{dim}(\underbrace{[\Rset^{d+1},[\Rset^{d+1},\dots ,[\Rset^{d+1},\Rset^{d+1}]]]}_{(j-1)\text{ brackets}}).
\end{align*}
For $\beta \in \Nset_0^{M}$, we write $\beta=(\beta_1,\dots,\beta_K)$, with $\beta_j\in \Nset_0^{M_j}$, and set
\begin{align}\label{3eq: beta weighted}|\beta|:=\sum_{j=1}^K  \sum_{i=1}^{M_j}\beta^i_j,\qquad \text{and}\quad |\beta|_{\mathfrak{g}^K(\Rset^{d+1})}:=\sum_{j=1}^K  \sum_{i=1}^{M_j}j\beta^i_j. 
\end{align}Furthermore, let  $\partial^{\beta_j}_{\bxi^{(j)}}$ denote the differential operator  such that 
$$\partial^{\beta_j}_{\bxi^{(j)}}f(\bxi):=\frac{\partial^{|\beta_j|}f}{\partial ^{\beta_j^1}\bxi^{(j),1}\dots\partial ^{\beta_{j}^{M_j}}\bxi^{(j),M_j} }(\bxi),$$ 
for $\bxi^{(j),l}$ denoting the $l$-th component of $\bxi^{(j)}$, $|\beta_j|:=\sum_{i=1}^{M_j}\beta^i_j$, and 
$f:\mathcal{U}\rightarrow \R$ some sufficiently regular map defined on some open set $\mathcal{U}\subseteq \mathfrak{g}^K(\Rset^{d+1})$,  with $\bxi\in \mathcal{U}$. 

\noindent Finally, denote by
$D^\beta_{\bxi}:=\partial_{\bxi^{(1)}}^{\beta_1}\dots \partial_{\bxi^{(K)}}^{\beta_K}$ the operator given by a consecutive application of $\partial_{\bxi^{(1)}}^{\beta_1},\dots, \partial_{\bxi^{(K)}}^{\beta_K}$.

\end{notation}

 \noindent Recall that $G((\Rset^{d+1}))$ denotes the set of group-like elements, introduced in equation~\eqref{3eq: G group inf}.

\begin{assumpt}\label{assumpt 1}
     $F:[0,1]\times D^{p}([0,1],G^{[p]}(\Rset^{d+1}))\rightarrow \Rset$ is a Marcus canonical  non-anticipative path functional such that for all $(t,\XX)\in [0,1]\times {C}^p([0,1],G^{[p]}(\Rset^{d+1}))$,
    \begin{align}\label{eq: F=g assumpt}
        F(t,\XX)=g(\X_t),
    \end{align}
    for some $g:G((\Rset^{d+1}))\rightarrow\Rset$, and $\X$ denoting the signature of $\XX$.
\end{assumpt}
\begin{remark}\label{remark F tilde F}
\begin{enumerate}
    \item 
    By the definition of Marcus canonical  path functionals and the signature of \cadlag rough paths, the equality~\eqref{eq: F=g assumpt} is valid in fact on $[0,1]\times {D}^{p}([0,1],G^{[p]}(\Rset^{d+1}))$.

 \item \label{remark F tilde F ii} Let $S$ be the map such that for all $\XX\in {C}^p([0,1],G^{[p]}(\Rset^{d+1})) $, $S(\XX):=\X_1$. For $t\in [0,1]$, let $\widehat{C}^p_t([0,1],G^{[p]}(\Rset^{d+1}))$ denote the space of time-extended paths stopped at time $t$. One can show that on  the set
 \begin{align*}
     \bigcup_{t\in [0,1]}\widehat{C}_t^p([0,1],G^{[p]}(\Rset^{d+1}))
 \end{align*}
 $S$ is an injective map (see e.g., the proof of Proposition 3.6 in \cite{CPS:22}). Therefore, the restriction of 
 the map $S$ to its image, denoted as $\mathcal{S}$, is a bijection.
 Furthermore, consider the map $\widetilde{F}:   \bigcup_{t\in [0,1]}\widehat{C}_t^p([0,1],G^{[p]}(\Rset^{d+1}))\rightarrow\Rset $ given by $ \widetilde{F}(\wXX^t):=  F(t,\wXX^t)$
and notice that $F$ admits a representation of the form
 \begin{align}\label{eq: bar g}
     F(t,\wXX^t)=\bar g(\wX_t),
 \end{align}
 for $\bar{g}:\mathcal{S}\rightarrow \Rset$ given by $\bar{g}:=\widetilde{F}\circ S^{-1}$ and every (stopped) time-extended continuous path $\wXX^t$. {\color{blue} Notice that the above reasoning applies in fact to the larger set of \cadlag paths}.
 However, in general, the set $\mathcal{S}\subsetneq G((\Rset^{d+1}))$ might not be big enough for the application of the reasonings in the proofs of Lemma~\ref{lemma: partial derivatives gk} and Theorem~\ref{3th: Nachbin for path functionals}, which involve some ideas from Lie group theory (see \cite{BLU:07}).
 \end{enumerate}
    
\end{remark}
\noindent The proof of the following lemma is given in Appendix~\ref{proof of lemma: partial derivatives gk}.  
\begin{lemma}\label{lemma: partial derivatives gk}
    Let $F\in \Mcalp^K$, for some $K\in \Nset$. Fix $\XX\in[0,1]\times  C^p([0,1],G^{[p]}(\Rset^d))$ and denote by $\X$ its signature. Under Assumption~\ref{assumpt 1}, the map
    \begin{align}\label{eq: gk assumpt lemma}
        g^{\XX,K}:&[0,1]\times  \mathfrak{g}^K(\Rset^{d+1})\rightarrow\Rset\\
          & \qquad (t,\bxi)\mapsto g^{\XX,K}(t,\bxi):= g(\X_t\otimes \exp(\bxi,0\dots,0)) \nonumber
    \end{align}
 is well defined and its derivatives at zero $D^\beta_{\bxi}g^{\XX,K}(t,\bxi)|_{\bxi=0}$ exist for all $|\beta|_{\mathfrak{g}^K(\Rset^{d+1})}\leq K$. 
\end{lemma}
\noindent For the development of the following results, we need the map $g^{\XX,K}$ to satisfy some stronger continuity conditions in the next assumption.
\begin{assumpt}\label{assumpt 2}
   For all $ \XX\in C^p([0,1],G^{[p]}(\Rset^d))$, the map
\begin{align*}
       [0,1]\times \mathcal{U}(0)\ni  (t,\bxi)\mapsto D^\beta_{\bxi} g^{\XX,K}(t,\bxi)
    \end{align*}
    is jointly continuous for all $|\beta|_{\mathfrak{g}^K(\Rset^{d+1})}\leq K$, and  some open neighborhood $\mathcal{U}(0)$ of $0\in \mathfrak{g}^K(\Rset^{d+1})$.

\end{assumpt}
    \begin{remark}
 Under Assumption~\ref{assumpt 2}, the map $\mathcal{U}(0)\ni \bxi\mapsto D^\beta g^{\XX,K}(t,\bxi)$ does not depend on the order of the partial derivatives (see Step 2 in the proof of Theorem~\ref{3th: Nachbin for path functionals}).
\end{remark}

\begin{definition}\label{3def: C^K path functionals}
    Let  $F:[0,1]\times D^{p}([0,1],G^{[p]}(\Rset^{d+1}))\rightarrow \Rset$. We write $F\in C^K$ if $F\in \Mcalp^K$ and Assumptions~\ref{assumpt 1} and~\ref{assumpt 2} are satisfied.
\end{definition}
\noindent An example of path functional satisfying the assumptions of Definition~\ref{3def: C^K path functionals} is given by the linear functionals of the signature $F^\u$, for some $\u \in T(\Rset^{d+1})$ (see Section~\ref{3sec: linear functions as Marcus}). Then, by Proposition~\ref{prop: derivative linear func} $F^\u\in C^K$ for all $K\in \Nset$, and the map~\eqref{eq: gk assumpt lemma} explicitly reads as
     \begin{align*}
          g^{\XX,K}(t,\bxi):= \langle \u,\X_t\otimes \exp(\bxi,0\dots,0)\rangle,
     \end{align*}
     for all $(t,\bxi)$ and $\XX$.

\begin{remark}\label{remark: CK functionals}
    For $m\in \Nset_0$, $K\in \Nset$ and a $(\Rset^d)^{\otimes m}$-valued
     path functional $F$, we write $F\in C^K$
     if its components
    are $C^K$. Notice that if $F\in C^K$, then for all $l=1,\dots,K$, $\nabla^lF\in C^{K-l}$.
\end{remark}

Now we state the main result of the section stating that any path functionals $F \in C^K$ when evaluated at a {tracking jumps}-extended path $\wXX$ (Definition~\ref{3def: tracking-jumps paths}), can be uniformly approximated in time, along with its derivatives, by linear functionals of the signature and its derivatives.

\begin{theorem}\label{3th: Nachbin for path functionals}
 Let $F:[0,1]\times D^{p}([0,1],G^{[p]}(\Rset^{d+1}))\rightarrow \Rset$. Assume $F\in C^K$.
Then, for all  $\widehat \XX\in \widehat{D}^p([0,1], G^{[p]}(\Rset^{d+1}))$ there exists $(\u_n)_{n\in \Nset}\in T(\Rset^{d+1})$ (possibly depending on $\widehat \XX$) such that
    \begin{align}\label{3eq: uat derivatives}
        \lim_{n\rightarrow\infty }\sup_{t\in [0,1]}\sum_{j=0}^K\|\nabla ^jF(t,\widehat \XX)-\nabla ^jF^{\u_n}(t,\widehat \XX)\|=0.
    \end{align}
\end{theorem}

 The proof of Theorem \ref{3th: Nachbin for path functionals} is given in Appendix~\ref{3sec: proof Nachbin}.

\section{Functional \Ito-formula}\label{3sec: functional Ito formula}
\subsection{The case \texorpdfstring{$p\in [1,2)$}{p in [1,2)}}
In this section, we present the functional \Ito-formula for maps of \cadlag rough paths of finite $p$-variation, for $p\in[1,2)$.

\begin{theorem}\label{3th: Ito formula p 12}
  Let $p\in[1,2)$ and $F:[0,1]\times D^{p}([0,1],G^{1}(\Rset^{d+1}))\rightarrow \Rset$ be a $C^2$-non-anticipative Marcus canonical  path functional such that $F$ and $\nabla ^{}F$ are $[p,2)$\text{-}var continuous. Then, for every $ \wXX\in \widehat{D}^{p}([0,1],G^{1}(\Rset^{d+1}))$ the path $\nabla F(\cdot,\wXX)$ is a \cadlag path of finite $p'$-variation, for all $p'> p$ such that 
  $\frac{1}{p'}+\frac{1}{p}>1$, 
  and for all
$t\in [0,1]$, 
    \begin{align}\label{3eq: ito formula p 12}
        F(t,  \wXX)-&F( 0, \wXX)=\int_0^t \nabla F(s^-,  \wXX)d \wXX_s\\
         &+\sum_{0<s\leq t}F(s,\wXX)-F(s^-,\wXX)-\nabla ^{}F(s^-,\wXX) \Delta \wXX_s. \nonumber 
    \end{align}
      The integral in~\eqref{3eq: ito formula p 12} is  a Young integral and the summation term is well defined as an absolutely summable series.
\end{theorem}
To prove Theorem~\ref{3th: Ito formula p 12} we make use of the following lemma, proved in Appendix \ref{proof lemma ito formula p12}.

\begin{lemma}\label{3lemma: ito formula uniformly bounded p 12}.
Let $p\in[1,2)$, $F:[0,1]\times D^{p}([0,1],G^{1}(\Rset^{d+1}))\rightarrow \Rset$ be such that $F\in \Mcal_1^2$ and $\widehat \XX\in  \widehat C^{1}([0,1],G^1(\Rset^{d+1}))$ a time-extended piecewise linear path. Assume that there exists $(\u_n)_{n\in\Nset}\subset T(\Rset^{d+1})$ such that 
\begin{align}\label{eq: approx F and derivatives p 12}
     \lim_{n\rightarrow\infty} \sup_{t\in [0,1]}\sum_{j=0}^2\|\nabla ^jF(t,\widehat \XX)-\nabla ^jF^{\u_n}(t,\widehat \XX)\|=0.
\end{align}
Then,  
\begin{enumerate}
    \item \label{3item: lemma 1 i} $\sup_{n\in \Nset}\|\nabla F^{\u_n}(\cdot,\widehat \XX)\|_{1\text{-}var}<\infty$;
    \item \label{3item: lemma 1 ii}  for all $p'>1$, $\lim_{n\rightarrow\infty}d_{p'}\big(\nabla F^{\u_n}(\cdot,\widehat \XX), \nabla F(\cdot,\wXX)\big)=0.$ 
\end{enumerate}
    \end{lemma}

\begin{proof}[Proof of Theorem~\ref{3th: Ito formula p 12}] We split the proof into three main steps. Step 1 proves the assertion for functionals evaluated at a time-extended piecewise linear path. Then, Step 2 extends it to the whole space of time-extended continuous paths by a density argument.  Finally, Step 3 extablishes the general result by exploiting that the functionals considered are of Marcus type and an adaptation of the proof of Theorem 38 in~\cite{FS:17}, which deals only with functional given as the solution of Marcus-RDE.

\paragraph{Step 1:} Let $\widehat \XX\in  \widehat C^{1}([0,1],G^1(\Rset^{d+1}))$ be a time-extended peicewise linear path. Since $F$ is a $C^2$-non-anticipative path functional, by  Theorem~\ref{3th: Nachbin for path functionals} there exists a sequence $(\u_n)_{n\in\Nset}\subset T(\Rset^{d+1})$ such that the convergence in equation \eqref{eq: approx F and derivatives p 12} holds true.
By Lemma~\ref{3lemma: ito formula uniformly bounded p 12}~\ref{3item: lemma 1 i} and Lemma 5.12 in~\cite{FV:10}, $  \|\nabla F(\cdot ,\wXX)\|_{1\text{-}var}<\infty$ and thus in particular 
\begin{align}\label{3eq: finite p' var }
  \|\nabla F(\cdot ,\wXX)\|_{p'\text{-}var}<\infty,
\end{align}
for all $p'\geq1$. 
Fix $p'>1$. By Lemma~\ref{3lemma: ito formula uniformly bounded p 12}~\ref{3item: lemma 1 ii} $\lim_{n\rightarrow\infty}d_{p'}\big(\langle \u_n^{(1)},\widehat{\mathbb{X}}_{\cdot}\rangle,\nabla F(\cdot ,\wXX)\big)=0.$ Moreover, by Proposition 6.11 in~\cite{FV:10}, for all $n\in \Nset$,
\begin{align*}
    d_1\bigg(&\int_0^\cdot \langle \u^{(1)}_n,\widehat{\mathbb{X}}_{s}\rangle d\widehat{\XX}_{s}, \int_0^\cdot \nabla F(s,\widehat{\XX}) d\widehat{\XX}_{s}\bigg) \\&
    \leq C \bigg( \dcc\big(\langle \u_n^{(1)},\widehat{\mathbb{X}}_0\rangle,\nabla F(0,\widehat{\XX})\big)+d_{p'}\big(\langle \u_n^{(1)},\widehat{\mathbb{X}}_{\cdot}\rangle,\nabla F(\cdot ,\wXX)\big)\bigg),
\end{align*}
for some $C>0$. Since for all $n\in \Nset$, $j=0,\dots,2$, $\nabla ^jF^{\u_n}(t,\widehat \XX)=\langle \u_n^{(j)},\widehat \X_t\rangle$, where $\wX$ denotes the signature of $\wXX$, and by  equation~\eqref{eq: u,signature p12} for every $n\in \Nset$ and $t\in [0,1]$,  $\langle \u_n,\widehat{\mathbb{X}}_{t}\rangle =  \langle \u_n,\widehat{\mathbb{X}}_{0}\rangle+\int_0^t \langle \u^{(1)}_n,\widehat{\mathbb{X}}_{s}\rangle  d\widehat{\XX}_{s}$, the claim follows.

\paragraph{Step 2:} Fix $\widehat \XX \in  \widehat C^{p}([0,1],G^{1}(\Rset^{d+1})) $. By Theorem 5.23 in~\cite{FV:10}, there exists a sequence  of piecewise linear time-extended paths $(\widehat \XX^M)_{M\in \N}$ such that 
    \begin{align}\label{3eq: approximation w.g.p 12}
        \lim_{M\rightarrow\infty }\sup_{t\in [0,1]}d_{CC}(\widehat \XX^M_t,\widehat \XX_t)=0, \qquad \text{and}\qquad    \sup_{M\in \Nset}\|\widehat \XX^M\|_{p\text{-}var}\leq C\|\widehat \XX\|_{p\text{-}var},
    \end{align}
    for some $C>0$. By Step 1, for every fixed $M$ and $t\in [0,1]$, it holds that 
        \begin{align*}
        F(t, \widehat \XX^M)-F(0, \widehat \XX^M)=\int_0^t \nabla^{} F(s, \widehat \XX^M)d\widehat \XX^M_s.
    \end{align*}
    By conditions~\eqref{3eq: approximation w.g.p 12} and interpolation (see Lemma 5.12 and Lemma 5.27 in~\cite{FV:10}), for all $p'>p$, $  \lim_{M\rightarrow\infty }d_{p'}(\widehat \XX^M, \widehat \XX)=0.$

Since by Step 1 (see equation~\eqref{3eq: finite p' var }) for all $M$, $\|\nabla ^{}F(\cdot,\widehat \XX^M)\|_{p'\text{-}var}<\infty$, for all $p'\geq1$, and  $\nabla ^{}F$ is  $[p,2)$\text{-}var continuous by assumption, $\|\nabla ^{}F(\cdot,\widehat \XX)\|_{p'\text{-}var}<\infty$ for all $p'>p$ such that $\frac{1}{p'}+\frac{1}{p}>1$.

 Fix $p''>p'>p$ such that  $\frac{1}{p''}+\frac{1}{p'}>1$ and observe that $\wXX,\wXX^M\in \widehat{C}^{p'}([0,1], G^1(\Rset^{d+1}))$,  for all $M\in \Nset$,  and since  $  \nabla F(\cdot,\widehat \XX),\nabla F(\cdot,\widehat \XX^M)\in {C}^{p''}([0,1], \Rset^{d+1})$ as well, the Young integral of the integrands $\nabla F(\cdot,\widehat \XX),\nabla F(\cdot,\widehat \XX^M)$ with respect to $\wXX,\wXX^M$, respectively, is well defined. Finally, by Proposition 6.11 in~\cite{FV:10}, for all $M\in \Nset$,
\begin{align*}
   d_{p'}\bigg(&\int_0^\cdot \nabla^{} F(s, \widehat \XX^M)d\widehat \XX^M_s,\int_0^\cdot  \nabla^{} F(s, \widehat \XX)d\widehat \XX_s\bigg)\\
    &\leq  {C}\bigg(d_{p'}(\widehat \XX^M, \widehat \XX)+d_{CC}\big(\nabla^{} F(0, \widehat \XX^M),\nabla F(0, \widehat \XX)\big)+d_{p''}\big(\nabla F(\cdot, \widehat \XX^M),\nabla F(\cdot, \widehat \XX)\big)\bigg),
\end{align*}
for some ${C}>0$.  The claim follows as in Step 1.

\paragraph{Step 3:} Let $\wXX\in \widehat{D}^p([0,1],G^1(\Rset^{d+1}))$. For notational convenience,
let  $\ZZ\in \widehat{C}^p([0,1],G^{1}(\Rset^{d+1})$ denote the time-extended Marcus-transformed path of $\wXX$ and  let $\mu_t\in [0,1]$ be such that $\ZZ_{\mu_t}=\wXX_{t}$ for all $t\in [0,1]$ (see Notation~\ref{3notation: mut}).  By Step 1 and Step 2, $\nabla F(\cdot,\ZZ)$ is a continuous path of finite $p'$-variation, for all $p'> p$ such that 
  $\frac{1}{p'}+\frac{1}{p}>1$. Since $\nabla F$ is a Marcus canonical  path functional, for all $t\in [0,1]$,
  \begin{align*}
      \nabla F(t,\wXX)=\nabla F(\mu_t,\ZZ),
  \end{align*}
  implying by definition of $\mu_t$ that $\nabla F(\cdot,\wXX)$ is a \cadlag path of finite $p'$-variation. 
  
\noindent Next, assume first that $\wXX$ admits only one jump at time $a\in (0,1]$. Suppose that $t<a$. By the arguments in Step 2 applied to $\ZZ$, 
   \begin{align*}
       F(\mu_t,\ZZ)= F(0,\ZZ)+\int_0^{\mu_t}\nabla F(s,\ZZ)d\ZZ_s.
   \end{align*}
Observe that 
$\nabla F\in \Mcal_1^1\subset  \Mcal_1^0 $ implies that $\nabla F$ is an invariant under reparametrization path functional (see Proposition~\ref{3prop: F Marcus canonical then F invariant}). This, combined with the property of the Young integral, yields that the functional
\begin{align*}
    G: [&0,1]\times C^{p}([0,1],G^{1}(\Rset^{d+1}))\rightarrow \Rset \\
    &(t,\YY)\mapsto G(t,\YY):=\int_{0}^t\nabla F(s,\YY)d\YY_s
\end{align*}
is also invariant under reparametrization. Thus, for $t<a$,   $\int_0^{\mu_t}\nabla F(s,\ZZ)d\ZZ_s=\int_0^{t}\nabla F(s^-,\wXX)d\wXX_s$, as the stopped path $\ZZ^{\mu_t}$ is nothing but a time-reparametrization of the continuous path $\wXX^t$.
Since $F(\mu_t,\ZZ)=F(t,\wXX)$, $F(0,\ZZ)=F(0,\XX)$, the claim follows. Next, suppose that $t\geq a$ and observe that
          \begin{align*}
       &F(\mu_t,\ZZ)- F(\mu_a,\ZZ)=\int_{\mu_a}^{\mu_t}\nabla F(s,\ZZ)d\ZZ_s=\int_{a}^t\nabla F(s^-,\wXX)d\wXX_s,\\
        &F(\mu_a,  \ZZ)-F(\mu_{a^-},  \ZZ)=F(a,  \wXX)-F({a^-},\wXX),\\
        &F(\mu_{a^-},  \ZZ)-F(0,  \ZZ)=\int_0^{\mu_{a^-}} \nabla^{} F(s,   \ZZ)d \ZZ_s=\int_{0}^{a} \nabla^{} F(s^-,   \wXX)d  \wXX_{s}-\nabla^{} F(a^-,   \wXX)\Delta \wXX_a.
   \end{align*}
        The claim follows by combining all three terms.
        The argument extends trivially to a \cadlag paths with finitely many jumps. Finally, let $\wXX$ be a \cadlag path with countable many jumps. By the arguments in Step 2 applied to $\ZZ$, for $\varepsilon>0$, there  exists  $\eta>0$ and a partition ${\pi}_{[0,\mu_t]}$ with $|{\pi}_{[0,\mu_t]}|<\eta $ such that 
         \begin{align*}
             \big\|F(\mu_t,  \ZZ)-F(0, \ZZ)-\sum_{s_{i}\in {\pi}_{[0,\mu_t]}}\nabla F( s_{i},\ZZ)\ZZ_{ s_{i}, s_{{i+1}}}\big \|<\frac{\varepsilon}{2}.
        \end{align*}
           Since $\wXX$ has finite $p$-variation, for $p\in [1,2)$, we can find a finite set $B(\varepsilon,t)\subset[0,t]$ of jump times of  $\wXX$ such that $ \sum_{s\leq t, \ s \notin B(\varepsilon,t)} \|\Delta \wXX_s\|^2<\frac{\varepsilon}{2}.$ 
         Without loss of generality, we may assume that if $t_j\in B(\varepsilon,t)$, then $\mu_{t_j^-},\mu_{t_j}\in {\pi}_{[0,\mu_t]}$. Thus, repeating earlier arguments, there exists a partition ${\pi}_{[0,t]}$  such that
           \begin{align}\label{3eq: ito cadlag p12 sum big jumps}
             \big \|&F(t,  \wXX)-F(0, \wXX)-\sum_{{s_i}\in {\pi}_{[0,t]}}\nabla F({s_i},{\wXX}){\wXX}_{{s_i},{s_{i+1}}}\\
             & \quad \quad -\sum_{ s\in B(\varepsilon,t)} F(s,\wXX)-F(s^-,\wXX)-\nabla ^{}F(s^-,\wXX) \Delta \wXX_s\big\|<\frac{\varepsilon}{2}.\nonumber 
        \end{align}
        Next, let $(\u_n)_{n\in \N}\subset T(\Rset^{d+1})$ be such that $ \lim_{n\rightarrow\infty }\sup_{t\in [0,1]}\sum_{j=0}^2\|\nabla ^jF(t,\widehat \XX)-\langle \u_n^{(j)},\widehat \X_t\rangle\|=0,$ from Theorem~\ref{3th: Nachbin for path functionals}. By Remark~\ref{3rem: after proof Nachbin}\ref{3rem: after proof Nachbin i}, for all $n\in \Nset$,
           \begin{align*}
            \sum_{s\notin B(\varepsilon,t)} &\|\langle \u_n,\wX_s\rangle -\langle \u_n,\wX_{s^-}\rangle-\langle \u^{(1)}_n,\wX_{s^-}\rangle \Delta \wXX_s\|\\
          &  \leq
             \sup_{n\in\Nset} \sup_{s\notin B(\varepsilon,t)}\sup_{\theta \in [0,1]}\|\langle\u_n^{(2)},\widehat \X_{s^-} \otimes \exp( \theta\log^{(1)}(\Delta  \widehat\XX_{s}))\rangle\|
             \sum_{s\notin B(\varepsilon,t)} \|\Delta \wXX_s\|^2<\infty,
        \end{align*}
       an application of the dominated convergence theorem yields that
        \begin{align*}
            \sum_{s\notin B(\varepsilon,t)} \|F(s,\wXX)-F(s^-,\wXX)-\nabla ^{}F(s^-,\wXX) \Delta \wXX_s\|<\infty.
        \end{align*}
       Thus, repeating the above argument, we can take $\varepsilon>0$  such that  
         \begin{align}\label{3eq: RRS ito p 12}
          F(t,  \wXX)-F(0, \XX)= &\lim_{(RRS)|\pi_{[0,t]}|\rightarrow 0}\sum_{{s_i}\in \bar{\pi}_{[0,t]}}\nabla F({s_i},{\wXX}){\wXX}_{{s_i},{s_{i+1}}}\\
              &-\sum_{0<s\leq t} F(s,\wXX)-F(s^-,\wXX)-\nabla ^{}F(s^-,\wXX) \Delta \wXX_s. \nonumber 
             \end{align}
             Since $\wXX$ is a \cadlag path, by Proposition 2.4 in~\cite{FZ:17}, the convergence in equation~\eqref{3eq: RRS ito p 12} holds in Mesh Riemann-Stieltracking-jumpses sense (see Definition~\ref{3def: MRS RRS}) and the claim follows.

\end{proof}

\begin{remark}
 \label{3rem: Ito formula p12} An inspection of the above proof shows that we can replace the hypothesis of $F$ being $[p,2)$-var continuous with the assumption that for every $(\XX^n, \XX)_{n\in\Nset}\subset C^p([0,1],G^{[p]}(\Rset^{d}))$,
 $$\lim_{n\rightarrow\infty}\sup_{t\in [0,1]}\dcc(\XX^n_t,\XX_t)=0 $$ implies $$\lim_{n\rightarrow\infty}\sup_{t\in[0,1]}\dcc(F(t,\XX^n),F(t,\XX))=0,$$

\end{remark}

\subsection{The case \texorpdfstring{$p\in [2,3)$}{p in [2,3)}}\label{sec: into formula p23}
We present the functional \Ito-formula for maps of \cadlag $p$-rough paths, for $p\in [2,3)$.

   \begin{theorem}\label{3th: Ito formula p 23}
     Let $p\in[2,3)$ and $F:[0,1]\times D^{p}([0,1],G^{2}(\Rset^{d+1}))\rightarrow \Rset$ be a $C^3$-non-anticipative Marcus canonical  path functional such that $F,\nabla F, \nabla^2 F$, and $   R^{\nabla F,\nabla^2F}$ are $[p,3)$\text{-}var continuous. Then, for every Marcus-like $\widehat \XX\in \widehat{D}^{p}([0,1],G^{2}(\Rset^{d+1})) $, it holds that $ \big( \nabla F(\cdot ,\wXX), \nabla^{2}F(\cdot ,\wXX)\big)\in \mathcal{V}^{p',r}_{\wXX}$, for all $p'>p$ such that $\frac{1}{p'}+\frac{2}{p}>1$ and $r\geq 1$ given by $ \frac{1}{r}=\frac{1}{p}+\frac{1}{p'},$  and for all
$t\in [0,1]$, 
    \begin{align}\label{3eq: ito formula p 23}
        F(t, \wXX)-&F(0, \widehat \XX)=\int_0^t \nabla F(s^-, \widehat \XX)d\widehat \XX_s\\
       & +\sum_{0<s\leq t}F(s,\wXX)-F(s^-,\wXX)-\nabla F(s^-,\wXX) \Delta \widehat{X}_s-\nabla ^{2}F(s^-,\wXX) \Delta \wX^{(2)}_s.\nonumber 
    \end{align}
    The integral in~\eqref{3eq: ito formula p 23} is  a rough integral and the summation term is well defined as an absolutely summable series.
\end{theorem}
The proof of Theorem~\ref{3th: Ito formula p 23} will make use of the following lemma, where algebraic properties of the signature of piecewise linear paths are exploited to infer some key analytical properties. Its proof is given in Appendix \ref{proof lemma ito p23}.

\begin{lemma}\label{3lemma: ito formula uniformly bounded p 23}
    Let $p\in[2,3)$, $F:[0,1]\times D^{p}([0,1],G^{2}(\Rset^{d+1}))\rightarrow \Rset$ be such that $F\in  \mathcal{M}^3_{2}$ and $\widehat \XX\in  \widehat C^{1}([0,1],G^2(\Rset^{d+1}))$ be the truncated signature at level $2$ of a time-extended piecewise linear path. Assume that there exists $(\u_n)_{n\in\Nset}\subset T(\Rset^{d+1})$ such that 
\begin{align}\label{3eq: approx F and derivatives p 23}
     \lim_{n\rightarrow\infty} \sup_{t\in [0,1]}\sum_{j=0}^3\|\nabla ^jF^{\u_n}(t,\widehat \XX)-\nabla ^jF(t,\widehat \XX)\|=0.
\end{align}
Then, 
\begin{enumerate}
    \item \label{3item lemma 2 i} $\sup_{n\in \Nset}\bigg(\| \nabla ^2F^{\u_n}(\cdot,\widehat \XX)\|_{1\text{-}var}+\|R^{\u_n^{(1)},\u_n^{(2)}}((\cdot,\cdot),\wXX)\|_{\frac{1}{2}\text{-}var}\bigg)<\infty$;
    \item \label{3item lemma 2 ii} for all $p'>1$, $r>\frac{1}{2}$, 
   \begin{align*}
       &\lim_{n\rightarrow\infty} d_{p'}\big(\nabla ^2F^{\u_n}(\cdot,\widehat \XX),\nabla^{2}F(\cdot,\wXX )\big)=0, \\
       &\lim_{n\rightarrow\infty}d_{r}\big(R^{\u_n^{(1)},\u_n^{(2)}}((\cdot,\cdot),\wXX),R^{\nabla F,\nabla^2 F}((\cdot,\cdot),\wXX)\big)=0.
   \end{align*}
\end{enumerate}
    \end{lemma}

\begin{proof}[Proof of Theorem~\ref{3th: Ito formula p 23}]The structure of the proof is very similar to the one of Theorem~\ref{3th: Ito formula p 12}. Therefore, we emphasize only the main differences.

\paragraph{Step 1:}  Let $\widehat \XX\in  \widehat C^{1}([0,1],G^2(\Rset^{d+1}))$ be the truncated signature at level $2$ of a time-extended piecewise linear path. Since $F$ is a $C^3$, by  Theorem~\ref{3th: Nachbin for path functionals}, there exists a sequence $(\u_n)_{n\in\Nset}\subset T(\Rset^{d+1})$ such that the convergence in equation \eqref{3eq: approx F and derivatives p 23} holds true.
By Lemma~\ref{3lemma: ito formula uniformly bounded p 23}~\ref{3item lemma 2 i} and Lemma 5.12 in~\cite{FV:10},
\begin{align}\label{3eq: ito F^2 and R good control}
     \|\nabla^{2}F(\cdot ,\wXX)\|_{1\text{-}var}+\|R^{\nabla F,\nabla^2 F}((\cdot,\cdot),\wXX)\|_{\frac{1}{2}\text{-}var}<\infty,
\end{align}
implying that $ \big( \nabla F(\cdot ,\wXX), \nabla^{2}F(\cdot ,\wXX)\big)\in \mathcal{V}^{p',r}_{\wXX}$, for all $p'>p$ such that $\frac{1}{p'}+\frac{2}{p}>1$. 

Fix such a $p'$ and notice that $r>\frac{p}{2}>\frac{1}{2}$. By Lemma~\ref{3lemma: ito formula uniformly bounded p 23}~\ref{3item lemma 2 ii},
 \begin{align*}
       \lim_{n\rightarrow\infty} d_{p'}\big(\nabla^{2}F(\cdot,\wXX ),\nabla ^2F^{\u_n}(\cdot,\widehat \XX)\big)=0, \qquad 
       \lim_{n\rightarrow\infty}d_{r}\big(R^{\nabla F,\nabla^2 F}((\cdot,\cdot),\wXX),R^{\u_n^{(1)},\u_n^{(2)}}((\cdot,\cdot),\wXX)\big)=0. \nonumber
   \end{align*}
Moreover, an application of Theorem 8.10 in~\cite{FV:10} and  Proposition 2.7 in~\cite{ALP:23} yields that for all $n\in \Nset$, 
\begin{align*}
    d_p\bigg(\int_0^\cdot &\langle \u^{(1)}_n,\wX_{s}\rangle d\wXX_{s}, \int_0^\cdot \nabla F(s,\wXX) d\wXX_{s}\bigg) \\
    \leq & C\bigg(\dcc\big(\langle \u_n^{(1)},\wX_0\rangle,\nabla F(0,\wXX)\big)+\dcc\big(\langle \u_n^{(2)},\wX_0\rangle,\nabla^2 F(0,\wXX)\big)\bigg.  \nonumber \\
    &  \ \bigg. \quad +d_{p'}\big(\langle \u_n^{(2)},\wX\rangle,\nabla^{2}F(\cdot ,\wXX)\big)+d_{r}\big(R^{\nabla F,\nabla^2 F}((\cdot,\cdot),\wXX),R^{\u_n^{(1)},\u_n^{(2)}}((\cdot,\cdot),\wXX)\big)\bigg)\nonumber
\end{align*}
for some $C>0$.
Since by equation~\eqref{eq: u,signature p23} for every $n\in \Nset$ and $t\in [0,1]$, it holds that $\langle \u_n,\widehat{\mathbb{X}}_{t}\rangle = \langle \u_n,\widehat{\mathbb{X}}_{0}\rangle+  \int_0^t \langle \u^{(1)}_n,\widehat{\mathbb{X}}_{s}\rangle  d\widehat{\XX}_{s},$ and the claim follows.

\paragraph{Step 2:} Fix $\widehat \XX \in  \widehat C^{p}([0,1],G^{2}(\Rset^{d+1})) $. By Theorem 8.12 in~\cite{FV:10}, there exists a sequence  of piecewise linear time-extended paths such that their truncated signature at level $2$, denoted by $(\widehat \XX^M)_M$, satisfy
    \begin{align}\label{3approximation w.g.p 23}
        \lim_{M\rightarrow\infty }\sup_{t\in [0,1]}\dcc(\widehat \XX^M_t,\wXX_t)=0, \qquad \text{and}\qquad    \sup_{M\in \Nset}\|\widehat \XX^M\|_{p\text{-}var}\leq C\|\widehat \XX\|_{p\text{-}var},
    \end{align}
    for some $C>0$. By Step 1, for every fixed $M$ and $t\in [0,1]$, 
        \begin{align*}
        F(t, \widehat \XX^M)-F(0, \widehat \XX^M)=\int_0^t \nabla F(s, \widehat \XX^M)d\widehat \XX^M_s.
    \end{align*}

 By conditions~\eqref{3approximation w.g.p 23} and interpolation (see Lemma 5.12 and Lemma 8.16 in~\cite{FV:10}), for all $q>p$, $  \lim_{M\rightarrow\infty }d_{q}(\widehat \XX^M, \widehat \XX)=0.$ Fix $p'$ and $r$ as in Step 1 and notice that since for all $M$, $ \big(  \nabla  F(\cdot,\widehat \XX^M), \nabla ^{2} F(\cdot,\widehat \XX^M)\big)\in \mathcal{V}^{p',r}_{\wXX},$ and $\nabla F, \nabla^2 F$ and $R^{\nabla F,\nabla ^2F}$ are by assumption $[p,3)$\text{-}var continuous, $ \big(  \nabla  F(\cdot,\widehat \XX), \nabla ^{2} F(\cdot,\widehat \XX)\big)\in \mathcal{V}^{p',r}_{\wXX}$ too.

Next, fix $p''>p'>p $ such that $\frac{2}{p'}+\frac{1}{p''}>1$. Observe that for all $M\in \Nset$, $\wXX^M,\wXX\in \widehat{D}^{p'}([0,1], G^2(\Rset^{d+1}))$ and by~\eqref{3eq: ito F^2 and R good control},  $ \big(  \nabla  F(\cdot,\widehat \XX^M), \nabla ^{2} F(\cdot,\widehat \XX^M)\big)\in \mathcal{V}^{p'',r}_{\wXX^M}$. Therefore, since  $\nabla F, \nabla^2 F$ and $R^{\nabla F,\nabla ^2F}$ are $[p,3)$\text{-}var continuous,  $ \big(  \nabla  F(\cdot,\widehat \XX), \nabla ^{2} F(\cdot,\widehat \XX)\big)\in \mathcal{V}^{p'',r}_{\wXX}.$ 
 Finally, applying Theorem 8.10 of ~\cite{FV:10} and Proposition 2.7 of ~\cite{ALP:23} yields that there exists a constant $C>0$, such that for all $M\in \Nset$,
 \begin{align}\label{3eq: proof  step 2 ito p23 }
    d_{p'}\bigg(\int_0^\cdot & \nabla F(s,\wXX^M) d\wXX^M_{s}, \int_0^\cdot \nabla F(s,\wXX) d\wXX_{s}\bigg) \\
    \leq & C\bigg(d_{p'}\big(\wXX^M,\wXX \big)+\dcc\big(\nabla F(0,\wXX^M),\nabla F(0,\wXX)\big)+\dcc\big(\nabla^2 F(0,\wXX^M),\nabla^2 F(0,\wXX)\big)\bigg.  \nonumber \\
    &  \ \bigg. \quad +d_{p''}\big(\nabla^{2}F(\cdot ,\wXX^M),\nabla^{2}F(\cdot ,\wXX)\big)+d_{r}\big(R^{\nabla F,\nabla^2 F}((\cdot,\cdot),\wXX^M),R^{\nabla F,\nabla^2 F}((\cdot,\cdot),\wXX)\big)\bigg).\nonumber
\end{align}
Since $F$ is $[p,3)$ var continuous, the claim follows.
\paragraph{Step 3:} 
        
Fix $\wXX\in \widehat{D}^p([0,1],G^1(\Rset^{d+1}))$, with $\wXX$ Marcus-like and set $\widehat{X}:=\pi_1(\wXX)$. For notational convenience, let us denote by $\ZZ\in \widehat{C}^p([0,1],G^{2}(\Rset^{d+1}))$ the time-extended Marcus-transformed path of $\wXX$ and  let $\mu_t\in [0,1]$ be such that $\ZZ_{\mu_t}=\wXX_{t}$ for all $t\in [0,1]$ (see Notation~\ref{3notation: mut}). It holds that $ \big( \nabla F(\cdot ,\wXX), \nabla^{2}F(\cdot ,\wXX)\big)\in \mathcal{V}^{p',r}_{\wXX}$.
Next, assume first that $\wXX$ admits only one jump at time $a\in (0,1]$. If $t<a$, the claim follows as in the Young case. If $t\geq a$, it suffices to notice that  for level $2$ rough integration, 
          \begin{align*}
        &F(\mu_{a^-},  \ZZ)-F(0,  \ZZ)=\int_{0}^{a} \nabla^{} F(s^-,   \wXX)d  \wXX_{s}-\nabla F(a^-,\wXX) \Delta \widehat{X}_a-\nabla ^{2}F(s^-,\wXX) \Delta \wX^{(2)}_a.
   \end{align*}
        The argument extends trivially to \cadlag paths with finitely many jumps. Finally, for $\wXX$ a \cadlag path with countable many jumps, since $\wXX$ has finite $p$-variation, for $p\in [2,3)$, there exists a finite set $B(\varepsilon,t)\subset[0,t]$ of jump times of  $\wXX$ such that $ \sum_{s\leq t, \ s \notin B(\varepsilon,t)} \|\Delta \widehat{X}_s
        \|^3<\frac{\varepsilon}{2},$ for $\varepsilon>0$. This very last step of the proof follows from Proposition 2.6 in~\cite{FZ:17}.

\end{proof}
\begin{remark}\phantomsection\label{3rem: Ito formula p23}
\begin{enumerate}
    \item \label{3rem: Ito formula p23 item i} An inspection of the above proof shows that we can replace the hypothesis of $F$, $\nabla F$ being $[p,3)$-var continuous with the assumption that for every sequence $(\XX^n, \XX)_{n\in\Nset}\subset C^p([0,1],G^{[p]}(\Rset^{d}))$,
 $$\lim_{n\rightarrow\infty}\sup_{t\in [0,1]}\dcc(\XX^n_t,\XX_t)=0 $$ implies $$\lim_{n\rightarrow\infty}\sup_{t\in[0,1]}\dcc(F(t,\XX^n),F(t,\XX))=0 ,\qquad \lim_{n\rightarrow\infty}\sup_{t\in[0,1]}\dcc(\nabla F(t,\XX^n),\nabla F(t,\XX))=0.$$ 
 \item \label{3rem: Ito formula p23 item ii} Instead of assuming $R^{\nabla^1F,\nabla^2 F}$ to be $[p,3)$-var continuous, one can also suppose that $F,\nabla F, \nabla^2 F$, and $ \nabla^3F$ are $[p,3)$\text{-}var continuous. Indeed,  for every (signature at level $2$ of) a piecewise linear time-extended path $\wXX^M$ as in Step 2, by Therem~\ref{3th: Nachbin for path functionals} and equation~\eqref{eq: reminder linear}, we get
 \begin{align*}
     \|R^{\nabla^1F,\nabla^2 F}((s,t),\wXX^M)\|&=\|\nabla^1F(t,\wXX^M)-\nabla^1F(s,\wXX^M)-\nabla^2F(s,\wXX^M)\widehat{X}^M_{s,t}\|\\
     &\leq C \sup_{s\in [0,1]}\|\nabla^3F(u,\wXX^M)\|\|\widehat{X}^M_{s,t}\|^2,
 \end{align*}
 for $\widehat{X}_{s,t}:=\pi_1(\wXX_{s,t})$, $(s,t)\in \Delta_1$ and $C>0.$ Therefore, if $\nabla^3F$ is $[p,3)$-var continuous, the sequence $ \|R^{\nabla^1F,\nabla^2 F}((\cdot,\cdot),\wXX^M)\|$ is uniformly bounded in $p/2$-variation. By interpolation  for all $p'>p$,
 \begin{align*}
     \lim_{M\rightarrow\infty }d_{\frac{p'}{2}}\big(R^{\nabla F,\nabla ^2F}((\cdot,\cdot),\wXX^M),R^{\nabla F,\nabla ^2F}((\cdot,\cdot),\wXX)\big)=0.
 \end{align*}
   
\item We highlight that the consideration of tracking-jumps-extended paths is necessary not only for deriving the approximation result in Theorem~\ref{3th: Nachbin for path functionals}, but also for expressing the dependence on a non-anticipative path functional on the parameter $t$. Consider, for instance, a path functional 
 $F$ such that for all $(t,\XX)\in [0,1]\times {D}^p([0,1],G^{2}(\Rset^{d+1}))$,
\begin{align*}
    F(t,\XX)=\sin(\langle \e_0,\XX_t\rangle).
\end{align*}
Then, for a tracking-jumps-extended path $\wXX$, $ F(t,\wXX)=\sin(Z(t))$, computing the vertical derivative as described in \eqref{eq: vertical derivativ} and \eqref{3eq: value vertical derivatives higher order}, we get 
$U^0F(t,\wXX)=\cos(Z(t))$, $U^0U^0F(t,\wXX)=-\sin(Z(t))$, and all the other (first and second order) vertical derivatives are equal to $0$. Since the hypothesis of Theorem~\ref{3th: Ito formula p 23} are satisfied, we get that
\begin{align*}
    \sin(Z_t)-\sin(Z_0)=\int_0^t\cos(Z_{s^-})dZ_s
     +\sum_{0<s\leq t} \sin (Z_{s})-\sin(Z_{s^-})-\cos(Z_{s^-})\Delta Z_s,
\end{align*}
which coincides with the change of variable formula for paths of finite variation. In particular, if $\wXX$ is continuous, the path functional explicitly reads as $F(t,\XX)=\sin(t)$, and the consideration of time-extended paths is necessary in order to embed the functional $F$ in the present framework.

    \item Theorems~\ref{3th: Ito formula p 12} and~\ref{3th: Ito formula p 23} provide an \Ito-formula for $C^2$ and $C^3$ non-anticipative Marcus canonical  path functional, respectively (see Definition~\ref{3def: C^K path functionals}). Therefore, our framework includes all the path functionals for which for all $i=0,\dots,2$ (or $i=0,\dots,3$) the maps $ [0,1]\ni t\mapsto \nabla^{i}F(t,\wXX)$ (or equivalently the map $t\mapsto \tilde{F}(\wXX^t)$ introduced in Remark \ref{remark F tilde F}\ref{remark F tilde F ii}) are continuous if $\wXX\in \widehat{C}^p([0,1],G^{[p]}(\Rset^{d+1}))$, and \cadlag if $\wXX\in \widehat{D}^p([0,1],G^{[p]}(\Rset^{d+1}))$ (see Definition~\ref{3def: C^K path functionals}). This in particular implies that functionals of the form $F(t,\XX):=\XX_{s\wedge t}$, for some fixed $s\in (0,1]$ are not included in our setup. A similar argument applies to the functional $F(t,\XX):=\XX_{t^-}$, as the map $[0,1]\ni t\mapsto F(t,\XX)$ is \caglad for $\XX$ \cadlag.  Notice furthermore that such functionals are not Marcus canonical  as condition \ref{item def Mcpf i} in Definition \ref{3def: Marcus canonical path functional} is not satisfied. Similarly, condition  \ref{item def Mcpf i} is satisfied for neither the functionals
    $F(t,\XX):= \XX^2_{\XX^1_t}1_{\{\XX^1_t\in[0,1]\}} $ nor   the delayed functional, $F(t,\XX):=\XX_{t-\delta}$, for some $\delta>0$. This is also the case in the setting of ~\cite{CRF:13} (see  Example 6 therein), however for different reasons.

\item\label{remark reduced} If in particular for each $\XX\in  D^p([0,1],G^2(\Rset^{d+1}))$, the functional $F$ and its derivatives depend only on $\pi_1(\XX)$, and the second order vertical derivative $\nabla ^2F(\cdot,\XX)$ is a path with values in the subspace of symmetric matrices, a rough functional \Ito-formula can be derived by neglecting the information provided by the area of $\XX$, i.e. $\text{Anti}(\X^{(2)})$, and considering instead a rough integral with respect to the \textit{canonical reduced rough path} (see Definition 5.3 in \cite{FH:14})
\begin{align}\label{eq: reduced rp canonical}
    \XX^{\mathcal{R}}:=(X,\mathbb{S}^X)
\end{align}
where $X:=\pi_1(\XX)$ and $\mathbb{S}^X:\Delta_1\rightarrow\text{Sym}((\Rset^{d+1})^{\otimes 2})$ is given by $\mathbb{S}^X_{s,t}:=\frac{1}{2}X_{s,t}^{\otimes 2}$ for each $(s,t)\in \Delta_1$. The proof technique is the same as the one in Theorem \ref{3th: Ito formula p 23}. Moreover, in this specific case, the \Ito-formula can be deduced by considering the  variation topology on the space $C^p([0,1],\Rset^d)$ instead of the one on $C^p([0,1],G^2(\Rset^d))$ as described in Definition \ref{3def: var continuous funcitonals}.

\item Following \cite{P:24}, an \Ito-formula for non-anticipative functionals of paths with values in $\Rset^d$ is derived subsequently  in \cite{B:24} using rough integrals with respect to \textit{continuous reduced rough paths} of $\Rset^d$-valued paths of arbitrary regularity. This alternative approach adopts the notion of functional vertical derivatives proposed in \cite{D:09}, and consequently the higher-order functional vertical derivatives always take value  in the space of symmetric tensors,  as highlighted in Remark \ref{rem: non commutativity derivate}. The corresponding   rough integrals are  such functionals with respect to the canonical reduced rough paths derived from the powers of the increments of $\pi_1(\XX)$. We remark that due to the symmetry of the vertical derivatives within this framework, a Taylor expansion based on the signature of non-anticipative path functionals cannot be derived (see also Remark \ref{remark: taylor p23}\ref{rem: valentin dupire}).

\end{enumerate}

\end{remark}

\subsection{Connections to the literature}\label{sec: connections litarature}

\paragraph{The \Ito-formula for rough paths.}  We will see that the functional \Ito-formula in Theorem~\ref{3th: Ito formula p 23} matches the existing  \Ito-formulas for rough paths in the literature. Throughout, we fix $p\in [2,3)$.
\begin{itemize}
    \item Let $f\in C^3(\Rset^{d+1})$ and consider the path functional given by
    \begin{align*}
        [0,1]\times \widehat{D}^{p}([0,1],G^{2}(\Rset^{d+1}))\ni (t,\wXX)\mapsto F(t,\wXX):=f(\widehat{X}_t),
    \end{align*}
    where we set $\widehat{X}:=\pi_1(\wXX)$.
    Observe that $F$ is a $C^3$-non-anticipative Marcus canonical  path functional such that for all $i=1,2,3$, $(t,\wXX)\in  [0,1]\times \widehat{D}^{p}([0,1],G^{2}(\Rset^{d+1}))$,
    \begin{align*}
        \nabla^i F(t,\wXX)=\nabla^i f(\widehat{X}_t),
    \end{align*}
where  $\nabla^i f$ denotes the $(\Rset^{d+1})^{\otimes i}$-valued maps given by the $i$-th (standard) partial derivatives of the map $f$. Moreover, from  properties of regular functions on $\Rset^{d+1}$, $F,\nabla F, \nabla^2 F$, and $   R^{\nabla F,\nabla^2F}$ are 
$[p,3)$\text{-}var continuous. Therefore, applying Theorem \ref{3th: Ito formula p 23}
yields that for all $ (t,\wXX)$ with $\wXX$-Marcus-like, 
\begin{align}\label{3eq: ito formula f non path dependent}
        f(\widehat{X}_t)-&f(\widehat{X}_0)=\int_0^t \nabla f(\widehat{X}_{s^-})d \wXX_s\\
       & +\sum_{0<s\leq t}f(\widehat{X}_s)-f(\widehat{X}_{s^-})-\nabla f(\widehat{X}_{s^-}) \Delta \widehat{X}_s-\nabla ^{2}f(\widehat{X}_{s^-}) \frac{1}{2}\Delta \widehat{X}_s^{\otimes2 }.\nonumber 
    \end{align}
      Notice that the considerations in Remark \ref{3rem: Ito formula p23}\ref{remark reduced} apply here. 
    The formula in equation~\eqref{3eq: ito formula f non path dependent} has been derived in Theorem 2.12 of ~\cite{FZ:17} by  
    exploiting the Taylor expansion of the regular function $f$ on $\Rset^{d+1}$, a completely different techniques from the one used in Theorem~\ref{3th: Ito formula p 23}.

    \item Let $f\in C^3(\Rset)$ and consider the path functional given by
    \begin{align}\label{eq: f compones sig}
        [0,1]\times \widehat{D}^{p}([0,1],G^{2}(\Rset^{d+1}))\ni (t,\wXX)\mapsto F(t,\wXX):=f(\langle \u,\wX_t\rangle),
        \end{align}
        for some $\u\in T(\Rset^{d+1})$. 
        Applying the rules of derivation for compound functions yields that $F$ is a $C^3$-non-anticipative Marcus canonical  path functional such that 
    \begin{align*}
        &\nabla F(t,\wXX)=f'(\langle \u,\wX_t\rangle )\langle \u^{(1)},\wX_t\rangle,\\
           &\nabla^2 F(t,\wXX)=f''(\langle \u,\wX_t\rangle )\langle \u^{(1)},\wX_t\rangle
           \langle \u^{(1)},\wX_t\rangle^\top+f'(\langle \u,\wX_t\rangle )\langle \u^{(2)},\wX_t\rangle ,\\
            &\nabla^3 F(t,\wXX)=f'''(\langle \u,\wX_t\rangle )\langle \u^{(1)},\wX_t\rangle\langle \u^{(1)},\wX_t\rangle^\top \langle \u^{(1)},\wX_t\rangle^\top\\
            &\qquad \qquad \qquad +3f''(\langle \u,\wX_t\rangle )\langle \u^{(2)},\wX_t\rangle\langle \u^{(1)},\wX_t\rangle^\top 
             +f'(\langle \u,\wX_t\rangle )\langle \u^{(3)},\wX_t\rangle,
    \end{align*}
    for all $(t,\wXX)\in  [0,1]\times \widehat{D}^{p}([0,1],G^{2}(\Rset^{d+1}))$, with $f',f'',f'''$ denoting the first, the second, and the third derivatives of $f$, respectively.  A further application of the properties of regular functions on $\Rset$ yields that $F,\nabla F, \nabla^2 F$ and $   R^{\nabla F,\nabla^2F}$ are $[p,3)$\text{-}var continuous. Therefore the assertion of Theorem~\ref{3th: Ito formula p 23} holds for all $ (t,\wXX)\in [0,1]\times \widehat{D}^{p}([0,1],G^{2}(\Rset^{d+1}))$, with $\wXX$-Marcus-like.   In particular, if  $\wXX\in (t,\wXX)\in [0,1]\times \widehat{C}^{p}([0,1],G^{2}(\Rset^{d+1}))$, the functional \Ito-formula in equation~\eqref{3eq: ito formula p 23} coincides with the \Ito-formula for controlled rough paths stated in Theorem 7.7 of ~\cite{FH:14}. 
Observe that the above reasoning can be easily generalized to path functionals of the form 
    \begin{align*}
         [0,1]\times \widehat{D}^{p}([0,1],G^{2}(\Rset^{d+1}))\ni (t,\wXX)\mapsto F(t,\wXX):=g(\langle \u_1,\wX_t\rangle,\dots,\langle \u_m,\wX_t\rangle),
    \end{align*}
 for some $g\in C^3(\Rset^m)$,  $\u_1,\dots\u_m\in T(\Rset^{d+1})$. 
    \end{itemize}

\paragraph{Föllmer, (RIE) and stochastic integration theories.} We elaborate on the connections of the results in Section~\ref{3sec: functional Ito formula} with Föllmer (\cite{F:81}), (RIE) (~\cite{PP:16},~\cite{ALP:23}), and stochastic integration (\cite{JS:87}) theories, whose main concepts have been summarized in Appendix~\ref{3appendix: integration theories}. The proof of the following corollary is given in Appendix \ref{proof 3coro: ito formula follmer}.
\begin{corollary}[Föllmer]\label{3coro: ito formula follmer}
     Let $p\in[2,3)$, $F:[0,1]\times D^{p}([0,1],G^{2}(\Rset^{d+1}))\rightarrow \Rset$ be a $C^3$-non anticipative Marcus canonical  path functional such that $F,\nabla F, \nabla^2 F$ and $   R^{\nabla F,\nabla^2F}$ are $[p,3)$\text{-}var continuous, and  $(\pi^n_{[0,1]})_{n\in\Nset}$ a sequence of partitions with vanishing mesh size. Let $\widehat{X}\in \widehat{D}^p([0,1],\Rset^{d+1})$ be a path with finite quadratic variation in the sense of Föllmer along  $(\pi^n_{[0,1]})_{n\in\Nset}$,  $[\widehat X,\widehat X]^c$ its continuous quadratic variation, and  $\widehat \XX\in \widehat{D}^{p}([0,1],G^{2}(\Rset^{d+1})) $ a Marcus-like weakly geometric rough path such that $\pi_1(\wXX)=\widehat{X}$.
     Assume that for all $t\in [0,1]$, 
     $\nabla^2 F(t,\wXX)=\text{Sym}(\nabla^2 F(t,\wXX))$.
     Then, the following limit exists,
     \begin{align}\label{3eq: coro follmer limit integral}
        \int_0^t\nabla F(s^-,\wXX)d\widehat{X}_s:=\lim_{n\rightarrow\infty }\sum_{s^n_i\in \pi^n_{[0,1]}}\nabla F(s^n_i,\wXX) \widehat{X}_{s^n_i,s^n_{i+1}},
    \end{align} 
and
    \begin{align}\label{3eq: ito formula follmer}
        F(t, \wXX)-F(0, \widehat \XX)=&\int_0^t\nabla F(s^-,\wXX)d\widehat{X}_s+\frac{1}{2}\int_0^t\nabla^2F(s^-,\wXX)d[\widehat{X},\widehat{X}]^c_s\\
       & +\sum_{0<s\leq t}F(s,\wXX)-F(s^-,\wXX)-\nabla F(s^-,\wXX) \Delta \widehat{X}_s.\nonumber 
    \end{align}
     The integral in~\eqref{3eq: ito formula follmer} with respect to the path $[\widehat{X},\widehat{X}]^c$ is understood as a Young integral, and the summation term is well defined as an absolutely summable series.
\end{corollary}

Next,  consider the special case of \cadlag rough paths over some $\widehat{X}\in \widehat{D}([0,1],\Rset^{d+1})$ which satisfies (RIE) with respect to some $p\in (2,3)$ and some sequence of nested partitions $(\pi^n_{[0,1]})_{n\in\Nset}$ (Property~\ref{3property RIE}). Recall from Proposition~\ref{3prop: RIE integal properties} that any path that satisfies (RIE) along a sequence of partition also has quadratic variation in the sense of Föllmer along the same partition. 
The proof of the next corollary is given in Appendix \ref{proof 3coro: ito RIE}.  Furthermore, we refer to Remark \ref{rem: ito=dupire}\ref{rem: ito=dupire i} for a comparison on its conditions with those in Corollary \ref{3coro: ito formula follmer}.
\begin{corollary}[(RIE) property]\label{3coro: ito RIE }

     Let $p\in(2,3)$ and $F:[0,1]\times D^{p}([0,1],G^{2}(\Rset^{d+1}))\rightarrow \Rset$ be a $C^3$-non-anticipative Marcus canonical  path functional such that $F,\nabla F, \nabla^2 F$, and $   R^{\nabla F,\nabla^2F}$ are $[p,3)$\text{-}var continuous, and  $(\pi^n_{[0,1]})_{n\in\Nset}$ be a sequence of nested partitions with vanishing mesh size. Let $\widehat{X}\in \widehat{D}([0,1],\Rset^{d+1})$ be a tracking-jumps-extended path which satisfies (RIE) with respect to $p$ and $(\pi^n_{[0,1]})_{n\in\Nset}$,  $[\widehat X,\widehat X]^c$ be its continuous quadratic variation, and $\widehat \XX\in \widehat{D}^{p}([0,1],G^{2}(\Rset^{d+1})) $ the rough path specified in Proposition \ref{3prop: RIE integal properties}\ref{i RIE} such that $\pi_1(\wXX)=\widehat{X}$. Then, for all $t\in [0,1]$,
    \begin{align}\label{3eq: ito formula RIE}
        F(t, \wXX)-F(0, \widehat \XX)=&\int_0^t\nabla F(s^-,\wXX)d\widehat{X}_s+\frac{1}{2}\int_0^t\nabla^2F(s^-,\wXX)d[\widehat{X},\widehat{X}]^c_s\\
       & +\sum_{0<s\leq t}F(s,\wXX)-F(s^-,\wXX)-\nabla F(s^-,\wXX) \Delta \widehat{X}_s.\nonumber 
    \end{align}
    The first integral of~\eqref{3eq: ito formula RIE} is interpreted as in~\eqref{3eq: coro follmer limit integral}, the second is a Young integral with respect to $[X,X]^c$, and the summation term is well defined as an absolutely summable series.
\end{corollary}

Now, let us  analyze path functionals evaluated at some random rough paths. More precisely,  consider the  Marcus lift of some \cadlag semimartingale,  given by $ \mathbf{X}_t:=(1,X_{t}-X_0,\mathbb{X}^{(2)}_{0,t}),$
for all $t\in[0,1]$, where,
\[
	\mathbb{X}^{(2)}_{0,t}:=\int_0^t X_{0,s^-}\otimes dX_s+\frac{1}{2}[X,X]^c_{t}+\frac{1}{2}\sum_{0< s \leq t}\Delta X_s\otimes\Delta X_s. 
	\]
    Here the integral is an \Ito-integral and $[X,X]^c$ denotes the continuous quadratic variation of $X$ (Proposition 16 in \cite{FS:17}). The proof of the following corollary follows directly from  Theorem~\ref{3th: Ito formula p 23}, Lemma 4.35 in~\cite{CF:19}.

\begin{corollary}[Semimartingales (stochastic)]\label{3coro: ito sotchastic}
 For $p\in(2,3)$, consider  a $C^3$-non-anticipative Marcus canonical  path functional  $F:[0,1]\times D^{p}([0,1],G^{2}(\Rset^{d+1}))\rightarrow \Rset$ such that $F,\nabla F, \nabla^2 F$, and $   R^{\nabla F,\nabla^2F}$ are $[p,3)$\text{-}var continuous. Let $\wXX$  denote the Marcus lift a tracking-jumps-extended \cadlag semimartingale $\widehat{X}$.  Assume that the processes $\nabla F(\cdot,\wXX):=(\nabla F(t,\wXX))_{t\in [0,1]},\nabla^2 F(\cdot,\wXX):=(\nabla^2 F(\cdot,\wXX))_{t\in [0,1]}$ are locally bounded and predictable. Then, for all $t\in [0,1]$, a.s. 
     \begin{align}\label{3eq: ito formula stochastic}
        F(t,\wXX)-F(0,\wXX)=&\int_0^t\nabla F(s^{-},\wXX)d\widehat{X}_s+ \frac{1}{2}\int_0^t\nabla^{2}F(s^-,\wXX)d[\widehat{X},\widehat{X}]^c_s\\
        & +\sum_{0<s\leq t}F(s,\wXX)-F(s^{-},\wXX)-\nabla F(s^-,\wXX)\Delta \widehat{X}_s.\nonumber  \end{align}
 The first integral   of~\eqref{3eq: ito formula stochastic} is an \Ito-integral, the second  is a Young integral with respect to $[\widehat{X},\widehat{X}]^c$, and  the summation term is well defined as an a.s. absolutely summable series.
        \end{corollary}
     We conclude this section with some remarks on the above corollaries.

\begin{remark}\phantomsection\label{rem: ito=dupire}

\begin{enumerate}

\item \label{rem: ito=dupire i} The \Ito-formula~\eqref{3eq: ito formula RIE} coincides with the formula derived in Corollary~\ref{3coro: ito formula follmer}. Observe that the former has been deduced under weaker assumptions on $\nabla^2 F$, but stronger ones on the path $\widehat{X}$. Specifically, in Corollary~\ref{3coro: ito formula follmer}, we assumed that for all $t\in [0,1]$, 
     $\nabla^2 F(t,\wXX)=\text{Sym}(\nabla^2 F(t,\wXX))$ to ensure the convergence of the limit in~\eqref{3eq: coro follmer limit integral}. In contrast, the latter convergence is achieved without any assumptions on the antisymmetric part of 
  $\nabla^2F$ when $\widehat{X}$ satisfies property (RIE). Property (RIE) is indeed a stronger requirement for a path than having finite quadratic variation in the sense of Föllmer (see Proposition~\ref{3prop: RIE integal properties}).

\item \label{item rem: ito=dupire i}

 The formula in~\eqref{3eq: ito formula stochastic} and the one derived in Theorem 31 of \cite{D:09} coincide within their common domain of validity, provided that $\nabla F$ and $\nabla ^2F$, as computed in our framework with respect to the direction $i=1,\dots,d$  and evaluated at continuous paths, are equal with their functional derivative representations. A necessary condition for this to be satisfied is that $F(\cdot,\wXX)=F(\cdot,\pi_1(\wXX))$ for all $\wXX$. Notice that in such a case, for all $t\in [0,1]$,
\begin{align}\label{eq: vertical time=horizontal}
    U^0F(t,\pi_1(\wXX))=\mathcal{D}F(t,\pi_1(\XX)),
\end{align}
where  the LHS is the vertical derivative of $F$ at $(t,\wXX)$ with respect to the time component of $\wXX$ as computed in \eqref{eq: vertical derivativ}, and   the RHS is the so-called \textit{horizontal derivatives} of $F$ at $(t,\pi_1(\XX))$ considered in \cite{D:09} and defined via
\begin{align}\label{eq: horizontal derivative}
        \mathcal{D}F(t,\pi_1(\XX)):=\lim_{h\searrow0}F(t+h,\pi_1(\XX)^t),
    \end{align}
    whenever the limit exists. Equality \eqref{eq: vertical time=horizontal} can be proved by matching the two formulas and noticing that by the consistency between Young and level 2 rough integration (see Section \ref{sec: rough=young}), for all $t\in [0,1]$,
    \begin{align*}
        \int_0^t U^0F(s,\pi_1(\wXX))ds=\int_0^t \mathcal{D}F(s,\pi_1(\XX))ds.
    \end{align*}

    Notice that $\mathcal{D}F$ is a functional of the \textit{non-time extended} path $\XX$ given by removing the time component of $\wXX$ (see Remark \ref{rem: extended paths}\ref{item 4 remark extended path}). Indeed, a key difference between the two approaches is  that  the current framework captures the dependence on the time of the functionals   by time-extending the path $\XX$ and considering functionals 
 of time-extended paths (see e.g., Remark \ref{remark: non marucs like}\ref{remark: non marucs like ii}). In \cite{D:09}  this dependence is instead expressed  via the horizontal derivative. Moreover, by conditions \ref{item def Mcpf i},\ref{item def Mcpf ii} in Definition \ref{3def: Marcus canonical path functional}, for every non-anticipative Marcus canonical  path functional $F$ and time-extended continuous path $\wXX$, the horizontal derivative \eqref{eq: horizontal derivative} is always  $0$. Indeed, for every $\varepsilon>0$ such that $ 1-\varepsilon\geq t$,    
 \begin{align*}
     F(t+\varepsilon,\wXX^t)=& F(t+\varepsilon,(\wXX^t)^{t+\varepsilon})\\=&F(t+\varepsilon,(\wXX^t)^{t+\varepsilon, \rhd})\\
     =&F(t+\varepsilon,\wXX^t_\phi)\\=&F(\phi(t+\varepsilon),\wXX^{t}),
 \end{align*}
 for some time-reparametrization $\phi$ such that $\phi([0,t+\varepsilon])=[0,t]$.

 Finally,  the framework in~\cite{D:09} for deriving a functional \Ito-formula requires continuity of both the functionals and their derivatives with respect to the supremum norm.
As pointed out in the introduction, this excludes several important examples of non-anticipative functionals, such as linear functionals of the signature, which are covered in the present framework. 
More generally, the regularity conditions on the functionals and their derivatives expressed with respect to a stronger topology ($p$-variation versus uniform topology) allow here for the consideration of a larger set of regular functionals.

\end{enumerate}
  
\end{remark}

\section{Functional Taylor expansion}\label{sec: functional taylor} 

Taylor expansion is  fundamental in classical calculus, offering explicit polynomial approximations of smooth functions. This section is to derive a \textit{functional Taylor expansion} of sufficiently regular Marcus canonical path functionals in terms of the signature. The core idea behind its derivation is to iteratively apply the \Ito-formulas  in Theorem \ref{3th: Ito formula p 12} and Theorem \ref{3th: Ito formula p 23}.

\subsection{The case \texorpdfstring{$p\in [1,2)$}{}}

\begin{theorem}\label{th: taylor 12}
   Let $p\in[1,2)$, $K\geq2$  and $F:[0,1]\times D^{p}([0,1],G^{1}(\Rset^{d+1}))\rightarrow \Rset$ be a $C^K$ non-anticipative Marcus canonical  path functional. Assume that   $F,\nabla F, \dots, \nabla ^{K-1}F $ are $[p,2)$\text{-}var continuous. Let $\widehat \XX\in \widehat D^{p}([0,1],G^{1}(\Rset^{d+1})) $, denote by $\wX$ its signature and by $\ZZ$ its time-extended Marcus-transformed path with respect to some pair $(R,\psi_R)$. Then, for all 
$t\in [0,1]$,
   \begin{align}\label{3eq: taylor expansion p12}
    F(t,\widehat \XX)=&\sum_{j=0}^{K-2}\nabla^{j}F (0,\widehat \XX)\widehat \X^{(j)}_{0,t}+ \mathcal{R}^F_{K-1}(t,\widehat \XX),
\end{align}
where, $$\mathcal{R}^F_{K-1}(t,\wXX):=\int_0^{\psi_{R}^{-1}(\tau_{\wXX,R}(t))}\int_0^{t_1}\dots \int_{0}^{t_{K-2}}\nabla^{ K-1} F(t_{K-1}, \ZZ)d \ZZ_{t_{K-1}} \dots  d \ZZ_{t_1}$$  is defined as iterated Young integrals.
\end{theorem}
The proof of the  theorem is given in the Appendix \ref{prof th: taylor 12}.
\begin{remark}
\begin{enumerate}
    \item Notice that for $\wXX$ continuous, the remainder term becomes  
    $$\mathcal{R}^F_{K-1}(t,\wXX):=\int_0^{t}\int_0^{t_1}\dots \int_{0}^{t_{K-2}}\nabla^{ K-1} F(t_{K-1}, \wXX)d \wXX_{t_{K-1}} \dots  d \wXX_{t_1}.$$

    \item 
    Since the proof of Theorems \ref{th: taylor 12} relies on Theorem \ref{3th: Ito formula p 12},  possible modifications of the conditions on the functional $F$ discussed in Remark \ref{3rem: Ito formula p12} also apply here.
    \end{enumerate}
\end{remark}

\subsection{The case \texorpdfstring{$p\in [2,3)$}{}}\label{sec: taylor p23}

\begin{theorem}\label{3th: taylor p23 continuous}
     Let $p\in[2,3)$, $K\geq 3$, and $F:[0,1]\times D^{p}([0,1],G^{2}(\Rset^{d+1}))\rightarrow \Rset$ be a $C^K$ non-anticipative Marcus canonical  path functional. Assume that  
     $F,\nabla ^{}F, \dots, \nabla ^{K-1}F $, and $R^{\nabla^{K-2} F,\nabla^{K-1} F}$ are $[p,3)$\text{-}var continuous. Let $\widehat \XX\in \widehat D^{p}([0,1],G^{2}(\Rset^{d+1})) $, denote by $\wX$ its signature, and by $\ZZ$ its time-extended Marcus-transformed path with respect to some pair $(R,\psi_R)$. Then, for all 
$t\in [0,1]$,  
   \begin{align}\label{3eq: taylor expansion p23}
    F(t,\widehat \XX)=&\sum_{j=0}^{K-3}\nabla^{j}F (0,\widehat \XX)\widehat \X^{(j)}_{0,t}+ \mathcal{R}^F_{K-2}(t,\widehat \XX),
\end{align}
where $$\mathcal{R}^F_{K-2}(t,\wXX):=\int_0^{\psi_R(\tau_{\wXX,R}(t))}\int_0^{t_1}\dots \int_{0}^{t_{K-3}}\nabla^{ K-1} F(t_{K-1}, \ZZ)d \ZZ_{t_{K-1}} \dots  d  \ZZ_{t_1}$$ is defined as iterated rough integral.
\end{theorem}

The proof of the  theorem is given in the Appendix \ref{prof 3th: taylor p23 continuous}.
\begin{remark}\phantomsection\label{remark: taylor p23}
\begin{enumerate}
\item Notice that for $\wXX$ continuous, the remainder term explicitly reads as 
    $$\mathcal{R}^F_{K-1}(t,\wXX):=\int_0^{t}\int_0^{t_1}\dots \int_{0}^{t_{K-2}}\nabla^{ K-1} F(t_{K-1}, \wXX)d \wXX_{t_{K-1}} \dots  d \wXX_{t_1}.$$

Moreover, the iterated rough integrals determining this remainder term are defined as follows. Set $G^{-1}(\cdot,\wXX):=\nabla ^{K-1}F(\cdot,\wXX)$, $G^0(\cdot,\wXX):=\nabla ^{K-2}F(\cdot,\wXX)$, and 
\begin{align}\label{3eq: Gj}
G^j(\cdot,\wXX):=\int_0^\cdot\int_0^{t_{1}}\dots \int_0^{t_{j-1}} \nabla ^{K-2}F(t_{j},\wXX)d\wXX_{t_{j}}\dots d\wXX_{t_{1}},
\end{align}
for $j=1,\dots,K-2$. Then,  $G^j(\cdot,\wXX)$ is the rough integral of $(G^{j-1}(\cdot,\wXX),G^{j-2}(\cdot,\wXX))$ with respect to $\wXX$.
    \item \label{rem: weaken assumptions TF p 23}  Since the proofs of Theorem~\ref{3th: taylor p23 continuous} rely on Theorem~\ref{3th: Ito formula p 23}, discussions in Remark~\ref{3rem: Ito formula p23}\ref{3rem: Ito formula p23 item i} regarding a possible modifications of the conditions on the functional $F$ also apply here.
    \item  \label{rem: valentin dupire}   Adopting and combining the pathwise framework pioneered by \cite{F:81}  with the (functional) differential calculus introduced in \cite{D:09},  a pathwise functional Taylor expansion in terms of signature for \textit{continuous one-dimensional time-extended} paths of finite quadratic variations  has been proved in \cite{DT:23}.  (See Theorem 3.10 therein).
However, recalling  Remark~\ref{rem: non commutativity derivate} on the commutativity of the derivation order, the dependence in \cite{DT:23}  with respect to the time and path component is captured via the horizontal and vertical derivatives, respectively. Indeed, in their framework, only the horizontal and vertical derivatives do not commute (see also the discussions at page ~33 in \cite{CR:12}  and on page~10 in \cite{BD}). Therefore, considering time-extended one-dimensional paths, they can indeed recover all terms in the signature that appear in the expansion. In a multidimensional setting, however, since the higher-order (mixed) vertical derivatives commute, any expansion in terms of the signature can not be obtained.

The Taylor expansion thus further emphasizes the importance of establishing a differential calculus on path functionals that allows a non-commutative order of differentiation, without which it would not be possible to achieve Taylor expansions in terms of the signature in higher dimension.

\end{enumerate}

\end{remark}

\subsection{Analytic non-anticipative path functionals }

Building on the Taylor expansions in Theorem \ref{3th: taylor p23 continuous},  we mimic here the classical analysis approach to the study of smooth functions and introduce a notion of analytical path functional. We will address these concepts in the case of \cadlag $p$-rough path for $p\in [2,3)$. (Note that for a  simpler setting of \cadlag $p$-rough path for $p\in [1,2)$, it can be similarly formulated). 

We start by  the definition of the concatenation path. 
\begin{definition}
 Let $\XX=\exp^{(2)}(X,A^X)$ and $ \YY=\exp^{(2)}(Y,A^Y)\in D([0,1],G^{2}(\Rset^d))$. We define the \textit{concatenation path} of $\XX$ and $\YY$ at time $t\in[0,1]$ as the path $$\XX \oplus_t \YY \in D([0,1],G^{2}(\Rset^d)),$$ defined as follows: for all $u\in[0,1]$
    \begin{align*}
        (\XX \oplus_t \YY)_u:=\exp^{([p])}((X_u,A^X_u)1_{\{u\leq t\}})\otimes \exp^{([p])}((Y_u-Y_t+X_t,A^Y_u-A^Y_t-A^X_t)1_{\{u>t\}}). 
    \end{align*}

\end{definition}
 Next, we use the concept of concatenation path to derive the Taylor expansion of a path functional within a neighborhood of a given path.  The proof  of the following corollary is exactly the same as the proof of Theorem~\ref{3th: taylor p23 continuous}, hence omitted.

\begin{corollary}
 Let $p\in[2,3)$, $K \geq 3$ and $F:[0,1]\times D^{p}([0,1],G^{2}(\Rset^{d+1}))\rightarrow \Rset$ be a $C^K$ non-anticipative path functional such that $F,\nabla ^{}F, \dots, \nabla ^{K-1}F $, and $R^{\nabla^{K-2} F,\nabla^{K-1} F}$ are $[p,3)$\text{-}var continuous. Then, for all $\wXX\in \widehat{C}^p([0,1],G^{2}(\Rset^{d+1}))$,  $\wYY\in \widehat{C}^p([0,1],G^{2}(\Rset^{d+1}))$, $t,s\in [0,1]$, $t\leq s$, it holds that
    \begin{align}\label{3eq: taylor p23 concatenated}
    F(s,\wXX \oplus_t \wYY )=&\sum_{j=0}^{K-3}\nabla^{j}F (t,\wXX)\wY^{(j)}_{t,s}+ \mathcal{R}^F_{K-2}(s,\wXX \oplus_t \wYY),
\end{align}
where $\mathcal{R}^F_{K-2}(s,\wXX \oplus_t \wYY):=\int_t^s\int_t^{s_1}\dots \int_{t}^{s_{K-4}}\nabla^{K-2} F(s_{K-2},\wXX \oplus_t \wYY)d\wYY_{s_{K-2}} \dots d\wYY_{s_1}$, and $\wY$ denotes the signature of $\wYY$. 
\end{corollary}
\begin{remark}
 The formula~\eqref{3eq: taylor p23 concatenated} for $F$ evaluated at a \cadlag concatenation path can be derived as in the Step 3 of the proof of Theorem \ref{3th: taylor p23 continuous}, and results in an expansion whose remainder term $\mathcal{R}^F_{K-2}$ is expressed in terms of the time-extended Marcus-transformed path of the concatenation path. 
\end{remark}

Recall from Remark~\ref{remark: taylor p23}\ref{rem: weaken assumptions TF p 23} that the formula~\eqref{3eq: taylor p23 concatenated} can be derived also under the assumption of  $F$ being a $C^K$-non-anticipative path functional such that $$F,\nabla ^{}F, \dots, \nabla ^{K-1}F, \nabla ^{K}F$$ are $[p,3)$\text{-}var continuous. Now, we use this slightly stronger assumption in order to introduce the notion of analytical and entire path functional. In the following, we say that a Marcus canonical path functional is $C^\infty$  if it is $C^K$, for all $K\in \Nset$.

\begin{definition}
    Let $p\in [2,3)$ and $F:[0,1]\times D^p([0,1],G^{2}(\Rset^{d+1}))\rightarrow \Rset$ be a $C^\infty$ non-anticipative Marcus canonical  path functional. Assume  that all its vertical derivatives  $\nabla ^{j}F$ are $[p,3)$\text{-}var continuous. Fix $(t,\wXX)\in [0,1]\times \widehat{D}^p([0,1],G^{2}(\Rset^{d+1}))$. 
     \begin{enumerate}
      \item We say that $F$
 is \textit{real analytic} at $(t,\wXX)$ if there exists $\delta>0$ such that for all $(s,\wYY)\in[t,1]\times  \widehat{D}^p([0,1],G^{2}(\Rset^{d+1}))$ with  
 $|s-t|+\|\wYY\|_{p\text{-}var[t,s]}<\delta$,
 \begin{align*}
     F(s,\wXX \oplus_t \wYY)=&\sum_{j=0}^{\infty }\nabla^{j}F (t,\wXX)\wY^{(j)}_{t,s}.
 \end{align*}
\item We say that $F$
 is \textit{entire} if for all $(t,\wXX)\in  \widehat{D}^p([0,1],G^{2}(\Rset^{d+1}))$, 
 \begin{align*}
     F(t,\wXX)=\sum_{j=0}^{\infty }\nabla^{ j}F (0,\wXX)\wX^{(j)}_{0,t}.
     \end{align*}
    \end{enumerate}

 \end{definition}

 \begin{example}
   Let $\u\in T^N(\Rset^{d+1})$, for some $N\in \N$ and consider the linear functional of the signature $F^\u$ (see Section~\ref{3sec: linear functions as Marcus}). Then, by Proposition~\ref{prop: derivative linear func}, $F^\u$ is an entire functional and for all $(t,\wXX)\in  \widehat{D}^p([0,1],G^{2}(\Rset^{d+1}))$, setting $F^\u(t,\wXX):=\langle \u,\wX_t\rangle$, it holds that
 \begin{align}\label{eq: linear func pw expansion}
\langle \u,\wX_t\rangle=\sum_{j=0}^{N}\langle \u^{(j)},\wX_0\rangle \wX^{(j)}_{0,t}.
     \end{align}
    Equation~\eqref{eq: linear func pw expansion} coincides with the Chen's relation.
     More generally, we can also consider path functional of the form $F(t,\wXX):=f(\langle \u,\X_t\rangle)$ for an entire function $f:\Rset\rightarrow\Rset$. Similar computations as in~\eqref{eq: f compones sig} show that $F$ is an entire path functional.

\end{example}

\appendix

\section{Proofs of Section \ref{sec: funcionals of RP}}

\subsection{Proof of Proposition \ref{3prop: F Marcus canonical then F invariant}}\label{proof 3prop: F Marcus canonical then F invariant}

    Let $\XX\in D([0,1],G^{[p]}(\Rset^d))$ and a time-reparametrization  $\phi$, one can show that there exists another time-reparametrization $\eta$ such that ${\ZZ}_\eta=\widetilde{\ZZ}$, where
${\ZZ}$ and $\widetilde{\ZZ}$
   denote the Marcus transformations of the \cadlag paths $\XX_{\phi}$ and $\XX$ with respect to some pairs $(R,\psi_R)$ and  $(\widetilde{R},\psi_{\widetilde{R}})$, respectively. Since $\eta$ can be chosen to satisfy
    \begin{align*}
 \eta(\psi^{-1}_{\widetilde{R}}(\tau_{\XX,\widetilde{R}}(\phi(t))))=\psi_R^{-1}(\tau_{\XX_\phi,R}(t)),
    \end{align*}
    for all $t\in [0,1]$, the claim follows as in equation~\eqref{3eq: Marcus functional independent on the pair}. 
\qed

\subsection{Proof of Proposition \ref{3prop: F Marcus canonical UiF Marcus canonical}}\label{proof 3prop: F Marcus canonical UiF Marcus canonical}

    Fix  $(t,\XX)\in [0,1]\times C([0,1],G^{[p]}(\Rset^{d}))$, $\phi$  a time-reparametrization and $\xi\in \Rset^d$. Since $F\in \mathcal{M}^1_{[p]}\subset \mathcal{M}^0_{[p]}$, by Proposition~\ref{3prop: F Marcus canonical then F invariant},
     
    \begin{align*}
            F(\phi(t),\XX\otimes \exp^{([p])}(\xi)1_{\{\cdot\geq \phi(t)\}})=&F(\phi(t),\XX_{\phi\circ\phi^{-1}}\otimes \exp^{([p])}(\xi)1_{\{\phi^{-1}(\cdot)\geq t\}})\\
            =&F(t,\XX_{\phi}\otimes \exp^{([p])}(\xi)1_{\{\cdot\geq t\}}),
        \end{align*}
        which implies condition \ref{item def Mcpf i} in Definition~\ref{3def: Marcus canonical path functional}.  The property \ref{item def Mcpf ii} is inherited by the functional $F$. Finally, condition \ref{item def Mcpf iii} follows by definition of the functional given in equation~\eqref{3eq: pre UiF path functional}. 
  \qed

  \subsection{Proof of Proposition \ref{prop: derivative linear func}}\label{proof prop: derivative linear func}

 ~\ref{item 1 Fu}: Property \ref{item def Mcpf i} of Definition \ref{3def: Marcus canonical path functional} follows from Corollary~\ref{coro: equation signature} and the property of the Young and rough integral, for which it holds that for all $\mathbf{X}\in C^p([0,1],G^{[p]}(\Rset^d))$  and $t\in [0,1]$,
    \begin{align*}
        \int_0^{\phi(t)}Y_sd\XX_s= \int_0^{t}Y_{\phi(s)}d\XX_{\phi(s)},
    \end{align*}
    for suitably chosen paths $Y$
    (or $(Y,Y')$) such that the integrals are well defined. Condition \ref{item def Mcpf ii} follows from Chen's relation.
 The property \ref{item def Mcpf iii}  of Definition \ref{3def: Marcus canonical path functional} follows from the definition of the signature of a weakly geometric~\cadlag rough path.

 ~\ref{item 2 Fu}: Fix $\XX \in D^p([0,1],G^{[p]}(\Rset^d))$ and let $\widetilde{\XX}$ denote its Marcus-transformed path with respect to some pair $(R,\psi_R)$, with signature $\widetilde{\X}$, and for $t\in [0,1]$, let $\mu_t\in [0,1]$ be defined as in Notation~\ref{3notation: mut}. Fix $i=1,\dots,d$, $h\in \Rset$, and consider the vertically perturbed path
\begin{align}\label{3eq: vertical perturbed path}
    \YY:=\widetilde{\XX}\otimes \exp^{([p])}(h\e_i)1_{\{\cdot \geq \mu_t\}}\in D^p([0,1],G^{[p]}(\Rset^d)).
\end{align}
By Chen's relation and the minimal jump extension property (see equation~\eqref{1eqn11}), it holds that the signature  of $\YY$, denoted as $\Y$, at time $\mu_t$ reads as follows:
\begin{align*}
    \Y_{\mu_t}= \Y_{\mu_t^-}\otimes \Delta\Y_{\mu_t}=\widetilde{\X}_{\mu_t}\otimes \exp(h\e_i).
\end{align*}
Notice that we use that the path $s\mapsto \widetilde{\X}_s$  has continuous components. We then get
\begin{align*} 
U^iF^\u(t,\XX):=&\frac{d}{dh}F^\u(\mu_t,\widetilde{\XX}\otimes \exp^{([p])}(h\e_i)1_{\{\cdot \geq \mu_t\}})|_{h=0}\\=&\frac{d}{dh}\langle  \u,\widetilde{\X}_{\mu_t}\otimes \exp(h\e_i)\rangle |_{h=0}.
\end{align*}
An explicit computation yields $ U^iF^\u(t,\XX)=\langle \u^{(1)}_{(i)},\widetilde{\X}_{\mu_t}\rangle$ for $\u^{(1)}_{(i)}\in T(\Rset^d)$ introduced in \eqref{3eq: shifts}. Moreover, by definition of the signature,
$\langle \u^{(1)}_{(i)},\widetilde{\X}_{\mu_t}\rangle=\langle \u^{(1)}_{(i)},\X_t\rangle$. Thus, computing the vertical derivative in all the directions $j=1,\dots,d$, we get $ \nabla F^\u(t,\XX)=\langle \u^{(1)},\X_t\rangle,$
and more generally, by iterating the same reasoning, for all $k\in \Nset$, $ \nabla^k F^\u(t,\XX)=\langle \u^{(k)},\X_t\rangle.$

~\ref{item 3 Fu}: It follows from Corollary 10.28 in~\cite{FV:10}. 
\qed

\section{Proofs of Section~\ref{sec: uat differential}}
 Before presenting the proofs of the main results of the section, we list some key notions and introduce a lighter notation. Recall from Section \ref{3sec: algebra} that
\begin{align*}
    G((\Rset^{d+1})):=\{\x\in T((\Rset^{d+1})) \ | \ \pi_{\leq N }(\x)\in G^N(\Rset^{d+1}) \text{ for all }N\in \Nset \},
\end{align*}
for $G^N(\Rset^{d+1}):=\exp^{(N)}(\mathfrak{g}^N(\Rset^{d+1}))$. Notice that $G^N(\Rset^{d+1})$ is a Carnot group (see Definition 2.2.1 in \cite{BLU:07}). 
Moreover, set
\begin{align*}
    \mathfrak{g}((\Rset^{d+1})):=\{ \x\in T((\Rset^{d+1})) \ | \ \pi_{\leq N }(\x)\in \mathfrak{g}^N(\Rset^{d+1}) \text{ for all }N\in \Nset \}.
\end{align*}
Then $G((\Rset^{d+1}))=\exp(\mathfrak{g}((\Rset^{d+1})))$ (see e.g., \cite{BBHRN:24}, \cite{S:22}) and is closed with respect to the tensor multiplication $\otimes$.

    Fix $F\in C^K$ and let $g$ be the map for which Assumption~\ref{assumpt 1} holds. For $\XX\in {C}^p([0,1], G^{[p]}(\Rset^{d+1}))$, let $g^{\XX,K}$ be the map introduced in equation~\eqref{eq: gk assumpt lemma}. To simplify the notation, we write $g^K:=g^{\XX,K}$ whenever no confusion arises, and write $\mathfrak{g}^K$  and $|\beta|_{\mathfrak{g}^K}$ in place of $\mathfrak{g}^K(\Rset^{d+1})$ and $|\beta|_{\mathfrak{g}^K(\Rset^{d+1})}$ (see  Notation~\ref{3notation derivatives lie algebra}), respectively. Finally, we denote with the index $0$ the first component of $\XX$.

\subsection{Proof of Lemma~\ref{lemma: partial derivatives gk}}\label{proof of lemma: partial derivatives gk}
 Fix $\XX\in {C}^p([0,1], G^{[p]}(\Rset^{d+1}))$ and let $(t,\bxi)\in [0,1]\times \mathfrak{g}^K$. By definition of signature and $G((\Rset^{d+1}))$,  
$$\X_t\otimes \exp(\bxi,0,\dots,0)\in G((\Rset^{d+1})).$$
Therefore, by Assumption~\ref{assumpt 1}, the map $g^K$ is well defined.  Next, fix $t\in [0,1]$. We show that for all $j=1,\dots,K$ and $|\beta|_{\mathfrak{g}^K}=j$, the derivatives $D^\beta_{\bxi} g^K(t,{\bxi})|_{{\bxi}=0}$ exist.

The claim follows if we prove that for all $j=1,\dots,K$, $i_1,\dots,i_j=0,1,\dots,d$, 
\begin{align}\label{eq: derivatives F equal deriv g}
  U^{i_{j}}\dots U^{i_1}F(t,\XX)=
  \frac{d^j}{dh_{j}\dots dh_1}g(\X_t\otimes \exp(h_{j}\e_{i_{j}})\otimes \dots \otimes \exp(h_1\e_{i_1}))|_{h_{j}=\dots=h_1=0}.
\end{align}
Indeed,  by assumption $F\in \Mcalp^K$, the quantities $ U^{i_{j}}\dots U^{i_1}F(t,\XX)$ are well defined.  Furthermore, by Propositions 20.1.7 (see in particular equation (20.20)) and 20.1.9 (see also Proposition 20.1.4) in \cite{BLU:07}, 
setting
\begin{align*}
    L(t,h_1\e_{i_1},\dots,h_j\e_{i_j}):=g(\X_t\otimes \exp(h_{j}\e_{i_{j}})\otimes \dots \otimes \exp(h_1\e_{i_1})),
\end{align*}
it holds that for all $|\beta|_{\mathfrak{g}^K}=j$, 
\begin{align*}
    D^\beta_{\bxi} g^K(t,{\bxi})|_{{\bxi}=0}\in\text{span}\Big\{\frac{d^j}{dh_{j}\dots dh_1}L(t,h_1\e_{i_1},\dots,h_j\e_{i_j})|_{h_1=\dots=h_j=0},\ i_1,\dots,i_j=0,\dots,d\Big\},
\end{align*}
from which we deduce that $ D^\beta_{\bxi} g^K(t,{\bxi})|_{{\bxi}=0}$ exists.
In order to show~\eqref{eq: derivatives F equal deriv g}, we
apply the iterative procedure to compute the higher order vertical derivatives of $F$.

Fix $\delta>0$ and for $i_K=0,1,\dots,d$, $h_K\in (-\delta,\delta)$ set
\begin{align*}
    &\XX^{[i_K]}(h_K):=\XX+h_{i_{K}} \e_{i_{K}}1_{\{\cdot \geq \mu_t^{[0]}\}},\\
    &\YY^{[i_K]}(h_K):=\reallywidetilde{\XX^{[i_K]}(h_K)},\\
    & \mu_t^{[i_K]}(h_K):=\psi_{R_K}^{-1}(\tau_{\XX^{[i_K]}(h_K),R_K}(t)),
\end{align*}
for some pair $(R_K,\psi_{R_K})$ that might depend on $h_K$.
For $j=K-1,\dots,1$, $i_j=0,1,\dots,d$, $h_j\in (-\delta,\delta)$, set
\begin{align*}
    &\XX^{[i_K,\dots,i_j]}(h_K,\dots,h_j):=\YY^{[i_K,\dots,i_{j+1}]}(h_K,\dots,h_{j+1})+h_{j} \e_{i_{j}}1_{\{\cdot \geq \mu_t^{[i_K,\dots,i_{j+1}]}\}},\\
    &\YY^{[i_K,\dots,i_j]}(h_K,\dots,h_j):=\reallywidetilde{\XX^{[i_K,\dots,i_j]}(h_K,\dots,h_j)},\\
    & \mu_t^{[i_K,\dots,i_j]}(h_K,\dots,h_j):=\psi_{R_j}^{-1}(\tau_{\XX^{[i_K,\dots,i_j]}(h_K,\dots,h_j),R_j}(\mu_t^{[i_K,\dots,i_{j+1}]}(h_K,\dots,h_{j+1}))),
\end{align*}
for $\reallywidetilde{\XX^{[i_K,\dots,i_j]}(h_K,\dots,h_j)}$ denoting the Marcus-transformed path of  $\XX^{[i_K,\dots,i_j]}(h_K,\dots,h_j)$ with respect to some pair $(R_j,\psi_{R_j})$ that might depend on $(h_K,\dots,h_j)$. 

Then, by definition of higher-order vertical derivatives (see Definition~\ref{def: higher order derivatives}), it holds that  for all $j=1,\dots,K$, 
\begin{align*}
    &U^{i_{j}}\dots U^{i_1}F(t,\XX)\\
    &=\frac{d^j}{dh_j\dots dh_1}F(\mu_t^{[i_K,\dots,i_1]}(0,\dots,0,h_j,\dots,h_1),\YY^{[i_K,\dots,i_j]}(0,\dots,0,h_j,\dots,h_1))|_{h_j=\dots=h_1=0}. 
\end{align*}
Finally, Assumption~\ref{assumpt 1} and a direct computation of the signature at time $\mu_t^{[i_K,\dots,i_1]}(h_K,\dots,h_1)$ of the continuous path $\YY^{[i_K,\dots,i_j]}(h_K,\dots,h_1)$ yield that
\begin{align*}
    F(\mu_t^{[i_K,\dots,i_1]}(h_K,\dots,h_1),\YY^{[i_K,\dots,i_j]}(h_K,\dots,h_1))=g(\X_t\otimes \exp(h_{j}\e_{i_{j}})\otimes \dots \otimes \exp(h_1\e_{i_1})).
    \end{align*}
The claim follows.\qed

\subsection{Proof of Theorem~\ref{3th: Nachbin for path functionals}}\label{3sec: proof Nachbin}
Before delving into the proof of Theorem~\ref{3th: Nachbin for path functionals} we outline its main steps. Recall the simplified notation introduced at the beginning of the present section.
\paragraph{Sketch of the proof.} Fix $\wXX\in \widehat{C}^p([0,1], G^{[p]}(\Rset^{d+1}))$.
\paragraph{Step 1:} We show that for $t\in[0,1]$, $i_1,\dots,i_j\in \{0,\dots,d\}^j$, $j=1,\dots,K$, 
        \begin{align*}
        &F(t,\wXX)=g^K(t,0),\\
           &U^{i_j}\dots U^{i_1}F(t,\wXX)\in \text{span}\{D^\beta_{\bxi} g^K(t,{\bxi})|_{{\bxi}=0} \ : \ \beta\in \Nset^M_0, \ |\beta|_{\mathfrak{g}^K}=j\}.
        \end{align*}
        
    \paragraph{Step 2:} Since $F\in C^K$, Assumption~\ref{assumpt 2}
     and an adaptation of the proof of the Weierstrass theorem (see Theorem 1.6.2 in~\cite{N:85}) to the present setting yield that there exists a sequence of polynomials $(p_n)_{n\in\Nset}$ on $[0,1]\times \mathcal{U}(0)$ such that  
 \begin{align*}
        \lim_{n\rightarrow\infty}\sup_{t\in [0,1]}\sup_{{\bxi}\in H}\sum_{\beta\in \Nset^M_0,\ |\beta|_{\mathfrak{g}^K}\leq K}\|D^\beta_{\bxi} g^K(t,{\bxi})-D^\beta_{\bxi} p_n(t,{\bxi})\|=0,
    \end{align*}
      for some compact set $H\subset \mathcal{U}(0)$ such that $0\in H$.

      \paragraph{Step 3:} We define $\YY\in {\widehat C}^p([0,1], G^{[p]}(\Rset^{d+2}))$ as the time-extended path of $\wXX$ (Section~\ref{sec: uat differential} and in particular Remark~\ref{rem: extended paths}~\ref{item 4 remark extended path}), and denote this auxiliary component by the index $-1$ and its signature by $\Y$. Let $p$  be a polynomial on $[0,1]\times \mathcal{U}(0)$. Notice that the map
      \begin{align*}
          [0,1]\ni t\mapsto  \tilde{p}(t)(\cdot):= \Big(\mathcal{U}(0)\ni {\bxi}\mapsto  p(t,{\bxi})\Big)
      \end{align*}
      belongs to $C\big([0,1],  C^{K}(\mathcal{U}(0))\big)$, where $C^K(\mathcal{U}(0))$ denotes the set of $K$ times differentiable functions on $\mathcal{U}(0)$, which we endow with the topology of uniform convergence on compacts of the function and all its derivatives up to order $K$. We exploit the Stone-Weierstrass theorem for vector-valued maps (see e.g., Theorem 3.3 in~\cite{CSTuat:23}) to show that there exists a sequence $(\v_n)_{n\in \Nset}\subset T(\Rset^{d+2})$ whose indices are in $\{-1,0,\dots,d\}$ such that
    \begin{align}\label{3eq: convergence VVSW}
        \lim_{n\rightarrow\infty }\sup_{t\in [0,1]}\sup_{{\bxi}\in H}\sum_{\beta\in \Nset^M_0,\ |\beta|\leq K}\|D^\beta_{\bxi} p(t,{\bxi})-D^\beta_{\bxi}\langle \v_n,\Y_t\otimes \exp(i({\bxi}),0,\dots,0)\rangle\|=0,
    \end{align}
    where $i:g^K(\Rset^{d+1})\rightarrow g^K(\Rset^{d+2})$ denotes the embedding of $g^K(\Rset^{d+1})$ into $g^K(\Rset^{d+2})$ obtained by setting equal to $0$ all the additional components.  

    \paragraph{Step 4:} A combination of Step 2 and Step 3 yields the existence of a sequence $(\v_n)_{n\in \Nset}\subset T(\Rset^{d+2})$ whose indices are in $\{-1,0,\dots,d\}$ such that
    \begin{align}\label{3eq: convergence step 4 }
        \lim_{n\rightarrow\infty }\sup_{t\in [0,1]}\sum_{\beta\in \Nset^M_0,\ |\beta|_{\mathfrak{g}^K}\leq K}\|D^\beta_{\bxi} g^K(t,{\bxi})|_{{\bxi}=0}-D^\beta_{\bxi}\langle \v_n,\Y_t\otimes \exp(i({\bxi}),0,\dots,0)\rangle|_{{\bxi}=0}\|=0.
    \end{align}

   \paragraph{Step 5:} We show that there exists  a sequence $(\u_n)_{n\in \Nset}\subset T(\Rset^{d+1})$ whose indices are in $\{0,\dots,d\}$ such that for all $n\in\Nset$, $t\in [0,1]$, $\beta\in \Nset^M_0$, $|\beta|_{\mathfrak{g}^K}\leq K$,
$$D^\beta_{\bxi}\langle \v_n,\Y_t\otimes \exp(i(\bxi),0,\dots,0)\rangle|_{{\bxi}=0}=D^\beta_{\bxi}\langle \u_n,\wX_t\otimes \exp(\bxi,0,\dots,0)\rangle_{|{\bxi}=0}.$$

\paragraph{\textbf{Step 6:}}For $\u\in T(\Rset^{d+1})$, consider the path functional given by
\begin{align*}
      [0,1]\times D^p([0,1],G^{[p]}(\Rset^d))\ni (t,\XX)\mapsto F^\u(t,\XX):=\langle \u,\X_t\rangle.
\end{align*}
Recall that for such functional, the map $g^K$ in~\eqref{eq: gk assumpt lemma} reads as $g^K(t,\bxi):=\langle \u,\wX_t\otimes \exp(\bxi,0,\dots,0)\rangle$. An application of Step 1 to the functional $F^\u$ yields that for all $t\in[0,1]$, $i_1,\dots,i_j\in \{0,\dots,d\}^j$, $j=1,\dots,K$, 
        \begin{align*}
        &F^\u(t,\wXX)=\langle \u,\wX_t\rangle,\\
           &U^{i_j}\dots U^{i_1}F^\u(t,\wXX)\in \text{span}\{D^\beta_{\bxi} \langle \u,\wX_t\otimes \exp(\bxi,0,\dots,0)\rangle|_{{\bxi}=0} \ : \ \beta\in \Nset^M_0, \ |\beta|_{\mathfrak{g}^K}=j\},
        \end{align*}
which concludes the first part of the proof.

\paragraph{\textbf{Step 7:}} For $\wXX\in \widehat{D} ^p([0,1], G^{[p]}(\Rset^{d+1}))$, we deduce the claim by exploiting the Marcus property of the functionals.

\begin{proof}
Fix $\wXX\in \widehat{{C}}^p([0,1], G^{[p]}(\Rset^{d+1}))$. 

\paragraph{\underline{Step 1}:} 
First of all notice that $F(t,\XX)=g(\X_t)=g^K(t,0)$ by Assumption~\ref{assumpt 1} and definition of $g^K$ in~\eqref{eq: gk assumpt lemma}. Then, an inspection of the proof of Lemma~\ref{lemma: partial derivatives gk} shows that for all $j=1,\dots,K$, $i_1,\dots,i_j=0,\dots,d$ 
\begin{align*}
  U^{i_{j}}\dots U^{i_1}F(t,\wXX)=\frac{d^j}{dh_{j}\dots dh_1}g(\wX_t\otimes \exp(h_{j}\e_{i_{j}})\otimes \dots \otimes \exp(h_1\e_{i_1}))|_{h_{j}=\dots=h_1=0}.
\end{align*}
Since by Propositions 20.1.4 and 20.1.5 (see in particular equation 20.19), in~\cite{BLU:07}, for all $j=1,\dots,K$,
 \begin{align}\label{3eq: span Dbeta g }
     \frac{d^j}{dh_{j}\dots dh_1} g(\wX_t\otimes \exp(h_{j}\e_{i_{j}})\otimes \dots& \otimes \exp(h_1\e_{i_1}))|_{h_{j}=\dots=h_1=0} \\
     &\in \text{span}\{D^\beta_{\bxi} g^K(t,{\bxi})|_{{\bxi}=0}
      \ : \beta\in \Nset^M_0,\ |\beta|_{\mathfrak{g}^K=j}\}.\nonumber 
 \end{align} The claim follows. 

\paragraph{\underline{Step 2}:} It follows directly by the assumption on $g^K$ and a simple adaptation of the proof of Theorem 1.6.2 in~\cite{N:85} to the present setting.

\paragraph{\underline{Step 3}:} Let $p$ be a polynomial on $[0,1]\times \mathcal{U}(0)$, $\YY\in \widehat{C}^p([0,1], G^{[p]}(\Rset^{d+2}))$ be the time-extended path of $\wXX$, and denote the  auxiliary time component by the index $-1$ and its signature by $\Y$. Consider the following set of maps
$$\mathcal{W}:=span\left\{ t \mapsto\left ( \mathcal{U}(0)\ni {\bxi}\mapsto  w(t)({\bxi}):=  \langle e_I, {\mathbb{Y}}_t \otimes \exp ( i({\bxi}),0,\dots,0)\rangle \right), |I | \geq 0\right\}.$$
We verify the hypothesis of the Stone -Weiestrass theorem for vector-valued maps to show the convergence in equation~\eqref{3eq: convergence VVSW}. We must prove that $\mathcal{W}$ satisfies the following conditions: \begin{enumerate}
	    \item \label{3VVSW 1} it is a $\mathcal{A}$-submodule of $C([0,1],C^K(\mathcal{U}(0)))$, for some point separating subalgebra $\mathcal{A }\subseteq C([0,1])$ that vanishes nowhere; 
 \item \label{3VVSW 2} for all $t\in [0,1]$, $\mathcal{W}(t):=\{w(t) \ \colon \ w\in \mathcal{W}\}$ is dense in $C^K(\mathcal{U}(0))$.
	\end{enumerate}
~\ref{3VVSW 1}:    Notice that $\mathcal{W} \subseteq C([0,1], C^{K} (\mathcal{U}(0))$, as the functions in $\mathcal{W}$ are simply linear combinations of polynomials in ${\bxi}$ whose coefficients are continuous in $t$. Let $\mathcal{A}$ denote the space of polynomials on $[0,1]$. Then, $\mathcal{W}$ is an $\mathcal{A}$-submodule since the map $[0,1]\ni t \mapsto p(t) w(t) \in  \mathcal{W}$ for every  $w\in \mathcal{W}$ and $p$ polynomial on $[0,1]$.
 
\noindent~\ref{3VVSW 2}: Fix $t\in [0,1]$. We show that the set $\mathcal{W}(t)$,  which explicitly reads as 
    \begin{align}\label{3eq: W(t)}
        \mathcal{W}(t)=\text{span}\{ \mathcal{U}(0)\ni{\bxi}\mapsto\langle \epsilon_I,{\mathbb{Y}}_t\otimes \exp(i({\bxi}),0,\dots,0)\rangle\colon |I|\geq0\},
    \end{align}
 is dense in $C^{K}( \mathcal{U}(0))$.  To this end, we apply the Nachbin theorem (see~\cite{N:49}). We need to verify that the set $\mathcal{W}(t)$ satisfies the following properties:
 \begin{enumerate}
	    \item\label{3Nachbin 1} it is a linear subspace of $C^K(\mathcal{U}(0))$;
	    \item\label{3Nachbin 2}it is a sub-algebra that contains a
	non-zero constant function and separates points;
	    \item\label{3Nachbin 3} for all $\bxi\in \mathcal{U}(0)$, $\y\in \mathfrak{g}^K(\Rset^{d+1})$ with $\y\neq 0$, there exists $f\in \mathcal{W}(t)$ such that $\nabla_{\bxi} f(\bxi) \y \neq 0$.
	\end{enumerate}
~\ref{3Nachbin 1}: It is clear.

 \noindent~\ref{3Nachbin 2}: The shuffle properties of group-like elements yields that $\mathcal{W}(t)$ is a sub-algebra (see equation \eqref{3eq:shuffle}). Moreover, it contains a map constantly equal to $1$:  $$\mathcal{U}(0)\ni{\bxi}\mapsto\langle \epsilon_\emptyset,{\mathbb{Y}}_t\otimes \exp(i({\bxi}),0,\dots,0)\rangle,$$ and separates point. Indeed, let ${\bxi}_1,{\bxi}_2\in \mathcal{U}(0)$ with ${\bxi}_1\neq {\bxi}_2$. If for some $i=0,\dots,d$, $\langle\e_i,{\bxi}_1\rangle \neq \langle\e_i,{\bxi}_2\rangle $, then $$\langle \e_i,\Y_t\otimes \exp(i({\bxi}_1),0,\dots,0)\rangle \neq \langle \e_i,\Y_t\otimes \exp(i({\bxi}_2),0,\dots,0)\rangle.$$ Otherwise, let $J:=(i_1,\dots,i_j)\in \{0,\dots,d\}^j$ $j=2,\dots,K$ be such that $\langle \e_J,{\bxi}_1\rangle \neq \langle \e_J,{\bxi}_2\rangle$ and $\langle \e_I,{\bxi}_1\rangle=\langle \e_I,{\bxi}_1\rangle$ for all $I\in \{0,\dots,d\}^i $ with $i<j$. Then,
     $$\langle \e_J,\Y_t\otimes \exp(i({\bxi}_1),0,\dots,0)\rangle \neq \langle \e_J,\Y_t\otimes \exp(i({\bxi}_2),0,\dots,0)\rangle.$$ 
\noindent~\ref{3Nachbin 3}: Let ${\bxi}\in \mathcal{U}(0)$ and $\y\in \mathfrak g^K(\Rset^{d+1})$ with $\y\neq 0$. 
If for some $i=0,\dots,d$, $\langle \e_i,\y\rangle\neq 0$, then 
\begin{align*}
    \nabla_{{\bxi}}\langle \e_i,\mathbb{Y}_t\otimes \exp(i({\bxi}),0,\dots,0)\rangle  \y =\langle \e_i,\y\rangle\neq 0.
 \end{align*}
If instead for all $i=0,\dots,d$, $\langle \e_i,\y\rangle\neq 0$, let $J:=(i_1,\dots,i_j)\in \{0,\dots,d\}^j$ $j=2,\dots,K$ be such that $\langle \e_J,\y\rangle \neq 0$ and $\langle \e_I,\y\rangle  =0$ for all $I\in \{0,\dots,d\}^i $ with $i<j$. Then,
\begin{align*}
    \nabla_{{\bxi}}\langle \e_J,\mathbb{Y}_t\otimes \exp(i({\bxi}),0,\dots,0)\rangle  \y =\langle \e_J,\y\rangle\neq 0.
 \end{align*}
The claim follows.
\paragraph{\underline{Step 4}:} A direct consequence of  Step 2 and Step 3.
\paragraph{\underline{Step 5}:} Let $(\v_n)_{n\in \Nset}\subset T(\Rset^{d+2})$ be the sequence  whose indices are in $\{-1,0,\dots,d\}$ and for which the convergence in~\eqref{3eq: convergence step 4 } holds. The assertion follows if we prove that
\begin{enumerate}
    \item \label{3lask K indices } we can extract a subsequence, that without loss of generality still denoted as $(\v_n)_{n\in\Nset}$, whose last $K$ indices of each of its terms differ from $-1$. 
\end{enumerate}
Recall that here $-1$ denotes the index corresponding to the auxiliary time component of $\mathbf{Y}$. Indeed, in such a case, one can find a sequence $(\u_n)_{n\in\Nset}\subset  T(\Rset^{d+1})$ whose indices are in $\{0,\dots,d\}$ such that for all $n\in \Nset$, $t\in [0,1]$, $\beta\in \Nset^M_0$, $|\beta|_{\mathfrak{g}^K}\leq K$,
\begin{align*}
    &\langle \mathbf{v}_n, \mathbb{Y}_t \rangle= \langle \mathbf{u}_n,  {\wX}_t \rangle,\\
     &D^\beta_{\bxi}\langle \v_n,\Y_t\otimes \exp(i({\bxi}),0,\dots,0)\rangle|_{{\bxi}=0}=D^\beta_{\bxi}\langle \u_n,\wX_t\otimes \exp({\bxi},0,\dots,0)\rangle_{|{\bxi}=0},
\end{align*}
concluding the proof.

In order to prove~\ref{3lask K indices }, assume first  that for all $j=1,\dots,K$, there exist $i_j\in \{0,\dots,d\}$ and $t_j\in [0,1]$ such that 
    \begin{align}\label{3eq: gradientneq0}
        \langle \e_{i_j}^{\otimes j},\nabla^j F(t_j, \wXX)\rangle \neq 0.
    \end{align}
Under this assumption, condition~\ref{3lask K indices } is verified. Indeed, otherwise, for some $j=1,\dots,K$, and some $N\in\Nset$, the last $j$-th component of each $\v_n$ with $n\geq N$ would be equal to $-1$. This would imply that for all $n\geq N$ and $t\in [0,1]$,
\begin{align*}
    D^\beta_{{\bxi}}\langle \mathbf{v}_n, \mathbb{Y}_t \otimes \exp (i({\bxi}),0,\dots,0)\rangle _{|{\bxi}=0}=0,
\end{align*}
for all  $\beta\in \Nset_0^M$ with $|\beta|_{\mathfrak{g}^K}=j$.
In particular, by Step 4, for all $t\in [0,1]$
  \begin{align}\label{eq: gbeta}
      D^\beta_{{\bxi}}g^K(t,{\bxi})|_{{\bxi}=0}=0.
    \end{align}
   Finally, by Step 1, \eqref{eq: gbeta} 
  contradicts~\eqref{3eq: gradientneq0} and the claim follows.

Assume  that condition~\eqref{3eq: gradientneq0} does not hold. This means that there exist $1 \leq m_1<\dots<m_l\leq K$ such that for all $i\in \{0,\dots,d\}$ and for all $t\in [0,1]$,  $  \langle \e_{i}^{\otimes m_k},\nabla^{m_k} F(t, \wXX)\rangle =0$, for all $k=1,\dots,l$. Let $M$ be the biggest of such $m_k$ and
consider the functional  ${F}_0:[0,1]\times D^{p}([0,1],G^{[p]}(\Rset^{d+1}))\rightarrow \Rset$ defined via $${F}_0(t,\XX):=F(t,\XX)+\langle \e_0^{\otimes M},\X_t\rangle,$$ 
 for all $(t,\XX)\in [0,1]\times D^{p}([0,1],G^{[p]}(\Rset^{d+1}))$, with $\X_t$ denoting the signature of $\XX$ at time $t$. Notice that  for all $t\in (0,1]$, and for all $m_k$, $  \langle \e_{0}^{\otimes m_k},\nabla^{m_k} {F}_0(t, \wXX)\rangle =\langle \e_0^{\otimes M-m_k},\wX_t\rangle\neq 0$.
Therefore, 
\begin{enumerate}
\item [a)] if $M=1$, condition~\eqref{3eq: gradientneq0} is verified for the functional ${F}_0$. An application of the previous steps to the functional $F_0$ yields that there exists a sequence $(\u^0_n)_{n\in \Nset}\subset T(\Rset^{d+1})$ whose indices are in $\{0,\dots,d\}$ such that
    \begin{align*}
        \lim_{n\rightarrow\infty }\sup_{t\in [0,1]}\sum_{\beta\in \Nset^M_0,\ |\beta|_{\mathfrak{g}^K}\leq K}\|D^\beta_{\bxi}g^K(t,{\bxi})|_{{\bxi}=0}&+D^\beta_{\bxi}\langle \e_0^{\otimes M},\wX_t\otimes \exp(\bxi,0,\dots,0)\rangle|_{\bxi=0}\\
       &- D^\beta_{\bxi}\langle \u^0_n,\wX_t\otimes \exp(\bxi,0,\dots,0)\rangle|_{{\bxi}=0}\|=0.
    \end{align*}
  The claim follows by considering the sequence $\u_n:=\u^0_n-\e_0^{\otimes M}$,  $n\in \Nset$;
    \item [b)] if $M>1$ and for all $1\leq j<M$, $i_j$ for which $F$ satisfies condition~\eqref{3eq: gradientneq0} at $\wXX$ is different from $0$, then condition~\eqref{3eq: gradientneq0} is verified for the functional ${F}_0$ and the claim follows as in a).
\item [c)] if $M>1$ and for some $1\leq j<M$, $i_j$ for which $F$ satisfies condition~\eqref{3eq: gradientneq0} is equal to $0$, denote such $js$ by $1\leq j_1<\dots<j_l<M$, $1\leq l<M$.
    If for all $k=1,\dots,l$
    \begin{align*}
        \langle \e_{0}^{\otimes j_{k}},\nabla^{j_{k}} F(t, \wXX)\rangle+\langle \e_0^{\otimes M-{j_{k}}},\wX_t\rangle \neq  0,
    \end{align*}
     for some $t\in [0,1]$, the claim follows as in a). Otherwise, let $J$ be the smallest of the $j_k$ for which  for all $t\in [0,1]$,
\begin{align*}
        \langle \e_{0}^{\otimes j_{k}},\nabla^{j_{k}} F(t, \wXX)\rangle+\langle \e_0^{\otimes M-{j_{k}}},\wX_t\rangle = 0.
    \end{align*}
  Consider the path functional $F_1(t,\XX):=F(t,\XX)+\alpha_1\langle \e_0^{\otimes M},\X_t\rangle$  for  $\alpha_1\in \Rset$, $\alpha_1\notin\{ 1,0\}$. Observe that for all $j_k\geq J$ and for all $t\in (0,1]$,
$  \langle \e_{0}^{\otimes j_k},\nabla^{j_k} F_1(t, \wXX)\rangle \neq 0.$

 Finally, if for some $j_{k}$, with $j_k<J$
   \begin{align*}
        \langle \e_{0}^{\otimes j_{k}},\nabla^{j_{k}} F(t, \wXX)\rangle+\alpha_1
        \langle \e_0^{\otimes M-{j_{k}}},\wX_t\rangle = 0,
    \end{align*}
  for all $t\in [0,1]$, let $J_1<J$ 
 be the smallest of such $j_k$. Consider the path functional $F_2(t,\XX):=F(t,\XX)+\alpha_2\alpha_1\langle \e_0^{\otimes M},\X_t\rangle$  for  $\alpha_2\in \Rset$, $\alpha_2\notin\{ 1,0,\frac{1}{\alpha_1}\}$ and observe  that for all $j_k\geq J_1$  and for all $t\in (0,1]$,
$  \langle \e_{0}^{\otimes j_k},\nabla^{j_k} F_2(t, \wXX)\rangle \neq 0.$ An iterative application of this reasoning yields that condition~\eqref{3eq: gradientneq0} is verified for a functional  of the form $\widetilde{F}(t,\XX):=F(t,\XX)+\alpha\langle \e_0^{\otimes M},\X_t\rangle$  for some  properly choosen $\alpha\in \Rset$, $\alpha\neq 0$. The claim follows as in a).
\end{enumerate}

\paragraph{\textbf{\underline{Step 6}:}} It is a simple application of Step 1 to the functional $F^\u$.

\paragraph{\textbf{\underline{Step 7}:}} Fix $\wXX\in \widehat{D}^p([0,1], G^{[p]}(\Rset^{d+1}))$ and denote by $\ZZ$ the Marcus-transformed path such that $\ZZ\in \widehat{C}^p([0,1], G^{[p]}(\Rset^{d+1}))$. Recall that by definition of tracking-jumps-extended paths, such  Marcus transformed path always exists (see Remark~\ref{3rem: Marcus tracking-jumps paths}). An application of the previous steps of the proof to $\ZZ$ yields that there exists $(\u_n)_{n\in \Nset}\in T(\Rset^{d+1})$ such that
    \begin{align}\label{3eq: Nachbin approximation Marcus cadlag }
        \lim_{n\rightarrow\infty }\sup_{t\in [0,1]}\sum_{j=0}^K\|\nabla ^jF(t, \ZZ)-\nabla ^jF^{\u_n}(t, \ZZ)\|=0. 
    \end{align}
    Thus, in particular, 
\begin{align*}
        \lim_{n\rightarrow\infty }\sup_{t\in [0,1]}\sum_{j=0}^K\|\nabla ^jF(\mu_t, \ZZ)-\nabla ^jF^{\u_n}(\mu_t, \ZZ)\|=0, 
    \end{align*}
    for $\mu_t$ introduced in Notation~\ref{3notation: mut}. Since we deal here with Marcus canonical  path functionals, the claim follows.
\end{proof}

\begin{remark}\phantomsection \label{3rem: after proof Nachbin}
\begin{enumerate}
    \item \label{3rem: after proof Nachbin i}Let $\wXX\in \widehat{D}^p([0,1], G^{[p]}(\Rset^{d+1}))$ and $(\u_n)_{n\in\Nset}\in T(\Rset^{d+1})$ be the sequence such that the convergence in equation~\eqref{3eq: Nachbin approximation Marcus cadlag } holds. Observe that by the construction of the Marcus transformed path, for all $j=0,\dots,K$,
    \begin{align*}
        \sup_{s\in [0,1]}\sup_{n\in \Nset}\sup_{\theta\in [0,1]}\|\langle \u_n^{(j)},\wX_{s^-}\otimes \exp(\theta \log^{([p])}(\Delta \wXX_s))\rangle\|<\infty.
    \end{align*}

\end{enumerate}

\end{remark}

\section{Proofs of Section \ref{3sec: functional Ito formula}}
\subsection{Proof of Lemma \ref{3lemma: ito formula uniformly bounded p 12}}\label{proof lemma ito formula p12}
Recall that for all $n\in \Nset$, $j=0,\dots,2$, $\nabla ^jF^{\u_n}(t,\widehat \XX)=\langle \u_n^{(j)},\widehat \X_t\rangle$, where $\wX$ denotes the signature of $\wXX$.

~\ref{3item: lemma 1 i}: Let $(t_i)_i\subset [0,1]$ be the grid such that $\widehat \XX$ is a linear path on $[t_i,t_{i+1}]$. Notice that by Chen's relation and the invariance property of the signature under reparametrization (see e.g., the proof of Proposition~\ref{prop: derivative linear func}) for every $u\in [t_i,t_{i+1}]$,  $$  \widehat \X_u=\widehat \X_{t_i}\otimes \exp(\theta \widehat  \XX_{t_i,t_{i+1}}),$$
    for some $\theta \in [0,1]$.  Fix $s_1,s_2\in [t_i,t_{i+1}]$, $s_1\leq s_2 $, for some $i$. Then,  
    \begin{align*}
        \langle\u_n^{(1)},\widehat \X_{s_2} \rangle-\langle\u_n^{(1)},\widehat \X_{s_1} \rangle=g(1)-g(0),
    \end{align*}
    for $g:[0,1]\rightarrow\Rset^d$ such that  $ g(\theta):= \langle\u_n^{(1)},\widehat \X_{s_1}\otimes \exp(\theta \widehat  \XX_{s_1,s_2})\rangle $, for all $\theta \in [0,1]$.
    Thus, a first-order Taylor expansion of each of the components of $g$ at $0$ yields that 
    \begin{align}\label{3eq: mean value}
        \|\langle\u_n^{(1)},\widehat \X_{s_2} \rangle-\langle\u_n^{(1)},\widehat \X_{s_1} \rangle\|\leq \sup_{\theta\in [0,1]}\|\langle\u_n^{(2)},\widehat \X_{s_1} \otimes \exp( \theta \widehat\XX_{s_1,s_2})\rangle\|\|  \widehat \XX_{s_1,s_2}\|.
    \end{align}
    Since by~\eqref{eq: approx F and derivatives p 12}, 
    \begin{align*}
       \sup_{n\in\Nset} \sup_{s_1,s_2\in [t_i,t_{i+1}]}\sup_{\theta \in [0,1]}\|\langle\u_n^{(2)},\widehat \X_{s_1} \otimes \exp( \theta \widehat\XX_{s_1,s_2})\rangle\|<\infty,
    \end{align*}
    equation~\eqref{3eq: mean value} and $\|\widehat \XX\|_{1\text{-}var}<\infty $ yield that $  \sup_{n\in \Nset}\|\langle \u_n^{(1)},\widehat \X_\cdot\rangle\|_{1\text{-}var[t_i,t_{i+1}]}<\infty.$
    Repeating the same reasoning on every $[t_i,t_{i+1}]$, the claim follows.

 ~\ref{3item: lemma 1 ii}: It follows from~\ref{3item: lemma 1 i}, \eqref{eq: approx F and derivatives p 12}, and interpolation (see Lemma 5.12 and Lemma 5.27 in~\cite{FV:10}).
\qed

\subsection{Proof of Lemma \ref{3lemma: ito formula uniformly bounded p 23}}\label{proof lemma ito p23}

     ~\ref{3item lemma 2 i}: Let $(t_i)_i\subset [0,1]$ be the grid such that $\widehat X:=\pi_1(\wXX)$ is a linear path on $[t_i,t_{i+1}]$, and notice that for every $u\in [t_i,t_{i+1}]$, $  \widehat \X_u=\widehat \X_{t_i}\otimes \exp(\theta \widehat{X}_{t_i,t_{i+1}}),$
    for some $\theta \in [0,1]$. Fix   $s_1,s_2\in [t_i,t_{i+1}]$, $s_1\leq s_2 $, for some $i$, and consider the maps $g:[0,1]\rightarrow(\Rset^{d+1})^{\otimes 2}$, $ g':[0,1]\rightarrow\Rset^d$  given by 
    \begin{align*}
        &g(\theta):= \langle\u_n^{(2)},\widehat \X_{s_1}\otimes \exp(\theta \widehat  X_{s_1,s_2})\rangle, \\
        &g'(\theta):= \langle\u_n^{(1)},\widehat \X_{s_1}\otimes \exp(\theta \widehat  X_{s_1,s_2})\rangle,
    \end{align*}
    for all $\theta \in [0,1]$.
    A first and second order Taylor expansion of the components of $g$ and $g'$ at $0$ respectively yields that 
    \begin{align}\label{eq: reminder linear}
        &\|\langle\u_n^{(2)},\widehat \X_{s_2} \rangle-\langle\u_n^{(2)},\widehat \X_{s_1} \rangle\|\leq\sup_{\theta\in [0,1]}\|\langle\u_n^{(3)},\widehat \X_{s_1} \otimes \exp( \theta \widehat X_{s_1,s_2})\rangle\|\|  \widehat X_{s_1,s_2}\|,\nonumber \\
        &  \|  R^{\u_n^{(1)},\u_n^{(2)}}((s_1,s_2),\wXX)\|= \|\langle\u_n^{(1)},\widehat \X_{s_2}\rangle- \langle\u_n^{(1)},\widehat \X_{s_1}\rangle- \langle\u_n^{(2)},\widehat \X_{s_1}\rangle \widehat X_{s_1,s_2}\|\\ &  \qquad \qquad  \qquad \qquad \qquad \leq  
        \frac{1}{2}\sup_{\theta\in [0,1]}\| \langle\u_n^{(3)},\widehat \X_{s_1}\otimes \exp( \theta \widehat  X_{s_1,s_2})\rangle\| \|\widehat  X_{s_1,s_2}\|^2. \nonumber 
    \end{align}
 The claim follows as in the proof of~\ref{3item: lemma 1 i} of Lemma~\ref{3lemma: ito formula uniformly bounded p 23}.

  ~\ref{3item lemma 2 ii}: The assertion follows from~\ref{3item lemma 2 i} and interpolation (see Lemma 5.12 and Lemma 5.27 in~\cite{FV:10}).
    \qed

    \subsection{Proof of Corollary \ref{3coro: ito formula follmer}}\label{proof 3coro: ito formula follmer}

Let $\widehat \XX\in \widehat{D}^{p}([0,1],G^{2}(\Rset^{d+1})) $ be a Marcus-like weakly geometric rough path such that $\pi_1(\wXX)=\widehat{X}$.  By Theorem~\ref{3th: Nachbin for path functionals},
\begin{align*}
    \{s\in(0,1] \ \colon \ \Delta \nabla ^2 F(s,\wXX)\neq 0\}=\{s\in(0,1] \ \colon \ \Delta \wXX_s\neq 0\},
\end{align*}
and $\widehat{\XX}$ is Marcus-like, 
by Remark~\ref{3rem: quadratic variation}, 
\begin{align}\label{eq: convergence integral qv follmer}
    \lim_{n\rightarrow\infty }\sum_{s^n_i\in \pi^n_{[0,1]}}\nabla^2 F(s^n_i,\wXX) \widehat{X}^{\otimes 2}_{s^n_i,s^n_{i+1}}=\int_0^t\nabla^2F(s,\wXX)d[\widehat{X},\widehat{X}]^c_s+\sum_{0<s\leq t}\nabla^2F(s,\wXX)\Delta \widehat{X}^{\otimes 2}_{s}.
\end{align}
Moreover,  by Theorem~\ref{3th: Ito formula p 23}, $\big( \nabla F(\cdot ,\wXX), \nabla^{2}F(\cdot ,\wXX)\big)\in \mathcal{V}^{p',r}_{\wXX}$, for all $p'>p$ such that $\frac{1}{p'}+\frac{2}{p}>1$. Since for all $t\in [0,1]$, 
     $\nabla^2 F(t,\wXX)=\text{Sym}(\nabla^2 F(t,\wXX))$, the existence of the limit in equation~\eqref{3eq: coro follmer limit integral} follows by definition of the rough integral. Finally,  formula~\eqref{3eq: ito formula follmer} is derived from equation~\eqref{3eq: ito formula p 23}.
\qed 

\subsection{Proof of Corollary \ref{3coro: ito RIE }}\label{proof 3coro: ito RIE}

    By definition $\wXX$ is a Marcus-like rough path, therefore 
    \begin{align*}
          \{s\in(0,1] \ \colon \ \Delta \wXX_s\neq 0\}=\{s\in(0,1] \ \colon \ \Delta \widehat{X}_s\neq 0\}.
    \end{align*}
    Moreover, by Proposition 2.14 in~\cite{ALP:23}, 
    \begin{align*}
        \{s\in(0,1] \ \colon \ \Delta \widehat{X}_s\neq 0\} \subseteq \bigcup_{n\in \N}\pi^n_{[0,1]}.
    \end{align*}
Since by Theorem~\ref{3th: Nachbin for path functionals}, for $j=1,2$, $ \{s\in(0,1] \ \colon \ \Delta \nabla ^j F(s,\wXX)\neq 0\}=\{s\in(0,1] \ \colon \ \Delta \wXX_s\neq 0\},$
the claim follows by Theorem~\ref{3th: Ito formula p 23} and Proposition~\ref{3prop: RIE integal properties}.
\qed

\section{Proofs of Section \ref{sec: functional taylor}}

\subsection{Proof of Theorem \ref{th: taylor 12}}\label{prof th: taylor 12}
We first prove the result for $F$ evaluated at some continuous path (Step 1) and then extend it to \cadlag paths by applying the Marcus transformation (Step 2).

\paragraph{Step 1}: Let $\widehat \XX\in  \widehat C^{p}([0,1],G^1(\Rset^{d+1}))$. We prove the assertion by induction. If $K=2$, then the expansion~\eqref{3eq: taylor expansion p12} is simply the statement of Theorem~\ref{3th: Ito formula p 12}. Assume $K>2$ and that the assertion holds  for all $l=2,\dots,K-1$. By assumption $F$ is $C^K$, and is in particular $C^{K-1}$, thus by the induction hypothesis, 
\begin{align*}
    F(t,\widehat \XX)=&\sum_{j=0}^{K-3}\nabla^{j}F (0,\widehat \XX)\widehat \X^{(j)}_{0,t}+ \mathcal{R}^F_{K-2}(t,\widehat \XX),
\end{align*}
for $\mathcal{R}^F_{K-2}(t,\widehat \XX):=\int_0^{t}\int_0^{t_1}\dots \int_{0}^{t_{K-3}}\nabla^{ K-2} F(t_{K-2}, \widehat \XX)d \widehat \XX_{t_{K-2}} \dots  d \widehat \XX_{t_1}$.
    Next, since $\nabla^{K-2}F$ is $C^2$ (see Remark~\ref{remark: CK functionals}) and $\nabla ^{K-2}F, \nabla ^{K-1}F $ are $[p,2)$\text{-}var continuous, a component-wise application of Theorem \ref{3th: Ito formula p 12} to $ \nabla^{K-2}F$ yields that for all $s\in [0,1]$,
    \begin{align}\label{3eq: Taylor p12 linear from  k-2 to k-1}
         \nabla^{K-2}F(s,\widehat \XX)=  \nabla^{K-2} F(0,\widehat \XX)+\int_0^s \nabla^{K-1} F(r,\widehat \XX)d\widehat \XX_r.
    \end{align}
    The claim follows by replacing the expression in equation~\eqref{3eq: Taylor p12 linear from  k-2 to k-1} into  
   $\mathcal{R}^F_{K-2}(t,\widehat \XX)$.

   \paragraph{Step 2}: Let $(t,\widehat \XX)\in [0,1]\times  \widehat D^{p}([0,1],G^{1}(\Rset^{d+1})) $ and $\ZZ$ the time-extended Marcus-transformed path of $\wXX$. By Step 1 the expansion~\eqref{3eq: taylor expansion p12} holds when evaluating the $F$ at $(\psi_{R}^{-1}(\tau_{\wXX,R}(t)),\ZZ)$. Since for all $j=0,\dots,K-2$, $\nabla^jF$, $\wX^{(j)}$ are Marcus canonical  path functional, the claim follows.
\subsection{Proof of Theorem \ref{3th: taylor p23 continuous}}\label{prof 3th: taylor p23 continuous}

We split the proof into three steps. We first derive the Taylor expansion for functionals evaluated at the truncated signature at level $2$ of a time-extended piecewise linear path (Step 1). 
Then, we extend it to continuous paths by a density argument  (Step 2), and finally, we address the general \cadlag case  (Step 3). Starting from functionals evaluated at the truncated signature at level $2$ of a time-extended piecewise linear path simplifies the proof development as it enables working with Young integrals (and thus Riemann sums) instead of truly rough integrals (and thus compensated Riemann sums, see Proposition~\ref{3prop: rough integral}) in the derivation of the remainder term $\mathcal{R}_{K-2}^F$. This is possible because any rough integral with respect to the level-2 signature of a $\Rset^d$-valued path of finite variation is simply a Young integral with respect to the path itself (see e.g., Lemma~\ref{lemma: consistency}).

\paragraph{Step 1:}   Let $\widehat \XX\in  \widehat C^{1}([0,1],G^2(\Rset^{d+1}))$ be the truncated signature at level $2$ of a time-extended piecewise linear path. We prove the assertion by induction. If $K=3$, then the expansion in equation~\eqref{3eq: taylor expansion p23} is simply the statement of Theorem~\ref{3th: Ito formula p 23}. Assume $K>3$ and that the assertion holds true for all $j=3,\dots,K-1$. Since $F$ is in $C^K$, it is in particular in $C^{K-1}$. By the induction hypothesis,  \begin{align*}
    F(t,\wXX)=&\sum_{j=0}^{K-4}\nabla^{j}F (0,\widehat \XX)\widehat \X^{(j)}_{0,t}+ \mathcal{R}^F_{K-3}(t,\widehat \XX),
\end{align*}
for $\mathcal{R}^F_{K-3}(t,\widehat \XX):=\int_0^{t}\int_0^{t_1}\dots \int_{0}^{t_{K-4}}\nabla^{ K-3} F(t_{K-3}, \widehat \XX) d \widehat \XX_{t_{K-3}} \dots  d \widehat \XX_{t_1}$. Thus, the claim follows
if we prove that 
\begin{align}\label{3eq: eq: eq3 step1 proof taylor p23}
    \mathcal{R}^F_{K-3}(t,\widehat \XX)=\nabla^{K-3}F(0,\wXX)\wX^{(K-3)}_{0,t}+\mathcal{R}^F_{K-2}(t,\widehat \XX).
\end{align} Since $\nabla^{K-3}F$ is a $C^3$ non-anticipative Marcus canonical  path functional (see Remark~\ref{remark: CK functionals}), and $\nabla ^{K-3}F,\nabla ^{K-2}F,\nabla ^{K-1}F $, and $R^{\nabla^{K-2} F,\nabla^{K-1} F}$ are $[p,3)$\text{-}var continuous, a component-wise application of Theorem~\ref{3th: Ito formula p 23} to $ \nabla^{K-3}F$ yields that for all $s\in [0,1]$,
    \begin{align}\label{3eq: Taylor p23 linear from  k-3 to k-2}
         \nabla^{K-3}F(s,\widehat \XX)=  \nabla^{K-3} F(0,\widehat \XX)+\int_0^s \nabla^{K-2} F(r,\widehat \XX)d\widehat \XX_r.
    \end{align}
Finally, since any rough integral with respect to $\wXX$ is nothing else than a Young integral with respect to $\pi_1(\wXX)$, (see e.g., Lemma~\ref{lemma: consistency}), the claim follows by replacing the expression in equation~\eqref{3eq: Taylor p23 linear from  k-3 to k-2} into $\mathcal{R}^F_{K-3}(t,\widehat \XX)$.

\paragraph{Step 2:}  Fix $\widehat \XX \in  \widehat C^{p}([0,1],G^{2}(\Rset^{d+1})) $. By Theorem 5.23 in~\cite{FV:10}, there exists a sequence of piecewise linear time-extended paths such that their truncated signature at level $2$, denoted by $(\widehat \XX^M)_{M\in \Nset}$, satisfies
    \begin{align*}
        \lim_{M\rightarrow\infty }d_{CC}\sup_{t\in[0,1]}(\widehat \XX^M_t,\widehat \XX_t)=0, \qquad \text{and}\qquad    \sup_{M\in \Nset}\|\widehat \XX^M\|_{p\text{-}var}\leq C\|\|\widehat \XX\|_{p\text{-}var},
    \end{align*}
    for some $C>0$. 
    By Step 1, for every fixed $M$ and $t\in [0,1]$, 
       \begin{align*}
    F(t,\widehat \XX^M)=&\sum_{j=0}^{K-3}\nabla^{j}F (0,\widehat \XX^M)\widehat \X^{M,(j)}_{0,t}+ \mathcal{R}^F_{K-2}(t,\widehat \XX^{M}).
\end{align*}
Thus, since for all $j=0,\dots,K-3$, $ \nabla ^{j}F $ and  $\wX^{(j)}$ are $[p,3)$\text{-}var continuous (see Proposition~\ref{prop: derivative linear func}), the claim follows if we prove that for every $t\in [0,1]$,
\begin{align*}
    \lim_{M\rightarrow\infty} \mathcal{R}^F_{K-2}(t,\wXX^{M})= \mathcal{R}^F_{K-2}(t,\wXX).
\end{align*}
To this end,  set $G^{-1}:=\nabla ^{K-1}F$, $G^0:=\nabla ^{K-2}F$ and define for $j=1,\dots,K-2$, 
\begin{align*}
G^j(t,\wXX):=\int_0^t\int_0^{t_{1}}\dots \int_0^{t_{j-1}} \nabla ^{K-2}F(t_{j},\wXX)d\wXX_{t_{j}}\dots d\wXX_{t_{1}}.
\end{align*}
We show that there exists $p:=p^{(0)}<p^{(1)}<\dots p^{(K-2)}<3$ such that for all $j=1,\dots,K-2$, 
\begin{align}\label{3eq: p^j taylor convergence}
&\lim_{M\rightarrow\infty }d_{p^{(j)}}(G^{j-2}(\cdot,\widehat \XX^M), G^{j-2} (\cdot,\widehat \XX))=0,\\
        & \lim_{M\rightarrow\infty }d_{p^{(j)}/2}\big(R^{G^{j-1},G^{j-2}}((\cdot,\cdot),\wXX^M),R^{G^{j-1},G^{j-2}}((\cdot,\cdot),\wXX)\big)=0,\nonumber
\end{align}
and $\lim_{M\rightarrow\infty}d_{CC}(G^{j-1}(0,\wXX^M),G^{j-1}(0,\wXX^M))=0.$
Then, a similar argument as in the Step 2 of the proof of Theorem~\ref{3th: Ito formula p 23} (see equation~\eqref{3eq: proof  step 2 ito p23 }) yields that  for all $j=1,\dots,K-2$, 
\begin{align}\label{3eq: convergence dpj final}
\lim_{M\rightarrow\infty}d_{p^{(j)}\text{-}var}\big( G^j(\cdot,\wXX^M), G^j(\cdot ,\wXX))=0.
\end{align}
Since $G^{K-2}=\mathcal{R}_{K-2}^F$, the claim follows. 

We reason by induction. Let $j=1$ and fix $p^{(1)}>p$.  By interpolation, $$\lim_{M\rightarrow\infty}d_{p^{(1)}\text{-}var}\big(\wXX^M,\wXX)=0.$$ Since by assumption $G^{-1}$, $G^0$, and $R^{G^{0},G^{-1}}$ are $[p,3)$-var continuous, the claim follows.

Assume that the assertion holds  for all $l=1,\dots,j$, with $K-3\geq j>1$, and fix $p^{(j+1)}>p^{(j)}>p^{(j-1)}$. By the induction hypothesis, 
\begin{align}\label{3eq: p3 taylor induction}
\lim_{M\rightarrow\infty}d_{p^{(j)}\text{-}var}\big( G^{j-1}(\cdot,\wXX^M), G^{j-1}(\cdot ,\wXX)\big)=0.
\end{align}
Moreover, an application of equation (2.3.3) in Theorem 31 in~\cite{FS:17} yields that 
$$  \sup_{M\in \Nset }  \|R^{G^{j},G^{j-1}} ((\cdot,\cdot),\wXX^M)\|_{p^{(j)}/2\text{-}var}<\infty.$$
Therefore, by interpolation
\begin{align*}
    \lim_{M\rightarrow\infty }d_{p^{(j+1)}/2}\big(R^{G^{j},G^{j-1}}((\cdot,\cdot),\wXX^M),R^{G^{j},G^{j-1}}((\cdot,\cdot),\wXX)\big)=0.
\end{align*}
Since the convergence in 
\eqref{3eq: p3 taylor induction} holds  also with respect to $d_{p^{(j+1)}}$, the claim follows.

\paragraph{Step 3:} It follows as the Step 2 in the proof of Theorem \ref{th: taylor 12}.
\qed

\section{Auxiliary remarks}\label{3appendix: integration theories}

\subsection{Young and (level 2) rough integration}\label{sec: rough=young}
In this section, we review the notion of the Young and the (level 2) rough integral and discuss their consistency. First, we review relevant notions of convergence.

\begin{definition}[\cite{FZ:17}] \label{3def: MRS RRS}For $\Xi:\Delta_1\rightarrow\Rset^m$ and  a partition $\pi_{[0,1]}$  of the interval $[0,1]$, we say that the quantity $\sum_{s_i\in \pi_{[0,1]}}\Xi_{s_i,s_{i-1}}$ converges to $K\in \Rset^d$, 
\begin{enumerate}
    \item in the Mesh Riemann-Stieltjes
(MRS) sense: if for any $\varepsilon>0$, there exists $\delta>0$ such that for all $\pi_{[0,1]}$ with $|\pi_{[0,1]}|<\delta$,  $|\sum_{s_i\in \pi_{[0,1]}}\Xi_{s_i,s_{i-1}}-K|\leq \varepsilon$.
\item in the Refinement Riemann-Stieltjes (RRS) sense: 
if for any $\varepsilon>0$ there exists a partition $\pi^\varepsilon_{[0,1]}$ such that for any refinement $\pi_{[0,1]}$ of $\pi^\varepsilon_{[0,1]}$,  $|\sum_{s_i\in \pi_{[0,1]}}\Xi_{s_i,s_{i-1}}-K|\leq \varepsilon$.
\end{enumerate}
\end{definition}

\noindent  For $\XX\in D^{p}([0,1],G^{[p]}(\Rset^{d}))$, $(s,t)\in \Delta_1$, set   $X_{s,t}:=:\pi_1(\XX_{s,t})$ and $\X^{(2)}_{s,t}:=\pi_2(\XX_{s,t})$ if $p\in[2,3)$.

\begin{proposition}\label{prop:young}(Young integration, \textit{Proposition} 2.4 in~\cite{FZ:17}).
       Fix $p\in [1,2)$ and let $\XX\in D^{p}([0,1],G^{1}(\Rset^{d}))$ and $Y\in D^{p'}([0,1],\mathcal{L}(\Rset^d,\Rset^m))$, for some $p'$ such that $\frac{1}{p}+\frac{1}{p'}>1$. Then, for all $t\in [0,1]$, the limit
     \begin{align}\label{3eq: Young integral}
         \int_0^tY_{s^-}d\XX_s:=\lim_{(MRS)|\pi_{[0,t]}|\rightarrow0}\sum_{s_i\in \pi_{[0,t]}} Y_{s_i}X_{s_i,s_{i+1}}
     \end{align}
     exists. Moreover,
    $\int_0^\cdot Y_{s^-}d\XX_s\in D^{p}([0,1],\Rset^{m})$. We call the limit in~\eqref{3eq: Young integral} the \textit{Young integral} of $Y$ with respect to $\XX$.

\end{proposition}

Recall from Definition \ref{3def: controlled path qr} that for $p\in [2,3)$, $\XX\in D^{p}([0,1],G^{2}(\Rset^{d}))$, $p'\geq p$ such that $  \frac{2}{p}+\frac{1}{p'}>1$,  and $r\geq 1$ given by $ \frac{1}{r}=\frac{1}{p}+\frac{1}{p'},$  $\mathcal{V}^{p',r}_{\XX}$ denotes the space of controlled rough paths with respect to $X=\pi_1(\XX)$.

\begin{proposition}\label{3prop: rough integral}(Level 2 rough integration, \textit{Proposition} 2.4 in~\cite{ALP:23})
     Let $\XX\in D^{p}([0,1],G^{2}(\Rset^{d}))$ and $(Y,Y')\in \mathcal{V}^{p',r}_{\XX}$, for some $p'\geq p$ such that $\frac{2}{p}+\frac{1}{p'}>1$ and $r\geq 1$ given by $ \frac{1}{r}=\frac{1}{p}+\frac{1}{p'}$ . Then, for all $t\in [0,1]$, the limit
     \begin{align}\label{3eq: rough integral}
         \int_0^tY_{s^-}d\XX_s:=\lim_{(MRS)|\pi_{[0,t]}|\rightarrow0}\sum_{s_i\in \pi_{[0,t]}}\big( Y_{s_i}X_{s_i,s_{i+1}}+Y'_{s_i}\X_{s_i,s_{i+1}}^{(2)}\big),
     \end{align}
     exists. Moreover,  $(\int_0^\cdot Y_{s^-}d\wXX_s,Y)\in \mathcal{V}^{p',r}_{\wXX}$.  We call the limit in~\eqref{3eq: rough integral} the \textit{rough integral} of $(Y,Y')$ with respect to $\XX$.
\end{proposition}

\begin{notation}
    Throughout the paper, whenever the involved paths are continuous, we write $\int_0^\cdot Y_sd\XX_s$ in~\eqref{3eq: Young integral} or~\eqref{3eq: rough integral}.
\end{notation}

Fix $X\in D^1([0,1],G^1(\Rset^d))$, and for all $t\in [0,1]$ let
\begin{align}\label{eq: lift young integral}
  \XX_t:=\left (1,X_t-X_0,\int_0^tX_{0,s^-}\otimes dX_{s}+\frac{1}{2}\sum_{0< s \leq t}(\Delta X_s)^{\otimes 2}\right ),
\end{align}
where  $\int_0^\cdot X_{0,s^-}\otimes dX_s$ denotes the 
Young integral of $X_{0,\cdot}$ with respect to $X$. By Proposition \ref{prop:young} and the geometric properties of the Young integral (e.g., Proposition 2.4 in~\cite{FZ:17}), $\XX\in D^1([0,1],G^2(\Rset^d))$. Since for any $p>1$, $D^1([0,1],G^2(\Rset^d))\subseteq D^p([0,1],G^2(\Rset^d))$, $\XX$ can be interpreted as a weakly geometric \cadlag $p$-rough path, for any $p\in [2,3)$.

\begin{lemma}\label{lemma: consistency}
	Let  
 $\XX$ be defined as in equation~\eqref{eq: lift young integral}, for some $X\in D^1([0,1],G^1(\Rset^d))$, and let $ (Y,Y')\in \mathcal{V}^{p',r}_{\wXX},$  with $p'>p$, for some $p\in[2,3)$ such that $\frac{1}{p'}+\frac{2}{p}>1$. Then, 
 \begin{align}\label{eq: Rough = Marcus}
     \int_0^tY_{s^-}d\mathbf{X}_s=\int_0^tY_{s^-}d{X}_s+\frac{1}{2}\sum_{0<s\leq t}Y_{s^-}'(\Delta X_s)^{\otimes 2},
 \end{align}
 where the integral on the LHS denotes the rough integral of $(Y,Y')$ with respect to $\XX$, while the integrals on the RHS are Young integrals of $Y$ with respect to $X$, and $Y'$ with respect to the path $[0,1]\ni t\mapsto \sum_{0<s\leq t}(\Delta X_s)^{\otimes 2}$, respectively.
\end{lemma}
\begin{proof}
Notice that the path $[0,1]\ni t\mapsto \sum_{0<s\leq t}(\Delta X_s)^{\otimes 2}$ is of finite variation, implying that the second integral on the RHS of equation~\eqref{eq: Rough = Marcus} is well defined. Next, set for all $(s,t)\in \Delta_1$, $ \Y_{s,t}:=\int_s^tX_{s,r^-}\otimes dX_r$
and fix $u\in (0,1]$. By the definitions of rough and Young integral (see equations~\eqref{3eq: rough integral},~\eqref{3eq: Young integral}, respectively),
we need to show that
\[
\lim_{(MRS) |\pi_{[0,u]}|\rightarrow0}\sum_{s_{i}\in \pi_{[0,u]}}Y'_{s_i}\Y_{s_{i},s_{i+1}}=0.
\]
The claim follows by a similar reasoning as in  Theorem~35 of ~\cite{FS:17}.
\end{proof}

\subsection{Paths of finite quadratic variation in the sense of Föllmer}
In this section, we briefly revisit the concept of path of finite quadratic variation in the sense of~\cite{F:81}.   

\begin{definition}\label{3def: quadratic variation}
   Let $\mathcal{B}([0,1])$ denote the Borel $\sigma$-algebra on $[0,1]$. Let $X\in D([0,1],\Rset)$ and $\pi^n_{[0,1]}$, $n\in \N$, a sequence of partitions with vanishing mesh size. We say that $X$ has finite \textit{quadratic variation along} $(\pi^n_{[0,1]})_{n\in\Nset} $ \textit{in the sense of \Follmer} if the sequence of measures $(\nu_n)_{n\in\N}$ on $([0,1],\mathcal{B}([0,1])$ defined by
    \begin{align*}
        \nu_n:=\sum_{s^n_i\in \pi^n_{[0,1]}}\|X_{s^n_i,s^n_{i+1}}\|\delta_{s^n_i}
    \end{align*}
    converges weakly to a measure $\nu$ such that the map
    \begin{align}\label{3eq: jumps quadratic variation}
       [0,1]\ni  t\mapsto [X]_t^c:=\nu([0,t])-\sum_{0<s\leq t}\|\Delta X_s\|^2
    \end{align}
    is continuous and increasing.  We then call the function $[X]^c$ and the function $[X]$ given by $[X]_t:=\nu([0,t])$ the \textit{continuous quadratic variation} and \textit{quadratic variation} of $X$ along $(\pi^n_{[0,1]})_{n\in\Nset}$, respectively. We say that a path $X\in D([0,1],\Rset^d)$, $d\geq1$, has finite quadratic variation along $(\pi^n_{[0,1]})_{n\in\Nset}$ if for all $i,j\in \{1,\dots,d\}$, the condition above holds for $X^i$ and $X^i+X^j$. In this case, we set
    \begin{align*}
        [X^i,X^j]:=\frac{1}{2}([X^i+X^j]-[X^i]-[X^j]),
    \end{align*}
 and similarly for $[X^i,X^j]^c$.
\end{definition}
\begin{remark}\label{3rem: quadratic variation}
  Notice that for all $i,j\in \{1,\dots,d\}$, $[X^i,X^j]=[X^j,X^i]$ and the map $[0,1]\ni t\mapsto [X,X]_t\in (\Rset^d)^{\otimes 2}$ is a \cadlag path of finite $1$-variation. Thus, for any $Z\in D([0,1],\mathcal{L}((\Rset^d)^{\otimes 2},\Rset^m))$, $m\in \Nset$, the Young integral of $Z$ with respect to $[X,X]$  is well defined. 
  \end{remark}
 A similar argument as in the proof of Théorème in~\cite{F:81} yields the following result. 

  \begin{proposition}
      Let $Z\in D([0,1],\mathcal{L}((\Rset^d)^{\otimes 2},\Rset^m))$, $m\in \Nset$ and $X\in D([0,1],\Rset)$ be a path of finite
quadratic variation along some sequence of partition $(\pi^n_{[0,1]})_{n\in\Nset} $ in the sense of \Follmer. Assume that 
    \begin{align*}
         \{s\in (0,1] \ \colon \ \Delta Z_s\neq 0\}\subseteq  \{s\in (0,1] \ \colon \ \Delta X_s\neq 0\}.
    \end{align*}
    Then, for all $t\in (0,1]$, 
    \begin{align}\label{3eq: integral qv=Young}
    \lim_{n\rightarrow\infty}\sum_{s^n_i\in \pi^n_{[0,1]} }Z_{s^n_i}X^{\otimes 2}_{s^n_i\wedge t ,s^n_{i+1}\wedge t}=\int_0^t Z_{s^-} d[X,X]_s,
    \end{align}
    where the intergral on the right-hand side is a Young integral.
    \end{proposition}

\begin{remark}
    In view of~\eqref{3eq: jumps quadratic variation},
\begin{align*}
    \int_0^t Z_{s^-} d[X,X]_s=\int_0^t Z_{s^-} d[X,X]^c_s+\sum_{0<s\leq t} Z_{s^-} \Delta X_s^{\otimes 2}.
\end{align*}
We refer also to Lemma 5.11 in~\cite{FH:14} for a proof of the equality~\eqref{3eq: integral qv=Young} when $X$ is continuous, and Lemma 2.6 in~\cite{H:17} when $X$ is \cadlag and $Z=f(X)$, for some suitably chosen $\R^m$-valued function $f$.
\end{remark}

\subsection{Paths satisfying (RIE) property}
Another line of research that lies between rough and Föllmer integration theory  exploits the so-called RIE property of a path with a given sequence of partition (see \cite{PP:16} and~\cite{ALP:23}). 

\begin{definition}
    A sequence of partitions $(\pi^n_{[0,1]})_{n\in\Nset}$  is called
nested if $\pi^n_{[0,1]}\subseteq \pi^{n+1}_{[0,1]}$ for all $n\in\Nset$.
\end{definition}

\begin{property}[RIE]\label{3property RIE}
Let $p\in (2,3)$ and $(\pi^n_{[0,1]})_{n\in\Nset}$ a sequence of nested partitions with vanishing mesh size. For $X\in D([0,1],\Rset^d)$, define $X^n:[0,1]\rightarrow\Rset^d$ by 
\begin{align*}
    X^n_u:=X_1 1_{\{u=1\}}+\sum_{s^n_i\in \pi^n_{[0,1]}}X_{s^n_i}1_{\{s^n_{i}\leq u<s^n_{i+1}\}},
\end{align*}
for all $u\in [0,1]$, $n\in\Nset$. Assume that 
\begin{enumerate}
    \item The sequence $(X^n)_{n\in\N}$ converges uniformly to $X$ as $n\rightarrow\infty$.
    \item The Riemann sums
    \begin{align}\label{riem sums rie}
        \int_0^tX^n_{u^-}dX_u:=\sum_{s^n_i\in \pi^n_{[0,1]}}X_{s^n_i}X_{s^n_i\wedge t, s^n_{i+1}\wedge t}
    \end{align}
    converges uniformly as $n\rightarrow\infty$ to a limit denoted by $\int_0^tX_{u^-} dX_u$, for all $t\in [0,1]$. 
    \item There exists a control function $w$  such that
    \begin{align}\label{3eq: def RIE iii}
        \sup_{(s,t)\in \Delta_1}\frac{\|X_{s,t}\|^p}{w(s,t)}+\sup_{n\in \N}\sup_{t^n_i<t^n_j\in \pi^n_{[0,1]}}\frac{\|\int_{t^n_i}^{t^n_j}X^n_{u^-}dX_u-X_{t^n_i}X_{t^n_i,t^n_j}\|^p}{w(t^n_i,t^n_j)}\leq 1,
    \end{align}
    with the convention of $\frac{0}{0}:=0$.
\end{enumerate}
\end{property}

\begin{definition}[\cite{ALP:23}]
      A path $X\in D([0,1],\R^d)$ is said to \textit{satisfy} (RIE) with respect to $p\in (2,3)$ and $(\pi^n_{[0,1]})_{n\in\Nset}$ if $p$, $(\pi^n_{[0,1]})_{n\in\Nset}$,  and $X$ together satisfy property (RIE).
\end{definition}

By Proposition 2.18 in~\cite{ALP:23} every path $X\in D([0,1],\R^d)$ that satisfies (RIE) property with respect to some $p\in (2,3)$ and some sequence of nested partition $(\pi^n_{[0,1]})_{n\in\Nset}$  has finite quadratic variation $[X,X]$ along $(\pi^n_{[0,1]})_{n\in\Nset}$ in the sense of Föllmer.

 \begin{proposition}\label{3prop: RIE integal properties}
      Let $p\in (2,3)$, $(\pi^n_{[0,1]})_{n\in\Nset}$ a sequence of nested partition with vanishing mesh size and $X\in D([0,1],\R^d)$  with $X_0=0$ satisfying (RIE) with respect to $p$ and $(\pi^n_{[0,1]})_{n\in\Nset}$. Then, the following conditions are satisfied.
      \begin{enumerate}
          \item \label{i RIE}\label{item RIE 2} For $t\in [0,1]$, take  $\int_0^t X_{u^-}dX_s$ the limit of the Riemann sums $(\int_0^tX^n_{s^-}dX_s)_{n\in \Nset}$  in Property (RIE)~\ref{3property RIE} (ii). For $(s,t)\in \Delta_1$, set $\int_s^t X_{u^-}dX_u:=\int_0^tX_{u^-}dX_u-\int_0^sX_{u^-}dX_u$, and define
    \begin{align*}
        \X^{(2)}_{s,t}:=\int_s^t X_{u^-}dX_u-X_sX_{s,t}+\frac{1}{2}[X,X]_{s,t},    \end{align*}
        and $\XX_t:=(1,X_t,  \X^{(2)}_{0,t})$. Then, $\XX\in D^p([0,1],G^2(\Rset^{d}))$.
   
   \item\label{item RIE 3}  Let $(Y,Y')\in \mathcal{V}^{q,r}_{\XX}$, for some $q\geq p$. Assume that 
          \begin{align*}
              &\{s\in (0,1] \ \colon \ \Delta Y_s\neq 0\}\subseteq \bigcup_{n\in \N}\pi^n_{[0,1]}, \quad \{s\in (0,1] \ \colon \ \Delta Y'_s\neq 0\}\subseteq  \{s\in (0,1] \ \colon \ \Delta X_s\neq 0\}.
          \end{align*}
          Then, for all $t\in [0,1]$, 
            \begin{align}\label{3eq: limit rough integral RIE}
              \int_0^t Y_{s^-}d\XX_s= \int_0^t Y_{s^-}dX_s+ \frac{1}{2}\int_0^t Y'_{s^-}d[X,X]_s,
          \end{align}
          where $  \int_0^t Y_{s^-} dX_s:= \lim_{n\rightarrow\infty }\sum_{s^n_i\in \pi^n_{[0,1]}}Y_{s^n_i} X_{s^n_i\wedge t, s^n_{i+1}\wedge t},$
is a well defined limit along the sequence of partition $(\pi^n_{[0,1]})_{n\in\Nset}$.
      \end{enumerate}
\end{proposition}

\begin{proof}

~\ref{item RIE 2}:  By Proposition 2.18 in~\cite{ALP:23} (in 
 particular equation (2.14)),  it holds that $\XX\in D([0,1],G^2(\Rset^{d}))$. Moreover, $\XX$ has finite $p$-variation by Lemma 2.13 in~\cite{ALP:23} and  Remark~\ref{3rem: quadratic variation}. 

~\ref{item RIE 3}: It follows by the definition of the rough integral with respect to $\XX$ (see equation~\eqref{3eq: rough integral}), Theorem 2.15 in~\cite{ALP:23}, and Remark~\ref{3rem: quadratic variation}.

\end{proof}

\end{document}